\documentclass[12pt]{amsart}
\usepackage{amssymb}
\usepackage{graphicx}

\input klein.def
\usepackage{hyperref}

\usepackage[all]{xy}

\numberwithin{equation}{section}

\begin{document}

\myitemmargin
\baselineskip =15pt plus 3pt

\title[Lagrangian Klein bottle and mapping class groups]%
{Lagrangian embeddings of the Klein bottle\\[4pt]
and\\[4pt] combinatorial properties of mapping class groups}
\author[V.~Shevchishin]{Vsevolod V.~Shevchishin}
\address{Mathematisches Institut der Universit\"at Bonn\\
Beringstra\ss e 1\\
 D-53115 Bonn \\
Germany} \email{sewa@math.uni-bonn.de}
\dedicatory{} \subjclass{} \thanks{This research was carried out with the
 financial support of the Heisenberg program of the Deutsche
 Forschungsgemeinschaft}
\begin{abstract}
 In this paper we prove the non-existence of Lagrangian embeddings of the
 Klein bottle $K$ in $\cp{}^2$. We exploit the existence of a special
 embedding of $K$ in a symplectic Lefschetz pencil $\pr: X \to S^2$ and study
 its monodromy. As the main technical tool, we develop the combinatorial
 theory of mapping class groups. The obtained results allow us to show that in
 the case when the class $[K] \in \sfh_2(X,\zz_2)$ is trivial the monodromy of
 $\pr: X \to S^2$ must have special form. Finally, we show that such a monodromy
 can not be realised in the case $\cp{}^2$.

\end{abstract}
\maketitle
\setcounter{tocdepth}{2}

\vskip-30pt

\setcounter{section}{-1} \pagebreak[1]

\newsection[intro]{Introduction}

In \cite{N-1}, Stefan Nemirovski proposed a proof of the 
non-existence of a Lagrangian embedding of the Klein bottle $K :=
\rr\pp^2 \# \rr\pp^2$ into $\cp{}^2$.  His main claim, that for any
Lagrangian embedding of the Klein bottle in a complex algebraic
surface the $\zz_2$-homology class $[K]$ is non-trivial, is false:
We construct a counterexample in \refsubsection{lagr-triv}. Another
attempt to prove this result was made by Klaus Mohnke in \cite{Mo}.
He used completely different techniques, but his proof also appears
to be incomplete.

The main result of the present paper is the following:

\newthma{thm1}{\slsf{Main theorem}} Any Lagrangian embedded Klein bottle $K$
in $\cp^2$ or in a compact ruled symplectic 4-manifold $(X, \omega)$ has non-trivial
$\zz_2$-homology class $[K] \in \sfh_2(X, \zz_2)$.

\smallskip%
\rm By \cite{MD-Sa}, a ruled symplectic manifold $X$ is (a blow-up of) a compact
\emph{complex} ruled surface, and $\omega$ is a K{\"a}hler form on $X$. The
characterisation of such $X$ from \cite{MD-Sa} will be used in the proof.\\[-20pt]
\end{thma}

According to \cite{AMP}, every Lagrangian embedding of $K$ in a symplectic $4$-manifold~$X$
can be modified into a \slsf{Morse-Lefschetz fibration} $\pr: (\wh X, K) \to (Y, \Gamma)$,
which is a Lefschetz fibration of a special form, see \refsubsection{sy-le-fi}.
The technical core of our proof is a study of the monodromy of such
Morse-Lefschetz fibrations. For this purpose, we obtain several new results
on the combinatorial structure of the mapping class group $\map_g$, see \refsubsection{braid-map}.

\smallskip%
As in \cite{N-1}, we can deduce from \refthm{thm1} the following result:

\newthma{thm2}{} There is no Lagrangian embedding of the Klein bottle $K$ in
$\cp^2$.
\end{thma}

Indeed, the $\zz_2$-self-intersection number of an embedded Lagrangian Klein bottle
in any $(X,\omega)$ equals its $\zz_2$-Euler class and hence vanishes, whereas $h^2 \neq0$
for the \emph{unique} non-trivial $h\in \sfh_2(\cp^2, \zz_2)$. Note that
Lagrangian embeddings of the Klein bottle in blown-up $\cp^2$ do exist, as an
easy example from \cite{N-1} shows: One realises $K$ as the \emph{real
 algebraic\/} blow-up of the real projective plane $\rr\pp^2$ and complexifies
this construction obtaining $\cp^2$ blown-up at one point.

\medskip
The following application was pointed out by St.~Nemirovski, see also 
\cite{N-1}, \cite{Aud}, and  \cite{Ga-La}, p.~489.

\newthma{thm4}{} Let $X$ be $\cc\pp^2$ or a compact ruled symplectic
$4$-manifold. Then two Lagrangian embeddings $\varphi_1, \varphi_2: \rr\pp^2 \to X$ with
transversal images representing the same $\zz_2$-homology class have at least
$3$ intersection points.
\end{thma}

Indeed, the $\zz_2$-self-intersection index of any Lagrangian embedding $\varphi:
\rr\pp^2 \to X$ 
coincides with the Euler characteristic of $\rr\pp^2$. In
particular, the class $[\varphi(\rr\pp^2)]$ is non-trivial and two transversal 
Lagrangian embeddings must have an odd number of intersection points. 
In the case of a single intersection point, one could perform
Lagrangian surgery on it producing a Lagrangian embedding of the Klein bottle
$K$. Its $\zz_2$-homology class is the sum of the classes of the
images, hence trivial which contradicts \slsf{Main Theorem}. Notice also that 
in $\cp{}^2$ the condition on the homology classes is fulfilled automatically
since there is a unique non-trivial element in $\sfh_2(\cp^2, \zz_2)$.

\smallskip%
Audin \cite{Aud} has established an additional restriction on the
class $[K]$, namely, $[K]^2 \equiv0 \mod4$. Here $[K]^2$ is
calculated as the Pontryagin square (of the Poincar{\'e} dual) of
$[K]$. For example, in the case when $X = S^2 \times S^2$, the class
$[K]$ can not be the sum of the vertical and the horizontal fibres,
since the Pontryagin square of such a sum is $2\mod4$. On the other
hand, it is not difficult to show that for any symplectic form of
product type $\omega = p_1^*\omega_1 + p_2^*\omega_2$ there exists a
Lagrangian embedding of the Klein bottle $K$ such that the
projection $p_1: S^2 \times S^2 \to S^2$ on the first factor,
restricted to $K$, is an $S^1$-bundle over a circle $\Gamma$ in the
first $S^2$. As we show in \propo{prop-mu-0} in this case $K$ is
homologous to the fibre of the projection $p_1$.

\newsubsection[skim]{Scheme of the proof} Let us recall Nemirovski's original
ideas.  He showed in \cite{N-1} that, given a Lagrangian submanifold $L$ in a
projective manifold $X$ and an $S^1$-valued Morse function $f: L\to S^1$, there
exist another submanifold $L'$ isotopic to $L$, a blow-up $\wh X \to X$ of $X$
along a complex submanifold $B \subset X$ with $B \bigcap L' = \varnothing$,
and a holomorphic Lefschetz fibration $\pr: \wh X \to\cp^1$ such that $(\wh X, \pr, L')$
form a Morse-Lefschetz fibration (see~\refdefi{df1}). Auroux, Mu{\~n}os, and Presas \cite{AMP}
generalised this result to the case of an arbitrary ambient compact symplectic manifold $X$,
see \refthm{thmMLF}. For our purposes we need a slight refinement of this generalisation;
it is established in~\lemma{lem-slf}.

In the case $\dim_\rr X =4$ the homology group $\sfh_2(X, \zz_2)$ is
naturally embedded in the group $\sfh_2(\wh X, \zz_2)$, and thus the
blow-up procedure does not affect the (non-)triviality of the
homology class of a surface avoiding the blow-up locus. Thus we may
assume that the Lefschetz fibration $\pr: X \to S^2$ is defined
already on $X$. As the Morse function $f$ on $K$, Nemirovski
chose a (topologically unique) non-singular fibration $\pr_K : K\to S^1$ with circle fibres.
As a result, he obtained a rather simple geometric configuration, described in \lemma{lem-slf} below.

The key idea of Nemirovski was to exploit the arising structure of a
Morse-Lefschetz fibration and construct a complex line bundle $\scrl$ on $X$ such that
$\lan c_1(\scrl), [K] \ran \equiv 1 \mod 2$, which would have implied the non-triviality
of $[K]$. Instead, in our approach we use this structure to describe the
necessary condition for the vanishing of $[K]$ in terms of the
Lefschetz fibration and its monodromy. First of all, we find an expression
for $[K]$ using the homology spectral sequence. This allows us to produce
an example of a Lagrangian embedding of the Klein bottle in a projective surface
with \emph{trivial} homology class $[K]$, see \refsubsection{lagr-triv}.
Next, we undertake a detailed study of the combinatorial structure of the mapping
class group $\map_g$. This technique (see Subsec.~\ref{braid-map}) allows us
to show that if the class $[K]$ is trivial, then the monodromy of our Lefschetz fibration
must have very special form, see \propo{prop-fine}. Then we observe that ``twisting''
appropriately a part of this special monodromy, we can construct a new symplectic
manifold $X'$ such that
\begin{equation}\eqqno(rank+1)
\rank\, \sfh_1(X', \zz_2)  = \rank\, \sfh_1(X, \zz_2) +1.
\end{equation}
On the other hand, using the classification of ruled symplectic manifolds from~\cite{MD-Sa}
we show that under the assumptions of the \slsf{Main Theorem} the manifold $X'$ must also be ruled,
and hence the ranks of both groups $\sfh_1(X, \zz_2)$ and $\sfh_1(X', \zz_2)$ must be even.
The obtained contradiction finishes the proof.

\newsubsection[0.1a]{Luttinger surgery and
totally real embeddings of the Klein bottle}
After the preprint version of this paper was finished and distributed in May
2006, Nemirovski has found an alternative proof of \refthm{thm1}, see
\cite{N-2}. He observes that the ``twisting'' of the manifold $X$
above is the \slsf{Luttinger surgery} along the Klein bottle $K$
(well-known in the case of (Lagrangian) embeddings of the torus $T^2$, see
\cite{Lutt,E-Po,ADK}). He shows also that equality \eqqref(rank+1)
can be deduced using the \slsf{Viro index} of curves on the Klein bottle
and holds for any totally real embedding of $K$ in an
\emph{almost complex} four-manifold~$(X,J)$.

\newsubsection[braid-map]{Braid combinatorial structure of mapping class
 groups} As we have mentioned above, for the proof of the \slsf{Main Theorem} we
need effective tools for calculations in the mapping class group $\map_g$ of a
surface $\Sigma$ of given genus $g\geq1$. This requires a better understanding of the
combinatorial structure of $\map_g$.

It was understood by Dehn \cite{De} that $\map_g$ is generated by certain simple
transformations called \slsf{Dehn twists}.  Seeking for a simple presentation
for $\map_g$ one first sees that the ``basic'' relation among Dehn twists
is the braid relation. So it was expected that $\map_g$ should resemble a braid group,
possibly generalised. And indeed, Wajnryb's presentation \cite{Waj} realises $\map_g$
as a quotient of a certain Artin--Brieskorn braid group $\br(\cals_g)$ with
a simple system of generators $\cals_g$ and relations.  Recently Matsumoto
\cite{Ma} gave a description of relations in Wajnryb's presentation in
terms of the so-called \slsf{Garside elements} $\idel(\cals)$ corresponding to
certain subsystems $\cals\subset\cals_g$. In \refsection{map-mu} we obtain some
combinatorial properties of the Wajnryb--Matsumoto presentation, which are used
in a crucial way in the proof of the \slsf{Main Theorem}.

Recall that every Artin--Brieskorn braid group $\br(\cals)$ is constructed from
a \slsf{Coxeter system} $(\cals, M_\cals)$ such that $\cals$ is the set of
generators of $\br(\cals)$ and $M_\cals$ is a matrix on $\cals$ encoding
the combinatorics of the defining set of relations. To every such system
$(\cals, M_\cals)$ one can also associate the so-called \slsf{Coxeter--Weyl
group} $\sfw(\cals)$ by adding the relation $s^2=1$ for each generator $s\in \cals$.
In particular, there exists a natural surjective homomorphism
$\br(\cals) \to \sfw(\cals)$ which is a bijection on generators. The
kernel of this homomorphism is called the \slsf{pure braid group}
associated with $\cals$ and will be denoted by $\sfp(\cals)$. Recall that
\slsf{quasireflections} in $\sfw(\cals)$ and \slsf{quasigenerators} in $\br(\cals)$
are the elements conjugate to the standard generators $s\in \cals$.
We denote by $\calt(\cals)$ the set of quasireflections in~$\sfw(\cals)$.

Our first result about these groups is the following:

\newthma{thm5}{\refthm{thm-gen-pure}} There exist natural lifts
$\sfw(\cals)\ni w \mapsto \hat w \in \br(\cals)$ and $\calt(\cals)\ni t \mapsto \ti t \in \br(\cals)$
such that the set $\wt\calt^2(\cals) := \{ \ti t^2 : t\in \calt(\cals)\}$
is a minimal generator set for the  pure braid group $\sfp(\cals)$.

Moreover, the abelianisation $\sfp_\ab(\cals)$ of the pure braid group $\sfp(\cals)$
is a free abelian group with the basis $\{ A_t : t\in \calt(\cals)\}$ indexed by the
set $\calt(\cals)$ of quasireflections such that $A_t$ is the projection
of $\ti t^2$ to $\sfp_\ab(\cals)$.
\end{thma}

\noindent
{\bf Remark.} This result is well-known for the usual braid group $\br_d$.

\smallskip
If the Coxeter--Weyl group $\sfw(\cals)$ is finite, then it can be realised as
a \slsf{reflection group}. The corresponding Coxeter systems $\cals$ are completely
classified. For such systems Deligne \cite{Del} and Brieskorn-Saito
\cite{Br-Sa} defined the so-called \slsf{Garside element} $\idel(\cals) \in
\br(\cals)$. Its square $\idel^2(\cals)$ lies in the centre of $\sfp(\cals)$,
and we show in \lemma{lem-gars-sq}\, that the projection of $\idel^2(\cals)$ to
$\sfp_\ab(\cals)$ is simply the sum $\sum_{t\in \calt(\cals)} A_t$ of all basis
elements of $\sfp_\ab(\cals)$. This immediately gives a description of the
projections of Matsumoto's relation elements to the group $\sfp_\ab(\cals_g)$.

\medskip
Our next important result is:

\newthma{thm6}{\refthm{thm-hurw-W}} Every factorisation $\mbft = t_1 \cdot t_2
\cdot \ldots \cdot t_l$ of the identity $1\in \sfw(\cals)$ into a product of
quasireflections $t_i\in \calt(\cals)$ is Hurwitz equivalent to a factorisation
into squares of quasireflections, \ie, of the form $\mbft'= t_1' \cdot t_2' \cdot \ldots \cdot t_l'$
with $t'_{2i-1} = t'_{2i}\in \calt(\cals)$.
\end{thma}

\noindent
{\bf Remark.}
This theorem is a generalisation of a classical result of Hurwitz. Namely,
the Coxeter--Weyl group corresponding to the usual braid group $\br_d$
is the symmetric group $\sym_d$, and the claim is equivalent in this case
to the irreducibility of the space of branched coverings $f: C \to\cp^1$ of fixed
degree and genus of the curve $C$. For the case of finite Weyl groups
corresponding to simple complex Lie algebras, this result was obtained by Kanev
\cite{Kan}.

\medskip
In its turn, \refthm{thm6} is used to describe the elements in $\sfp(\cals)$ that
can be represented as products of quasireflections $t_i$ projecting 
into a given subset $\calt' \subset \calt(\cals)$. Clearly, one can
replace $\calt'$ by the subgroup $G\subset \sfw(\cals)$ generated by it.
For such $G$, we denote by $\sfp_\ab(\cals)_G$ the group of \slsf{coinvariants}
(see, \eg, \cite{Bro}). Since $\sfw(\cals)$ acts on $\sfp_\ab(\cals)$
permuting the basis $A_t$, $t\in \calt(\cals)$, it follows that $\sfp_\ab(\cals)_G$
has a natural basis formed by the set of $\sfw(\cals)$-orbits in $\calt(\cals)$.
We denote the elements of this basis by $A_{G\cdot t}$.

\newthma{thm7}{\refthm{thm-subgr}} Let $G$ be a subgroup of $\,\sfw(\cals)$
and $x$ an element of $\,\sfp(\cals)$ which can be represented as a product $x
= \prod_i \hat t_i^{\epsilon_i} \cdot \prod_j [x_{2j-1}, x_{2j}]$ such that $\hat t_i$ are
quasigenerators, $\epsilon_i=\pm1$, and $[x_{2j-1}, x_{2j}]$ are commutators of some
$x_j \in \br(\cals)$.  Assume that the projections of $x_i$ and $\hat t_j$ to
$\sfw(\cals)$ lie in $G$.  Then the projection $[x]_G$ of $x$ to
$\sfp_\ab(\cals)_G$ lies in the free abelian group generated by basis elements
$A_{G\cdot t}$ with $t \in G \bigcap \calt(\cals)$.
\end{thma}

\noindent
{\bf Remark.}
The point is that basis elements $A_{G\cdot t}$ with $t \in\calt(\cals) \bs G $ do not appear
in the expansion of~$[x]_G\in \sfp_\ab(\cals)_G$.

\medskip
Next, we focus on the study of the structure of the braid group $\br(\cals_g)$
involved in the Wajnryb--Matsumoto presentation of the group $\map_g$ of mapping
classes for a surface $\Sigma$ of genus $g$.

In the case $g\geq 4$ the set $\calt(\cals_g)$ is infinite, and it is hard to
describe it geometrically in terms of $\Sigma$. We show that there exists a
$\sfw(\cals_g)$-equivariant map of this set onto the set of all non-zero homology
classes $v \in \sfh_1(\Sigma, \zz_2)$. The latter set is denoted by $\calh_g$ and we obtain a
$\sfw(\cals_g)$-equivariant homomorphism $\zz\langle\calt(\cals_g)\rangle \to \zz\langle\calh_g\rangle$
of free abelian groups.  Thinking of $\map_g$ as a generalised braid group,
the corresponding Coxeter--Weyl group is $\Sp(2g, \zz_2)$ and the elements of
$\calh_g$ are the corresponding quasireflections. Note that in the cases
$g=1,2,3$ we have the equalities $\calt(\cals_g) =\calh_g$ and $\sfw(\cals_g)
=\Sp(2g,\zz_2)$ (up to the centre $\zz_2$ in $\sfw(\cals_3)$ in the case
$g=3$).

Our next result describes generators of the kernel of the natural homomorphism
$\br(\cals_g)\to \Sp(2g, \zz_2)$. This gives us generators of the kernels of
the homomorphisms $\map_g\to \Sp(2g, \zz_2)$ and $\sfw(\cals_g)\to \Sp(2g, \zz_2)$. We
call the latter kernel the \slsf{Weyl--Torelli group} of $\Sigma$ and denote it by
$\IW_g$.

\newthma{thm8}{\refcorol{cor-IW}, \propo{prop-IW-ab}}\strut

\sli The kernel of $\map_g\to \Sp(2g, \zz_2)$ is generated by squares $T_\delta^2$ of
Dehn twists along non-separating curves $\delta \subset \Sigma$.

\slii The abelianisation $\IW_{g,\ab}$ of\/ $\IW_g$ is a $\zz_2$-vector space
isomorphic to $\land^6 \sfh_1(\Sigma, \zz_2)$.

\sliii There exists a natural lattice extension
\[
0 \to \zz\langle \calh_g \rangle \to \Lambda_g \to \IW_{g,\ab} \to 0
\]
such that $\Lambda_g$ can be realised as a sublattice in $\frac12 \zz\langle \calh_g \rangle$
and the image of the induced embedding $\IW_{g,\ab} \subset \frac12 \zz\langle
\calh_g \rangle /\zz\langle \calh_g \rangle \cong \zz_2\langle \calh_g \rangle$ is generated by the sums
$L_V := \frac12 \sum_{v \neq 0 \in V} A_v$ where $V$ is a symplectic $6$-dimensional subspace of\/
$\sfh_1(\Sigma, \zz_2)$.

\sliv There exists a group extension
\[
0 \to \Lambda_g \to \wh\Sp(2g, \zz_2)  \to \Sp(2g, \zz_2) \to 1
\]
and a homomorphism of the extension $1 \to \sfp(\cals_g) \to \br(\cals_g) \to
\sfw(\cals_g) \to1$ onto this extension such that the homomorphism
$\sfw(\cals_g) \to \Sp(2g, \zz_2)$ has the usual meaning and $\sfp(\cals_g)
\to\Lambda_g$ is the composition $\sfp(\cals_g) \to \sfp_\ab(\cals_g) \to \zz\langle \calh_g \rangle \subset
\Lambda_g $.
\end{thma}

As a result, we can perform our calculations in much smaller and geometrically
transparent groups $\Lambda_g$ and $\wh\Sp(2g, \zz_2)$, rather than in $\map_g$ or
in $\br(\cals_g)$. Basically, this means working in $\Lambda_g$ considered as
an $\Sp(2g, \zz_2)$-module.

\smallskip%
Next we describe the image of the Wajnryb--Matsumoto relations in $\Lambda_g$. For
this purpose, for every $\mu \in \sfh_1(\Sigma, \zz_2)$ we denote by $\varphi_\mu: \sfh_1(\Sigma,
\zz_2)\to \zz_2$ the homomorphism given by $v \mapsto v \cap \mu$ and define a
homomorphism $\hat \varphi_\mu: \zz\langle \calh_g \rangle \to \zz$ by setting $\hat \varphi_\mu(A_v) := 1$ for
every $v\in \calh_g \subset \sfh_1(\Sigma, \zz_2)$ with $\varphi_\mu(v) = 1$ and $\hat \varphi_\mu(A_v) :=
0$ otherwise. Thus $\hat \varphi_\mu$ is the ``counting function'' for the
set-theoretic map $\varphi_\mu: \calh_g \to \zz_2$. Our result describes the mapping
class group $\map_{g,1}$ of the surface $\Sigma$ with one point fixed. It is known
that for all $g\geq2$ the kernel of the natural homomorphism $\map_{g,1} \to
\map_g$ is the fundamental group $\pi_1(\Sigma)$.

\newthma{thm9}{\refcorol{cor-lant-ext}} The homomorphism $\hat \varphi_\mu:
\zz\langle\calh_g\rangle \to \zz$ extends naturally to a homomorphism from the
stabiliser group $\br(\cals_g)_\mu$ of the element $\mu$. The latter induces a
homomorphism $\hat \varphi_\mu: \map_{g,1,\mu} \to \zz_4$. Moreover, for every $\gamma \in
\pi_1(\Sigma)$ and the corresponding element $x_\gamma \in \map_{g,1}$ one has $\hat
\varphi_\mu(x_\gamma) \equiv (2g-2)\cdot \varphi_\mu(\gamma) \mod 4$.
\end{thma}

This result is used to obtain

\newthma{thm10}{} Every topological Lefschetz fibration $\pr: X\to Y$ with even
genus $g= 2g'$ of the fibre $F:= \pr\inv(y)$ admits a $\zz_2$-section, \ie, a
class $[\sigma] \in \sfh_2(X, \zz_2)$ with $[\sigma] \cap [F] \equiv 1 \mod2$.
\end{thma}

\noindent
{\bf Remark.}
This result is proved in \refthm{thm-mu0} in the form of a certain
factorisation problem in $\map_{g,1}$.

\medskip
Besides, we prove  the following in \refcorol{cor-p3}.

\newthma{thm11}{} For $g\geq2$ there exists an embedding
$\sfh_1(\Sigma, \zz_2) \hrar\sftt^3 \big( \sfh_1(\Sigma, \zz_2) \big) =: \sftt^3$
of the $\zz_2$-homology space into its third tensor power
and a homomorphism $\wp$ from the kernel $\ker \big( \map_{g,1} \to
\Sp(2g,\zz_2) \big)$ to $\sftt^3$ such that the composition of $\wp$ with the
natural embedding $ \pi_1(\Sigma) \hrar \map_{g,1}$ equals the composition of the
projection $\pi_1(\Sigma) \to \sfh_1(\Sigma, \zz_2)$ with the above embedding $\sfh_1(\Sigma, \zz_2)
\hrar \sftt^3$.
\end{thma}

Now let us give an explanation of how \refthm{thm7} is involved in the proof of the main result.
Recall that by Gervais (\cite{Ger}, \slsf{Theorem C}) the group $\map_g$ with $g\geq3$ has a universal
central extension $\wt\map_g$ whose kernel is $\zz$. There exists a natural lift of the monodromy
of a topological Lefschetz pencil $\pr: X \to Y$ to the extension $\wt\map_g$ whose evaluation
is an element $c$ in the kernel, i.\,e., an integer. This extension can be obtained
as the pull-back of the universal covering
\[
0 \to \zz \to \wt\Sp(2g, \rr) \to \Sp(2g, \rr) \to 0
\]
and one can show that the number $c$ is --- roughly speaking --- the Chern number $c_1$
of the symplectic vector sheaf $R^1\pr_*\underline{\rr}{}_X$, the first derived image
of the constant sheaf $\underline{\rr}{}_X$ over $X$ (compare with \cite{ABKP,BKP,Bo-Tsch}).
Recall that the stalk of $R^1\pr_*\underline{\rr}{}_X$ over a regular value $y \in Y$ is simply
the cohomology group $\sfh^1(X_y, \rr)$ of the fibre $X_y :=\pr\inv(y) \cong \Sigma_g$
and the symplectic structure on $\sfh^1(X_y, \rr)$ is given by the cup product.
The description of the extension $\wt\map_g$ given by Gervais shows that the natural
generator of the kernel $\zz$ is the \slsf{chain relation element} (see \refsubsection{map-mu-pp})
so that $c_1$ is the algebraic number of times that the chain relation is involved in the evaluation
of the monodromy. Moreover, the universality of the extension $\wt\map_g$ implies that one cannot find
any further invariants.

Now let $G$ be a subgroup of $\map_g$. Then there could exist a non-trivial
central extension $0 \to A \to \wh G \to G \to0$ with a ``larger'' abelian group $A$,
not induced by the extension $\wt\map_g$. Then an appropriate homomorphism $\psi:
A \to \zz$ or $\psi: A \to \zz_m$ would give a (new) ``characteristic class'' for
Lefschetz fibrations with monodromy in $G$. \refthm{thm7} explains how one
can construct such an extension $0 \to A \to \wh G \to G \to0$ in the case when $G$
is the pre-image of some subgroup $H< \Sp(2g, \zz_2)$ with respect to
the natural homomorphism $\map_g \to \Sp(2g, \zz_2)$. Then the homomorphism
$\hat\varphi_\mu$ from \refthm{thm9} is essentially an example of such a~$\psi$.

\medskip%

\newsubsection[nota]{Notation} Let $X$ be a 4-manifold and $Y$ an oriented
surface.  For a fibration $\pr: X \to Y$  and a subset $U \subset Y$,
we set $X_U := \pr\inv(U)$. Thus $X_y = \pr\inv(y)$ is the fibre over
a point $y\in Y$.

Let $\Sigma$ be a CW-complex and $f:\Sigma \to \Sigma$ a continuous map. The \slsf{torus of the
 map $f$} is the quotient space $Z$ of the cylinder $\Sigma \times [0,1]$ with respect
to the equivalence relation $(z, 1) \sim (f(z), 0)$. In this paper $\Sigma$ will be a
real surface and $f: \Sigma\to \Sigma$ a diffeomorphism. Observe that the map $\pr: [z, t]
\in \big(\Sigma \times [0,1] / {\sim}\big) \mapsto t \in S^1 = [0,1] /(0\sim1)$ is a well-defined
projection $\pr: Z \to S^1$ of a fibre bundle structure on $Z$ with fibre $\Sigma$.
We consider it a part of the map torus structure on $Z$. Vice versa, for any
fibre bundle structure $\pr: Z \to S^1$ with fibre $\Sigma$ there exists a continuous
map $f:\Sigma \to \Sigma$ such that $(Z, \pr)$ is the torus of $f:\Sigma \to \Sigma$. In this case
$f:\Sigma \to \Sigma$ is called the \slsf{monodromy map of $\pr: Z \to S^1$} or simply the
\slsf{monodromy of $\pr: Z \to S^1$}. The monodromy is defined up to isotopy and
two fibre bundles $\pr: Z \to S^1$ and $\pr': Z' \to S^1$ are isomorphic iff their
monodromies are isotopic.  More generally, let $\pr: Z \to Y$ be a fibre bundle with
fibre $\Sigma$ and $\gamma: S^1 \to Y$ a loop. Then the monodromy of the pulled-back
bundle $\gamma^*Z := Z \times_{Y, \gamma} S^1$ is called the \slsf{monodromy of $\pr: Z \to Y$
 along $\gamma$}. The monodromy defines the \slsf{monodromy homomorphism} between
$\pi_1(Y, y_0)$ and the mapping class group of $\Sigma$.

Let $\pr: Z \to S^1$ be a bundle with fibre $\Sigma$, $y_0 \in S^1$ a base point, and
$\Gamma \subset Z$ a section of $\pr$, so that $\pr: \Gamma \to S^1$ is a homeomorphism. Set $\Sigma=
\pr\inv(y_0)$ and $z_0 := \Gamma \bigcap \Sigma$. Then $(Z, \pr)$ can be realised as the torus
of a map $f\in \diff(\Sigma, z_0)$ such that $\Gamma \subset Z$ becomes the torus of the map $f:
\{z_0\} \to \{z_0\}$. Moreover, such $f$ is unique up to isotopy in $\diff(\Sigma, z_0)$,
and we call $f\in \diff(\Sigma, z_0)$ the \slsf{section monodromy of $\pr: Z \to S^1$}.
Similar notation applies to bundles $\pr: Z \to Y$ with fibre $\Sigma$.  In our case
$\Sigma$ is oriented and $f$ is orientation preserving, $f\in \diff_+(\Sigma, z_0)$.

\medskip
An embedded circle $\delta$ on a connected surface $\Sigma$ is \slsf{non-essential} if
it bounds a disc in $\Sigma$ and \slsf{essential} otherwise. \label{ess} All the
embedded circles used in the proofs will be essential. An embedded circle $\delta$ on
$\Sigma$ is \slsf{separating} if $\Sigma \bs \delta$ has two connected components and
\slsf{non-separating} otherwise. If $\Sigma$ is closed and orientable, then $\delta$ is
non-separating if and only if $[\delta] \neq 0 \in \sfh_1(\Sigma, \zz)$ or, equivalently,
if and only if $[\delta] \neq 0 \in \sfh_1(\Sigma, \zz_2)$

\newsubsection[thanks]{Acknowledgments} The author would like to thank Stefan
Nemirovski for numerous discussions.
Dmitry Akhiezer,
Denis Auroux,
Kai Cieliebak,
Yasha Eliashberg,
Ursula Hamenst{\"a}dt,
Sergei Ivashkovich,
Viatcheslav Kharlamov,
Dieter Kotschick,
Stepan Orevkov,
Leonid Polterovich,
and
Wolfgang S{\"o}rgel
have made many valuable comments on earlier drafts of the paper.

\newsection[prel]{Topology of Lefschetz fibrations}

\newsubsection[sy-le-fi]{Lagrangian submanifolds in Lefschetz fibrations}
We refer to \cite{Go-St} for the definition of
\slsf{topological Lefschetz fibrations} $\pr: X \to Y$, \slsf{topological
 Lefschetz pencils} $\pr: X \to S^2$, and \slsf{(symplectic) blow-ups} of such
$X$, and to \cite{AMP} and \cite{Do} for a general theory of symplectic
Lefschetz fibrations and their Lagrangian submanifolds. The paper \cite{N-1}
contains a discussion of Lagrangian submanifolds in projective algebraic manifolds.
In this paper only the special case of \emph{4-dimensional} ambient
manifolds $X$ will be considered. We shall always assume that the base surface $Y$
of a Lefschetz fibration is oriented.

\newdefi{df1} Let $X$ be a closed $4$-manifold and $L\subset X$ a closed embedded
surface.  A \slsf{Morse--Lefschetz fibration} of $(X, L)$ is a Lefschetz
fibration $\pr: X \to Y$ with the following properties:
\begin{enumerate}
\item[\sli] the image $\pr(L)$ is a smooth embedded curve $\Gamma$ in $Y$;
\item[\slii] all critical points of the restricted projection $\pr_L :=
 \pr|_L: L \to \Gamma$ are of Morse type;
\item[\sliii] a point $x \in \Gamma$ is a critical value of $\pr_L : L \to \Gamma$ if and only if it
 is a critical value of $\pr: X \to Y$.
\end{enumerate}
\end{defi}

Such a curve $\Gamma$ must be an arc or a circle. 
It turns out that any given symplectic manifold $(X, \omega)$ and any given
Lagrangian embedded closed surface $L\subset X$ can be included in a SMLF
over the two-sphere $S^2$. Namely, the result proved in \cite{AMP} ensures the following.

\newthm{thmMLF} Let $(X, \omega)$ be a closed symplectic $4$-manifold such that
$[\omega] = c_1(\scrl)$ for some line bundle $\scrl$ on $X$, $L \subset X$ a closed
Lagrangian submanifold, and $f: L \to S^1$ a Morse function. Then there exists a
symplectic blow-up $\wh X \to X$ at a finite locus $B = \{ b_1,\ldots, b_N\} \subset X
\bs L$, a Morse--Lefschetz fibration $\pr: \wh X \to S^2$ of $(X, \omega, L)$ and an
embedding $\gamma: S^1 \hrar S^2$, such that $\pr|_L = \gamma\circ f$.
\end{thm}

Observe that in the case of a usual Morse function $f:L \to\rr$ the image $f(L)$
is a compact interval which can be embedded in the circle $S^1$. So this case
reduces to the case of circle valued Morse functions.

\smallskip%
Let us establish the version of \refthm{thmMLF} needed for our purpose.

\newlemma{lem-slf} For any embedded Lagrangian Klein bottle $K$ in a
symplectic $4$-manifold $(X, \omega)$ there exists a deformation $\ti \omega$ of the
symplectic form $\omega$, a symplectic blow-up $\wh X$ of $X$ at several points
$x_1,\ldots,x_N \in X$ outside $K$, and a symplectic Lefschetz fibration $\pr: \wh X
\to S^2$ with the following properties:
\begin{itemize}
\item[(ML1)] the projection $\pr(K)$ is an embedded circle $\Gamma \subset S^2$;
\item[(ML2)] no critical value $y \in S^2$ of $\pr: X \to S^2$ lies on $\Gamma$;
\item[(ML3)] the fibres of the restricted projection $\pr|_K:K \to \Gamma$ are
 circles;
\item[(ML4)] there exists a section of the projection $\pr: X \to S^2$ {\rm (}\ie, a
 surface $S \subset X$ whose projection $\pr|_S:S \to S^2$ is a diffeomorphism\/{\rm )}
 such  that $S \bigcap K = \emptyset$.
\end{itemize}
\end{lem}

\proof The deformation of $\omega$ is needed to make its cohomology class $[\omega]$
rational. It is easy to show that for any fixed Riemannian metric $g_X$ on $X$
and any $\varepsilon>0$ there exists a closed 2-form $\ti \omega$ such that the class
$[\ti\omega]$ is rational and $\norm{\ti\omega -\omega}_{C^0(X)} \leq \varepsilon$. The latter condition
ensures that $\ti\omega$ is also a symplectic form.

Since $\sfh^2(K, \rr)=0$, the restriction $\ti\omega|_K$ is exact. Moreover,
by Hodge theory, there exists a 1-form $\alpha$ on $K$ with $d\alpha = \ti\omega|_K$ and
\[
\norm{\alpha}_{C^0(K)} \leq C\cdot \norm{\ti\omega|_K}_{C^0(K)} = C\cdot \norm{(\ti\omega
 -\omega)|_K}_{C^0(K)} \leq C\cdot \varepsilon
\]
with the constant $C$ depending only on the embedding $K\subset X$ and the metric
$g_X$. By the (generalised) Darboux theorem (see, \eg, \cite{A-G}), there
exists a neighbourhood $U$ of $K$ in $X$ which is symplectomorphic relative $K$
to a neighbourhood $U'$ of the zero section of the cotangent bundle $T^*K$
equipped with the canonical form $\omega_{T^*K}$. The condition $\norm{\alpha}_{C^0(K)}
\leq C\cdot\varepsilon$ ensures that the graph $L_\alpha\subset T^*K$ of the form $\alpha$ lies in $U'$.
Observe that $\omega_{T^*K}|_{L_\alpha} = d\alpha = \ti\omega$. It follows that there exists a
diffeomorphism $\Phi: U \to U$ such that $\Phi^*\ti\omega|_K = \omega|_K=0$ and such that $\Phi$ is
trivial near the boundary $\dd U$. Replacing $\ti\omega$ by $\Phi^*\ti\omega$,
we see that $K$ is $\wt\omega$-Lagrangian.

Now, the desired symplectic Lefschetz fibration is obtained by \refthm{thmMLF}
applied to the map $\varphi: K\to S^1$ which realises $K$ as an $S^1$-fibre bundle over
$S^1$. The number of single blow-ups in $\wh X \to X$ must be large for the
following reason. By the construction of the symplectic Lefschetz fibration
$\pr: \wh X\to S^2$ in \cite{AMP}, the blow-up centres $x_1,\ldots,x_N$ are the
common zeroes of a pair of sections $s_0, s_1$ of a certain line bundle $\scrl$
on $X$ with $c_1(\scrl) = k [\ti \omega]$ and $k\gg0$. Moreover, the intersection
index of the zero loci $s_0\inv(0)$ and $s_1\inv(0)$  at each $x_1,\ldots,x_N$ is equal to~1,
so that $N = k^2[\ti \omega]^2 \gg0$.
\qed

\newsubsection[loc-slf]{Topology of Lefschetz fibrations at
 singular fibres} Let $\Sigma$ be an oriented surface and $\delta \subset \Sigma$ an embedded
circle. Fix an annular neighbourhood $U$ of $\delta$ and realise it as the annulus
$\{(\rho, \theta) : \frac12 <\rho<2, 0\leq \theta \leq2\pi \}$ with the orientation given by $d\rho \land d\theta$.
Let $\chi(\rho)$ be a non-decreasing function with $\chi(\rho) \equiv 0$ for $\rho \leq\half$ and $\chi(\rho)\equiv 2\pi$ for
$\rho \geq 2$.

\newdefi{def-Dehn} The \slsf{positive Dehn twist} of $\Sigma$ along $\delta$ is the
diffeomorphism $T_\delta: \Sigma \to \Sigma$ which is identical outside the neighbourhood $U$
above and is given by the formula
\[
T_\delta(\rho,\theta) = (\rho, \theta - \chi(\rho))
\]
inside $U$. The \slsf{negative Dehn twist} is the inverse diffeomorphism, it is given
by $(\rho,\theta) \mapsto (\rho, \theta + \chi(\rho))$.
\end{defi}

The Dehn twist of a prescribed sign is unique up to isotopy. The sign of
a Dehn twist is defined  so  that symplectic and usual algebraic Lefschetz fibrations
have only \emph{positive} Dehn twists in their monodromy.

The positivity of the Dehn twists in symplectic Lefschetz fibrations plays no r{\^o}le in the proof.
This is the reason why the result can be generalised to topological Lefschetz fibrations.

The notation $T_\delta$ will be used both for the Dehn twist along a prescribed embedded
circle $\delta \subset \Sigma$ with a given sign and some specified neighbourhood $U$ and
coordinates $(\rho, \theta)$ on it and for the isotopy class of such a twist inside the
corresponding mapping class group, see \refsubsection{ext} below.

\smallskip
Let $\Delta$ be a disc with the boundary $\partial\Delta =: \gamma$, $Z$ a $4$-manifold, $\pr: Z \to
\Delta$ a proper Lefschetz fibration with a unique critical point $z^*$ over the
origin $0\in \Delta$. Denote by $\Sigma$ the generic fibre of $\pr: Z \to \Delta$, say over $y_0
\in \gamma$.

\newlemma{lem-loc}
\begin{enumerate}
\item The manifold $Z$ can be deformationally retracted on the singular fibre
 $\Sigma^* := \pr\inv(0)$.
\item The monodromy along the boundary circle $\dd\Delta$ acts on $\Sigma$ as the
 Dehn twist along a certain embedded circle $\delta \subset \Sigma$, and the
 singular fibre $\Sigma^*$ is obtained from $\Sigma$ by contracting $\delta$ to
 the \slsf{nodal point} on $\Sigma^*$.
\item If $\delta$ is non-separating, then the homology group $\sfh_2(\Sigma^*, \zz)$
 is isomorphic to\/ $\zz$ and generated by the image of the fundamental class $[\Sigma] \in
 \sfh_2(\Sigma, \zz)$; otherwise, the group $\sfh_2(\Sigma^*, \zz)$ is $\zz \oplus \zz$
 generated by the fundamental cycles of the two irreducible components of $\Sigma^*$.
\item The homology group $\sfh_1(\Sigma^*, \zz)$ is the quotient of $\sfh_1(\Sigma,
 \zz)$ by the subgroup $\zz\langle [\delta]\rangle$ generated by $[\delta]$, and the
 cohomology group $\sfh^1(\Sigma^*, \zz)$ is the orthogonal submodule $\zz\langle
 [\delta]\rangle^\perp \subset \sfh_1(\Sigma, \zz)$ with respect to the natural pairing $\sfh_1(\Sigma,
 \zz) \times \sfh_1(\Sigma, \zz)\to \zz$.
\item \label{ler} The homology groups of the boundary $\dd Z = Z_\gamma$ can be
 included into exact sequences
\[
\xymatrix@C-7pt@R-18pt{
0 \ar[r] &
\sfh_0( \gamma, \scrh_1(Z_y, \zz)) \ar[r] &
\sfh_1(Z_\gamma, \zz) \ar[r] &
\sfh_1( \gamma, \scrh_0(Z_y, \zz)) \ar[r] &
0
\\
0 \ar[r] &
\sfh_0( \gamma, \scrh_2(Z_y, \zz)) \ar[r] &
\sfh_2(Z_\gamma, \zz) \ar[r] &
\sfh_1( \gamma, \scrh_1(Z_y, \zz)) \ar[r] &
0
}
\]
in which $\scrh_p(Z_y, \zz)$ denotes the locally constant sheaf with the stalk
$\sfh_p(Z_y, \zz)$ over $y\in \gamma$. In particular, we have natural isomorphisms
\[
\sfh_1(\gamma, \zz) \xrar{\;\cong\;} \sfh_1(\gamma, \scrh_0(Z_y, \zz))
\quad\text{and}\quad
\sfh_0( \gamma, \scrh_2(Z_y, \zz)) \xrar{\;\cong\;} \sfh_2(\Sigma, \zz),
\]
the group $\sfh_1( \gamma, \scrh_1(Z_y, \zz))$ is the subgroup $[\delta]^\perp$ of\/
$\sfh_1(\Sigma, \zz)$, \ie, it consists of $\lambda \in \sfh_1(\Sigma, \zz)$ with $\lambda\cap  [\delta]
 =0$; $\sfh_0( \gamma, \scrh_1(Z_y, \zz))$ is the quotient group of $\sfh_1(\Sigma,
\zz)$ with respect to the subgroup $\zz\langle [\delta] \rangle$ generated by  $[\delta]$.
\end{enumerate}
\end{lem}

For the proof we refer to \cite{AGV}. 
The exact sequences in (\ref{ler}) are obtained from the Leray spectral
sequence of the projection $\pr: Z_\gamma \to\gamma$.  The circle $\delta$ is called the
\slsf{vanishing cycle} of the monodromy at $z^*$, and its homology class $[\delta] \in
\sfh_1(\Sigma, \zz)$ is called the \slsf{vanishing class}.

\newsubsection[ext]{Homotopy type of groups of surface diffeomorphisms}
Let $\Sigma$ be a closed oriented surface of genus $g=g(\Sigma)$ with the base point
$z_0 \in \Sigma$. Denote by $\diff_+(\Sigma)$ the group of orientation preserving
diffeomorphisms of $\Sigma$, by $\diff_+(\Sigma,z_0)$ the subgroup of diffeomorphisms
preserving the base point $z_0$, and by $\diff_+(\Sigma, [z_0])$ the subgroup of
diffeomorphisms preserving the base point $z_0$ and acting trivially on the
tangent plane $T_{z_0}\Sigma$. Denote by $\map_g$, $\map_{g,1}$, and $\map_{g,[1]}$
the corresponding mapping class groups, \ie, the groups of connected components
of $\diff_+(\Sigma)$, $\diff_+(\Sigma, z_0)$, and $\diff_+(\Sigma, [z_0])$. Observe that the
natural action of $\diff_+(\Sigma)$ on $\Sigma$ induces the principle fibre bundle
$\ev_{z_0}: \diff_+(\Sigma) \to \Sigma$ given by $f \in \diff_+(\Sigma) \mapsto f(z_0) \in \Sigma$ with the
structure group $\diff_+(\Sigma, z_0)$. In this way we obtain the long exact
sequence of homotopy groups:
\begin{equation}\eqqno(hmtpy)
\begin{split}
\cdots \to \pi_{k+1}(\Sigma)& \xrar{\;\dd\;}
   \pi_k(\diff_+(\Sigma, z_0)) \to \pi_k(\diff_+(\Sigma)) \to \pi_k(\Sigma) \xrar{\;\dd\;}  \cdots\\
\cdots  \to \pi_1(\Sigma)& \xrar{\;\dd_\pi\;} \pi_0(\diff_+(\Sigma, z_0)) \to \pi_0(\diff_+(\Sigma)) \to 1.
\end{split}
\end{equation}
Similarly, the natural action of $\diff_+(\Sigma, z_0)$ on $T_{z_0}\Sigma$ by means of
the differential $D_{z_0}f: T_{z_0} \to T_{z_0}$ of a given $f \in \diff_+(\Sigma,z_0)$
induces the principle fibre bundle $D_{z_0}: \diff_+(\Sigma, z_0) \to \Gl_+(T_{z_0}\Sigma)$
with the structure group $\diff_+(\Sigma, [z_0])$. In this way we obtain a similar long
exact sequence of homotopy groups:
\begin{equation}\eqqno(hmtpy-frame)
\begin{split}
\cdots \to \pi_{k+1}(\Gl_+(2,\rr))& \xrar{\;\dd\;}
   \pi_k(\diff_+(\Sigma, [z_0])) \to \pi_k(\diff_+(\Sigma, z_0)) \to   \cdots\\
\cdots  \to \pi_1(\Gl_+(2,\rr))& \xrar{\;\dd\;} \pi_0(\diff_+(\Sigma, [z_0])) \to \pi_0(\diff_+(\Sigma,z_0)) \to 1.
\end{split}
\end{equation}

The following facts about mapping class groups are well-known, see
\cite{Bi-2} or \cite{Iva}.

\newprop{prop-hm1}  
\sli In the cases $g=0$ and $g=1$ the map $\pi_0(\diff_+(\Sigma, z_0))
\to \pi_0(\diff_+(\Sigma))$ is an isomorphism.

\slii In the case $g\geq 2$ the connected component $\diff_0(\Sigma)$ of $\diff_+(\Sigma)$
is contractible and the sequence \eqqref(hmtpy) induces the exact sequence of
groups
\begin{equation}\eqqno(ext-g1)
1 \to \pi_1(\Sigma) \xrar{\;\dd_\pi\;} \map_{g,1} \to \map_g \to 1.
\end{equation}
Moreover, the group $\map_{g,1}$ is the subgroup of index $2$ of the group
$\aut(\pi_1(\Sigma))$ and consists of those automorphisms of $\pi_1(\Sigma)$ which act
trivially on $\sfh_2(\pi_1(\Sigma), \zz)$, the map $\pi_1(\Sigma) \xrar{\;\dd_\pi\;} \map_{g,1}$
associates with each $\gamma \in \pi_1(\Sigma)$ the inner homomorphism $i_\gamma: \alpha \in \pi_1(\Sigma) \mapsto
\gamma\alpha\gamma\inv \in \pi_1(\Sigma)$, and the group $\map_g$ is the corresponding group of outer
automorphisms of $\pi_1(\Sigma)$.

\sliii For $g\geq1$, the long exact sequence \eqqref(hmtpy-frame) induces the
central extension
\begin{equation}\eqqno(ext-g2)
1 \to \zz=\pi_1(\Gl_+(2,\rr))  \xrar{\;\dd\;} \map_{g,[1]} \to \map_{g,1} \to 1.
\end{equation}

\sliv $\pi_1(\diff_+(S^2)) = \zz_2$, $\pi_1(\diff_+(T^2)) \cong \pi_1(T^2) \cong  \zz^2$.
\end{prop}

The embedding $\pi_1(\Sigma, z_0) \xrar{\;\dd_\pi\;} \map_{g,1}$ (for $g\geq2$) is easy to
describe geometrically. Namely, for any $\gamma \in \pi_1(\Sigma, z_0)$ we choose an isotopy
$F_t: \Sigma \to \Sigma$ with $t\in [0,1]$ such that $F_0 \equiv \id_\Sigma$ and such that the curve
$t\in [0,1] \mapsto F_t(z_0)$ represents the class $\gamma \in \pi_1(\Sigma, z_0)$. It follows
immediately from the definition that $F_1$ represents the image $\dd_\pi(\gamma) \in
\map_{g,1}$. In the special case when $\gamma \in \pi_1(\Sigma, z_0)$ is represented by a
smooth embedded curve, still denoted by $\gamma$, one can give a more explicit
description of $\dd_\pi(\gamma)$. Namely, let $\gamma_+$ and $\gamma_-$ be the
right and left boundary components of a collar neighbourhood $U_\gamma$ of $\gamma$ in
$\Sigma$. Then the Dehn twists $T_{\gamma_\pm}$ of $\Sigma$ along $\gamma_\pm$ are well-defined
elements of the group $\map_{g,1}$ and, as it can be easily seen,
\begin{equation}\eqqno(d-gam)
\dd_\pi(\gamma) = T_{\gamma_+} \circ T_{\gamma_-}\inv.
\end{equation}
It follows that $\map_{g,1}$ is generated by  Dehn twists along smooth
embedded curves $\gamma$ on $\Sigma$ which do not pass through $z_0$.

\newsubsection[hmlg-TLF]{Homology of Lefschetz fibrations}
In this subsection, we describe the necessary conditions for the homological triviality
of a Morse--Lefschetz embedding of the Klein bottle in a Lefschetz fibration $\pr: X \to Y$.
We shall consider a somewhat more general situation than the one obtained by applying~\lemma{lem-slf}.

Let $\pr: X \to Y$ be a topological Lefschetz fibration with generic fibre $\Sigma$
and oriented base $Y$. Assume that $K \subset X$ is an embedded Klein bottle such that
\begin{enumerate}
\item[(T1)] the image $\pr(K)$ is an embedded  circle $\Gamma \subset Y$
 containing no critical value of $\pr$;
\item[(T2)] the restricted projection $\pr|_K: K \to \Gamma$ is a bundle with the
 fibre $S^1$.
\end{enumerate}
These conditions ensure that $(X,\pr, K, \Gamma)$ is a Morse--Lefschetz fibration of a special form.

\smallskip
Identify $\Sigma$ with the fibre $\pr\inv(y_0)$ over a fixed point $y_0 \in \Gamma$.
Let $\mbfm$ be the fibre of $K$ over $y_0$ so that $\mbfm = \Sigma \bigcap K$.
We call $\mbfm$ the \slsf{meridian circle} of $K$ but consider it mostly as a curve
on $\Sigma$. Let $F_\Gamma: \Sigma \to\Sigma$ be the monodromy of the
projection $\pr: X_\Gamma \to \Gamma$. Recall that $F_\Gamma$ is defined up to
isotopy. Consequently, we can choose $F_\Gamma$ so that $F_\Gamma(\mbfm)=\mbfm$.

\medskip
Let us consider an easier special case first.

\newprop{prop-mu-0} Assume that $[\mbfm] = 0 \in \sfh_1(\Sigma, \zz_2)$ so that
$\mbfm$ separates $\Sigma$. Then the genus $g$ of $\Sigma$ is even and $K$ is
$\zz_2$-homologous to the fibre $\Sigma = \pr\inv(y_0)$.
\end{prop}

\noindent
{\bf Remark.}
In the situation of the \slsf{Main Theorem} (once we have applied~\lemma{lem-slf}), the fibre
of the Lefschetz fibration is $\zz_2$-homologically non-trivial because it has intersection
index $1$ with any exceptional section. Thus, \propo{prop-mu-0} allows one to prove the
\slsf{Main Theorem} in the (easy) case when the meridian circle $\bfm$ is homologically trivial
in the fibre of a Morse--Lefschetz fibration for~$(X,K)$. In fact, it will be shown in
\refsubsection{map-g-fac} that the fibre is $\zz_2$-homologically non-trivial for any topological
Lefschetz fibration with even fibre genus (see~\refthm{thm-mu0} and~\refthm{thm10}).

\proof By assumption, the meridian circle $\mbfm$ divides $\Sigma$ into two pieces,
$\Sigma \bs \mbfm = \Sigma' \sqcup \Sigma''$. Since the monodromy $F_\Gamma$
preserves the orientation on $\Sigma$ and inverts the orientation of $\mbfm$, $F_\Gamma$
must interchange the pieces $\Sigma'$ and $\Sigma''$. Realising $X_\Gamma$ as the torus of
the monodromy $F_\Gamma :\Sigma \to\Sigma$, we see that the boundary of the $\zz_2$-chain $\Sigma' \times
[0,1]$ in $X_\Gamma$ is $[K] + [\Sigma]$. Further, the pieces $\Sigma'$ and $\Sigma''$ have
the same genus $g'$, hence the genus of $\Sigma$ is even, $g =2g'$.
\qed

\bigskip
Let us now tackle the more complicated case when the meridian is homologically non-trivial
in the fibre of the Lefschetz fibration. Denote by $\mu :=[\mbfm] \in \sfh_1(\Sigma, \zz_2)$
the homology class of the meridian and by $y_i^*$, $i=1,\ldots,n$,
the critical values of $\pr:X \to Y$ (there may be none). For
every $y^*_i$ fix a small disc $D_i$ containing $y^*_i$ and set $\Gamma_i
:= \dd D_i$. The latter is an embedded curve surrounding
$y^*_i$. Denote $Y^\circ := Y \bs \big( \cup_i D_i \big)$.  Clearly, $Y^\circ$ is
homotopy equivalent to the complement of the set of critical values.
For any subset $A \subset Y$, let $X_A := \pr\inv(A)$ be the part of $X$ lying
over $A$.  The set $X_{Y^\circ} = \pr\inv(Y^\circ)$ is denoted by $X^\circ$.  In
the case of separating $\Gamma$ we denote by $Y_+$, $Y_-$  the resulting
pieces of $Y$ and by $g_Y^+$, $g_Y^-$ their genera. Set $Y^\circ_\pm := Y_\pm
\bigcap Y^\circ$ and $X^\circ_\pm := X_{Y^\circ_\pm} = \pr\inv(Y^\circ_\pm)$,
the latter are the parts of $X$ lying over $Y^\circ_\pm$.

The homological spectral sequence for the fibre bundle $\pr: X^\circ \to Y^\circ$
degenerates at the term $E^2_{p,q}$ and yields the exact sequence
\[
0 \to \sfh_0( Y^\circ, \scrh_2(X_y, \zz_2)) \to \sfh_2(X^\circ, \zz_2)
\to \sfh_1( Y^\circ, \scrh_1(X_y, \zz_2)) \to 0,
\]
where $\scrh_p(X_y, \zz_2)$ denotes the locally constant sheaf with the stalk
$\sfh_p(X_y, \zz_2)$ at $y \in Y^\circ$. Since the monodromy map $F_\Gamma$ inverts the
orientation of $\mbfm$ we cannot have the case $Y^\circ =S^2$. Consequently, the
surface $Y^\circ$ is an Eilenberg-MacLane $\pi_1$-space and $\sfh_p( Y^\circ,
\scrh_q(X_y, \zz_2)) = \sfh_p( \pi_1(Y^\circ), \scrh_q(X_y, \zz_2))$.

We use a special cell decomposition of $Y^\circ$ and the induced presentation of
$\pi_1(Y^\circ)$ to calculate the homology groups $\sfh_p( \pi_1(Y^\circ), M)$ with
coefficients in a given $\pi_1(Y^\circ)$-module $M$. Let $g_Y$ denotes the genus of
$Y$. The construction we use depends on whether the curve $\Gamma$ separates
$Y$ or not.

In the non-separating case we consider a polygon $C_Y$ with $2g_Y + n$
sides, $n$ being the number of critical values of $\pr: X \to Y$ and thus the
number of holes in $Y^\circ$. Mark and orient the sides of $C_Y$ according to the
relation word
\begin{equation}\eqqno(R)
R:= [\xi_1, \eta_1] \cdot \ldots \cdot [\xi_{g_Y}, \eta_{g_Y}] \cdot \Gamma_1 \cdot \ldots \cdot\Gamma_n.
\end{equation}
Gluing the first $2g_Y$ sides of $C_Y$ pairwise according to the marking and
orientation we obtain a surface $Y'$ of genus $g_Y$ with one boundary circle
divided in subsequent segments $\Gamma_1,   \ldots ,\Gamma_n$. To obtain the surface $Y^\circ$ we
divide each $\Gamma_i$ into three pieces, say $\Gamma_i'$, $\Gamma_i''$, and  $\Gamma_i'''$, and
glue   $\Gamma_i'$ to  $\Gamma_i'''$ reversing the orientation.

On the constructed surface $Y^\circ$ all the end points of the sides $\xi_i, \eta_i, \Gamma_j$
are identified into a single point. We take this  point as the
base point $y_0$ on $Y^\circ$ and obtain the presentation
\begin{equation}\eqqno(pi1Yo)
\pi_1(Y^\circ, y_0) = \big\langle \xi_1,\eta_1, \ldots, \xi_{g_Y}, \eta_{g_Y}; \Gamma_1 , \ldots ,\Gamma_n \;\big|\;
[\xi_1, \eta_1] \cdot \ldots \cdot [\xi_{g_Y}, \eta_{g_Y}] \cdot \Gamma_1 \cdot \ldots \cdot\Gamma_n =1  \big\rangle
\end{equation}
for the group $\pi_1(Y^\circ, y_0)$. Here we identify the curves $\xi_i,\eta_i,\Gamma_j$ with the corresponding elements in
$\pi_1(Y^\circ, y_0)$. Further, we identify the curve $\Gamma$ with $\xi_1$.

In the case of separating $\Gamma$ we modify the construction as follows.
Order the critical values $y^*_i$ so that the first $n^+$ of them lie in $Y_+$ and the
last $n^- := n-n^+$ in $Y_-$. Set
{\small%
\[ R^+ := [\xi_1, \eta_1] \cdot \ldots \cdot [\xi_{g^+_Y}, \eta_{g^+_Y}] \cdot \Gamma_1 \cdot \ldots \cdot\Gamma_{n^+}
\quad
R^- := [\xi_{g^+_Y+1},  \eta_{g^+_Y+1}] \cdot \ldots \cdot [\xi_{g_Y}, \eta_{g_Y}] \cdot\Gamma_{n^++1}\cdot \ldots \cdot\Gamma_n.
\]}%
Take polygons $C^+_Y$ and $C^-_Y$ with $2g^\pm_Y + n^\pm +1$ sides marked and
oriented according to the relation words $R^+ \cdot \Gamma^+$ and $R^- \cdot \Gamma^-$,
respectively. Glue the polygons $C^+_Y$ and $C^-_Y$ along the sides $\Gamma^+$
and $\Gamma^-$ and denote by $\Gamma$ the resulting interval. We obtain a polygon $C_Y$
with $2g_Y + n$ sides marked and oriented according to the word $R^+ \cdot R^-$.
Construct $Y^\circ$ from $C_Y$ according to this word as it was done above. As the
presentation of $\pi_1(Y^\circ)$ we take
\[
\pi_1(Y^\circ, y_0) = \big\langle \xi_1,\eta_1, \ldots, \xi_{g_Y}, \eta_{g_Y}; \Gamma_1 , \ldots ,\Gamma_n; \Gamma  \;\big|\;
R^+ = \Gamma\inv, R^- = \Gamma \big\rangle.
\]
Observe that the images of the polygons $C^+_Y, C^-_Y$ on $Y^\circ$ are glued
into pieces $Y^\circ_+, Y^\circ_-$ according to the relation words $R^+ \cdot \Gamma^+, R^- \cdot
\Gamma^-$, respectively. One can see that the fundamental groups $\pi_1(Y^\circ_\pm)$
are embedded in $\pi_1(Y^\circ)$ by identification of generators.

\smallskip
Now let $\pi$ be any group, $\pi = \langle\calx\;|\;\calr \rangle$ its presentation with the
set of generators $\calx$ and the set of relations $\calr$, and $M$ a
$\pi$-module. Let $\fr(\calx)$ be the free group generated by $\calx$. Set
$\calx \otimes M$ to be the direct sum of ``$\calx$ copies'' of $M$, so that the
elements of $\calx \otimes M$ are finite formal sums $\sum_i x_i \otimes m_i$ with $x_i \in
\calx$ and $m_i \in M$ satisfying the distributive law for the first argument.
Define $\calr \otimes M$ and $\fr(\calx) \otimes M$ in the same way. We consider
$\calx \otimes M$ and $\calr \otimes M$ as subgroups of $\fr(\calx) \otimes M$, and to
distinguish between two copies of $\fr(\calx) \otimes M$ we use the notation $x
\otimes_1 m$ and $x \otimes_2 m$, respectively.  Now set $\scrr_0 := M$, $\scrr_1 :=
\calx \otimes M$, $\scrr_2 := \calr \otimes M$ and define the \slsf{chain homomorphisms}
$\dd_1: \scrr_1 \to \scrr_0$ and $\dd_2: \scrr_2 \to \scrr_1$ by setting $\dd_1(x \otimes_1 m)
:= xm -m$ and imposing the following \slsf{derivation properties}:
\begin{itemize}
\item $\dd_2(x \otimes_2 m) = x \otimes_1 m$ if  $x\in\calx$;
\item $\dd_2(ab \otimes_2 m) = \dd_2(a \otimes_2 bm) + \dd_2(b \otimes_2 m)$ for arbitrary $a,b\in
 \fr(\calx)$.
\end{itemize}
It follows that there exists a unique $\dd_2: \scrr_2 \to \scrr_1$ with these
properties. In particular, $\dd_2(1 \otimes_2 m) = \dd_2(1\cdot 1 \otimes_2 m) = 2\dd_2(1 \otimes_2 m)$
and hence $\dd_2(1 \otimes_2 m) = 0$ for the neutral element $1 \in \fr(\calx)$.
Similarly one concludes that $\dd_2(a\inv \otimes_2 m) = -\dd_2(a \otimes_2 a\inv m)$.

The geometric meaning of these homomorphisms is as follows. From the
presentation $\pi = \langle\calx\;|\;\calr \rangle$ one constructs the \slsf{associated
 2-dimensional cell complex} $Z$ with a single $0$-dimensional cell $*$, and
with 1- and 2-dimensional cells indexed by $\calx$ and $\calr$. The
construction of $Z$ simply repeats the construction of $Y^\circ$
above. Then $\scrr_2 \xrar{\dd_2} \scrr_1 \xrar{\dd_1} \scrr_0 \to 0$ is the
cellular chain complex associated with the locally constant coefficient system
on $Z$ induced by $M$.  In particular $\dd_1 \circ \dd_2 =0$. Moreover, $Z$ is
$2$-equivalent to the Eilenberg-MacLane space $B(\pi, 1)$ and hence for $p=0$ and
$p=1$ the $p$-th homology group of the complex $(\scrr_\bullet, \dd_\bullet)$ are the
desired homology groups $\sfh_p(\pi, M)$.

By the very construction, $Y^\circ$ is homotopy equivalent to the 2-dimensional
cell complex $Z$ associated with both presentations of $\pi_1(Y^\circ)$ above. As
the coefficient module $M$ we use the homology group $\sfh_1(\Sigma,\zz_2)$. The
image of the homology class $[K]$ under projection $\sfh_2(X^\circ, \zz_2) \to
\sfh_1( Y^\circ, \scrh_1(X_y, \zz_2))$ is easy to describe: It is represented by
the chain $\Gamma \otimes_1 \mu\in \scrr_1$.

\medskip
Next, we describe the image of $\sfh_2(X^\circ, \zz_2)$ in $\sfh_2(X, \zz_2)$.
For this purpose we use the Mayer-Vietoris exact sequence associated with the
decomposition of $X$ into $X^\circ$ and $\sqcup_i X_{D_i}$. Recall that $X_{D_i} =
\pr\inv(D_i)$ and $X^\circ \bigcap X_{D_i} =X_{\Gamma_i} = \pr\inv(\Gamma_i)$. The relevant
segment of the exact sequence is
\begin{equation}\eqqno(h2-X)
\cdots \to \sfh_2(\sqcup_i X_{\Gamma_i}, \zz_2) \to
\sfh_2(X^\circ, \zz_2) \oplus \sfh_2(\sqcup_i X_{D_i}, \zz_2)
\to \sfh_2(X, \zz_2) \to  \cdots.
\end{equation}

\newlemma{lem-h2-circ} The kernel of the homomorphism $\sfh_2(X^\circ, \zz_2)
\to\sfh_2(X, \zz_2)$ is generated by classes $\Lambda_i$ whose projections to $\sfh_1(
Y^\circ, \scrh_1(X_y, \zz_2))$ have the form $\Gamma_i \otimes_1 \lambda_i$
where $\lambda_i \in \sfh_1(X_y,\zz_2)$ are $\cap$-orthogonal to $\delta_i$.
\end{lem}

\noindent
{\bf Remark.}
The condition on $\lambda_i$ is vacuous if $[\delta_i] =0 \in \sfh_1(X_y, \zz_2)$.

\proof Clearly, $\sfh_2(\sqcup_i X_{\Gamma_i}) = \oplus_i \sfh_2( X_{\Gamma_i})$.
As in the case of $X^\circ$, the homological spectral sequence for each of these summands
reduces to the exact sequence
\begin{equation}\eqqno(hml-ga-i)
0 \to \sfh_0( \Gamma_i, \scrh_2(X_y, \zz_2)) \to \sfh_2(X_{\Gamma_i}, \zz_2)
\to \sfh_1( \Gamma_i, \scrh_1(X_y, \zz_2)) \to 0.
\end{equation}
The group $\sfh_0( \Gamma_i, \scrh_2(X_y, \zz_2))$ is just $\zz_2$ generated by
$[\Sigma]$. Since  $\pi_1(\Gamma_i) =\zz$, it follows that $\sfh_1( \Gamma_i, \scrh_1(X_y,
\zz_2))$ equals $\sfh_1( \zz, M_i)$ where $M_i $ is the group $\sfh_1(X_y,\zz_2)$ with
the action of $\pi_1(\Gamma_i) =\zz$ given by the monodromy along $\Gamma_i$,
which is the Dehn twist $T_{\delta_i}$ along the vanishing cycle $\delta_i$.  Thus
$\sfh_1( \Gamma_i, \scrh_1(X_y, \zz_2))$ is simply the space $[\delta_i]^\perp$ of those
$\lambda \in \sfh_1(X_y, \zz_2)$ that are $\cap$-orthogonal to $\delta_i$. On the other hand, the
group $\sfh_0( \Gamma_i, \scrh_2(X_y, \zz_2))$ is $\zz_2$ generated by the fibre class,
and so the composition $\sfh_0( \Gamma_i, \scrh_2(X_y, \zz_2)) \to
\sfh_2(X_{\Gamma_i}, \zz_2) \to \sfh_2(X_{D_i}, \zz_2)$ is an embedding.
Consequently, the sequence \eqqref(hml-ga-i) splits and $\sfh_1( \Gamma_i,
\scrh_1(X_y, \zz_2))$ can be considered as a subgroup of $\sfh_2(X_{\Gamma_i},
\zz_2)$. Moreover, the image of $\sfh_1( \Gamma_i, \scrh_1(X_y, \zz_2))$ in
$\sfh_2(X_{D_i}, \zz_2)$ is trivial. The lemma follows.
\qed

\medskip%
Let us establish necessary conditions for the vanishing of $[K]$ in $\sfh_2( X, \zz_2)$.

\newprop{prop-K-van} Assume that $\mu = [\mbfm]$ is non-trivial but the
projection of the class $[K]$ to $\sfh_1( Y^\circ, \scrh_1(X_y, \zz_2))$ vanishes.
\begin{itemize}
\item[\sli] In the case of separating $\Gamma$ there exists a decomposition $\mu =
 \mu_+ + \mu_-\in \sfh_1(\Sigma,\zz_2)$ such that for the pieces $Y^\circ_+$ and $Y^\circ_-$
 the monodromy action of the fundamental group $\pi_1(Y^\circ_+)$ (resp.,
 $\pi_1(Y^\circ_-)\,)$ preserves $\mu_+$ (resp., $\mu_-$). In particular, both $\mu_+$ and
 $\mu_-$ are $F_\Gamma$-invariant.
\item[\slii] In the case of non-separating $\Gamma$ the monodromy group of $\pr:
 X^\circ \to Y^\circ$ acts trivially on $\mu \in \sfh_1(\Sigma,\zz_2)$. Moreover, there exists a
 class $\nu \in \sfh_1(\Sigma,\zz_2)$ such that $\mu = (1 + F_{\eta_1}\inv)\nu$ for the
 monodromy $F_{\eta_1}$ along $\eta_1$ and such that $\nu$ is invariant under the
 action of the remaining generators in the above presentation of $\pi_1(Y^\circ)$.
\end{itemize}
\end{prop}

\proof To simplify notation, we denote by $\xi_i, \eta_i, \Gamma_j$ the curves on $Y^\circ$,
the corresponding elements in $\pi_1(Y^\circ)$, the generators of the associated free group,
and the monodromy transformations along these curves.

\smallskip\noindent
\sli \ It follows from  \lemma{lem-h2-circ} that  the element $\Gamma \otimes_1 \mu$ can be
represented in the form
\[\textstyle
\Gamma \otimes_1 \mu = \dd_2 \big( R^+ \cdot \Gamma \otimes_2 \mu_+\big) + \dd_2 \big( R^- \cdot \Gamma\inv \otimes_2 \mu_-\big) + \sum_i \Gamma_i \otimes_1 \lambda_i
\]
for some $\mu_\pm, \lambda_i \in \sfh_1(\Sigma,\zz_2)$ such that $\lambda_i$ are $\cap$-orthogonal to
$\delta_i$. Let us expand $\dd_2 ( R^+ \cdot \Gamma \otimes_2 \mu_+)$ preserving the commutators
$[\xi_i, \eta_i]$ and denoting by $w^\Gamma_i$ the final subword of $R^+ \cdot \Gamma$ after the
letter $\Gamma_i$ and by $w^{[]}_i$ the final subword of $R^+ \cdot \Gamma$ after the
commutator $[\xi_i, \eta_i]$. In particular, $w^\Gamma_{n^+} = \Gamma$, $w^\Gamma_{n^+-1} =
\Gamma_{n^+} \cdot\Gamma$, $w^\Gamma_1 = \Gamma_2\cdot \ldots \cdot \Gamma_{n^+} \cdot\Gamma$, and similarly for $w^{[\,]}_i$.
This gives
\[
\begin{split}
 \dd_2 \big( R^+ \cdot \Gamma \otimes_2 \mu_+\big) &=
 \dd_2 \big( [\xi_1, \eta_1]\cdot \ldots \Gamma_{n^+-1} \cdot \Gamma_{n^+} \cdot \Gamma \otimes_2 \mu_+\big) \\
 &= \Gamma\otimes_1 \mu_+ + \Gamma_{n^+}\otimes_1 \Gamma \mu_+ + \Gamma_{n^+-1}\otimes_1 \Gamma_{n^+} \cdot\Gamma \mu_+ + \cdots+ \Gamma_1 \otimes_1
 w^\Gamma_1\cdot \mu_+ \\
 & \hspace{20pt} + \dd_2 \big([\xi_{g^+_Y}, \eta _{g^+_Y}\big] \otimes_2 w^{[\,]}_{g^+_Y} \mu_+\big)
 +\dd_2 \big([\xi_{g^+_Y-1}, \eta _{g^+_Y-1}\big] \otimes_2 w^{[\,]}_{g^+_Y-1}\mu_+ \big)
+\cdots\\
 & \hspace{20pt} + \dd_2 \big([\xi_1, \eta_1\big] \otimes_2 w^{[\,]}_1\mu_+ \big)
\end{split}
\]
and a similar formula for $\dd_2 \big( R^- \cdot \Gamma\inv \otimes_2 \mu_-\big)$, with the only difference
that the first summand will be $-\Gamma \otimes_2 \Gamma\inv \mu_-$.  The expansion of the
commutators gives
\begin{equation}\eqqno(exp-comm)
\begin{split}
\dd_2 \big( [\xi, \eta] \otimes_2 \nu \big) &= \dd_2 ( \xi\eta\xi\inv \eta\inv \otimes_2 \nu ) \\
 &= \xi \otimes_1 \eta\xi\inv \eta\inv\nu  + \eta\otimes_1 \xi\inv \eta\inv\nu
 - \xi \otimes_1 \xi\inv \eta\inv\nu -\eta \otimes_1\eta\inv\nu \\
&=\xi \otimes_1 ( \eta-1)\xi\inv \eta\inv\nu + \eta\otimes_1 (1-\xi) \xi\inv \eta\inv\nu.
\end{split}
\end{equation}%
Collecting similar summands (and ignoring the signs since we are working with
$\zz_2$-spaces) we obtain the desired decomposition $\mu = \mu_+ + \Gamma\inv \mu_-$ from
the coefficient of $\Gamma$, the equality $\lambda_i = w^\Gamma_i \mu_\pm$ as the coefficient of
$\Gamma_i$, and the equalities $( \eta_i-1)\xi_i\inv \eta_i\inv w^{[\,]}_i \mu_\pm =0$,
$(1-\xi_i)\xi_i\inv \eta_i\inv w^{[\,]}_i \mu_\pm =0$ from the coefficients of $\xi_i$ and~$\eta_i$.
Now observe that the equality $\lambda_i = w^\Gamma_i \mu_\pm$ together with the
$\cap$-orthogonality $\lambda_i \cap \delta_i \equiv 0 \mod2$ implies that $w^\Gamma_i \mu_\pm$ is invariant
under the action of the Dehn twist $T_{\delta_i}$ which is the action of $\Gamma_i$.
This yields the identity $w^\Gamma_0 \mu_+ = w^\Gamma_1 \mu_+ = \cdots =\Gamma \mu_+$ for all $i=0; 1,\ldots,
n^+$. Here we have set $w^\Gamma_0 := \Gamma_1 \cdot w^\Gamma_1$, this element coincides with $w^{[\,]}_{g^+_Y}$.

The triviality of the action of $\xi_i, \eta_i$ is proved in the same way. Indeed, the
equalities above mean that both $\xi_i$ and $\eta_i$ act trivially on $\xi_i\inv
\eta_i\inv w^{[\,]}_i \mu_+$, so that $ w^{[\,]}_i \mu_+$ also remains invariant, and
then we conclude $w^{[\,]}_{i-1} \mu_+ = w^{[\,]}_i \mu_+$. Summing up, we obtain
the $\pi_1(Y^\circ_+)$-invariance of $\Gamma \mu_+$. Finally, since $\Gamma$ can be expressed in
$\pi_1(Y^\circ_+)$ as a product of already treated elements, $\Gamma \mu_+ = \mu_+$. The
argument for the case of $\mu_-$ is the same.

\medskip\noindent
\slii \ This time we have $\Gamma = \xi_1$ in $\pi_1(Y^\circ)$ and the unexpanded relation reads
\[\textstyle
\xi_1 \otimes_1 \mu = \dd_2 \big( R \otimes_2 \nu\big)  + \sum_i \Gamma_i \otimes_1 \lambda_i.
\]
The expansion provides the same equalities as above, except for the coefficient of
$\xi_1$ which now is $\mu + ( \eta_i-1)\xi_i\inv \eta_i\inv w^{[\,]}_1 \nu$. As above, we can
conclude the invariance of $\nu$ with respect to $\xi_i, \eta_i$ for $i=2, \ldots, n$ and
with respect to all $\Gamma_i$'s. So the remaining equalities are $\xi_1\inv \eta_1\inv \nu
= \eta_1\inv \nu$ and $(1 + \eta_1)\xi_1\inv \eta_1\inv \nu =\mu$. Since the commutator $[\xi_1,
\eta_1]$ is equal in $\pi_1(Y^\circ)$ to an expression in the remaining generators, the
actions of $\xi_1$ and $\eta_1$ on $\nu$ commute. Thus the first remaining equality
is equivalent to the $\xi_1$-invariance of $\nu$, and the second yields $(1 +
\eta_1\inv) \nu =\mu$, as desired.
\qed

\smallskip
\newlemma{lem-h1} The first homology group of $X$ can be included into the
exact sequence
\begin{equation}\eqqno(h1-x)
0 \to \sfh_0( Y^\circ , \scrh_1(X_y, \zz_2)) \to \sfh_1(X, \zz_2)
\to \sfh_1(Y, \zz_2) \to 0.
\end{equation}
\end{lem}

\proof The claim is trivial in the special case when $Y =S^2$ and there are no
critical points of $\pr: X\to Y =S^2$. In the remaining cases we consider the
homology of $X^\circ$. The Leray spectral sequence reduces to the exact sequence
\begin{equation}\eqqno(h1-x0)
0 \to \sfh_0( Y^\circ , \scrh_1(X_y, \zz_2)) \to \sfh_1(X^\circ, \zz_2)
\to \sfh_1(Y^\circ, \zz_2) \to 0.
\end{equation}
To obtain \eqqref(h1-x), we use the Mayer-Vietoris sequence corresponding
to the covering $X = X^\circ \cup \big( \sqcup_i X_{D_i} \big)$ with $X^\circ \bigcap \big( \sqcup_i
X_{D_i} \big) = \sqcup_i X_{\Gamma_i}$. This gives us the realisation of $\sfh_1(X,
\zz_2)$ as the cokernel of the homomorphism
\begin{equation}\eqqno(h1-MV-x)\textstyle
\bigoplus_i  \sfh_1(X_{\Gamma_i}, \zz_2)  \to  \sfh_1(X^\circ, \zz_2) \oplus \bigoplus_i  \sfh_1(X_{D_i}, \zz_2)
\end{equation}
For each $X_{\Gamma_i}$ we have
\begin{equation}\eqqno(h1-gga)
0 \to \sfh_0( \Gamma_i , \scrh_1(X_y, \zz_2)) \to \sfh_1(X_{\Gamma_i}, \zz_2)
\to \sfh_1(\Gamma_i, \zz_2) \to 0.
\end{equation}
Now observe that the embedding $X_{\Gamma_i} \subset X_{D_i}$ induces an isomorphism
$\sfh_0( \Gamma_i , \scrh_1(X_y, \zz_2)) \break \xrar{\,\cong\,} \sfh_1(X_{D_i}, \zz_2)$ and
hence a natural splitting of  \eqqref(h1-gga). ``Inserting'' this splitting in
\eqqref(h1-MV-x) we see that each $\sfh_1(\Gamma_i, \zz_2) \cong \zz_2$ ``kills'' the
class in $\sfh_1(Y^\circ, \zz_2)$ from \eqqref(h1-x0) represented by the curve
$\Gamma_i \subset Y$.
\qed

\smallskip%

\newsubsection[z2-sec]{$\zz_2$-sections of Lefschetz fibrations}
If $K$ is $\zz_2$-homologically non-trivial under the assumptions of \propo{prop-K-van},
then it must be $\zz_2$-homologous to the generic fibre $\pr\inv(y_0)$. (The same
situation occurred in \propo{prop-mu-0} above.) In turn, the homological
non-triviality of the generic fibre is equivalent to the existence of a class
$[\sigma] \in \sfh_2(X,  \zz_2)$ with $[\sigma] \cap [K] \neq0$. In this paragraph
we study the properties of such $[\sigma]$. This will give us the possibility to construct
examples of Lagrangian embeddings of the Klein bottle in symplectic manifolds $(X, \omega)$
with trivial homology class.

\medskip%
We maintain the
notation from the previous paragraph with a minor modification. In particular,
$Y$ is a compact surface with a smooth boundary $\dd Y$, possibly empty or not
connected, and $\pr: X\to Y$ a Lefschetz fibration with a generic fibre $\Sigma =
\pr\inv(y_0)$ of genus $g$ such that all critical values $y^*_i$ do not lie on
$\dd Y$.  Denote by $Y^\circ$ the surface obtained from $Y^\circ$ by cutting small
pairwise disjoint discs $D_i$ centred at $y^*_i$. Then the restriction $\pr:
\dd X \to \dd Y$ is a fibre bundle with fibre $\Sigma$.

\newdefi{def-z2-sec} A \slsf{$\zz_2$-section} of\/ $\pr: X\to Y$ is a relative
homology class $\sigma \in \sfh_2(X, \dd X; \zz_2)$ which has non-trivial
$\zz_2$-intersection index with the class $[\Sigma]$ of the generic fibre. We
denote by $\dd \sigma \in \sfh_1(\dd X, \zz_2)$ the class of the boundary of $\sigma$ and by
$\sigma^\lor \in \sfh^2(X, \zz_2)$ the Poincar{\'e} dual class. For an open set $U \subset Y$ with
smooth boundary $\dd U$ which avoids the critical values $y^*_i$ we define the
\slsf{restriction of $\sigma$} on the piece $X_U = \pr\inv(U)$ by means of the
restriction of $\sigma^\lor$, so that $(\sigma|_{X_U})^\lor := \sigma^\lor|_{X_U}$.
\end{defi}
The latter definition can also be used to define $\dd \sigma$, since $(\dd \sigma)^\lor =
\sigma^\lor|_{\dd X}$.

\medskip
The cohomology group $\sfh^2(X, \zz_2)$ where $\sigma^\lor$ ``lives'' admits a
description dual to the one given above for the homology group $\sfh^2(X, \zz_2)$.
In particular, we have the  exact sequence
\[
0 \to \sfh^1( Y^\circ, \scrh^1(X_y, \zz_2)) \to \sfh^2(X^\circ, \zz_2) \to
\sfh^0( Y^\circ, \scrh^2(X_y, \zz_2)) \to 0,
\]
where $\scrh^p(X_y, \zz_2)$ denotes the locally constant sheaf with the stalk
$\sfh^p(X_y, \zz_2)$ at $y \in Y^\circ$, $\sfh^0( Y^\circ, \scrh^2(X_y, \zz_2))$ is
naturally isomorphic to $\zz_2\langle[\Sigma]\rangle$, and $\sfh^1( Y^\circ, \scrh^1(X_y, \zz_2)) =
\sfh^1( \pi_1(Y^\circ), \scrh^1(X_y, \zz_2))$. Moreover, for any $\pi_1(Y^\circ)$-module
$M$ we can calculate the groups $\sfh^1( \pi_1(Y^\circ), M)$ using the
resolvent $0\to \scrr^0 \xrar{\; d^{(1)}\;} \scrr^1 \xrar{\; d^{(2)}\;} \scrr^2
\to 0$, dual to the resolvent $(\scrr_\bullet, \dd_\bullet)$ from the previous subsection. In
the case of a finite presentation $\langle \calx | \calr \rangle$ of $\pi_1(Y^\circ)$ the
cochain groups $\scrr^\bullet$ coincide with the corresponding chain groups, so that
$\scrr^0 =M$, $\scrr^1 =\calx \otimes^1 M$, and $\scrr^2 = \calr \otimes^2 M$, and the
index in $\otimes^k$ indicates that the element lies in $\scrr^k$. The
differentials $d^\bullet$ are dual to $\dd_\bullet$ and are given by the transposed
matrices. One sees immediately the formula $d^{(1)} : m \mapsto \sum_i x_i \otimes^1 (x_i -1)
\cdot m$ for $d^{(1)}: \scrr^0 =M \to \scrr^1$. The formal expression for $d^{(2)}:
\scrr^1 \to \scrr^2$ is $d^{(2)} \big(\sum_i x_i \otimes^1 m_i \big) = \sum_{ij} R_j \otimes^2
\frac{\dd R_j}{\dd x_i} \cdot m_i$ where the sums are taken over all $x_i \in \calx$
and $R_j \in \calr$, and $\frac{\dd R_j}{\dd x_i}$ denotes the \slsf{Fox free
 differential} of the word $R_j$ in the free group $\fr(\calx)$ generated by
$\calx$, see \cite{Bi-2}.  For a given $R_j$, the coefficient $\sum_i \frac{\d
 R_j}{\dd x_i} \cdot m_i$ equals the value on $R_j\in \fr(\calx)$ of the unique
\slsf{cross-homomorphism} $\phi: \fr(\calx) \to M$ with $\phi(x_i) =m_i$. Recall that
cross-homomorphisms are characterised by the property that $\phi(a\cdot b) = \phi(a) + a\cdot \phi(b)$
for every $a,b \in\fr(\calx)$, which is simply the dual of the derivation
properties of the differential $\dd_2$.

\smallskip%

For $g\geq2$, denote by $\wtmap_g$ the group in the extension
\begin{equation}\eqqno(wtmap)
1 \to \sfh_1(\Sigma, \zz_2) \to \wtmap_g \to \map_g\to1,
\end{equation}
so that $\wtmap_g$ is the quotient of $\map_{g,1}$ by the image with respect
to $\dd_\pi$ of the kernel of the homomorphism $\pi_1(\Sigma, z_0)\to \sfh_1(\Sigma, \zz_2)$. A
\slsf{Dehn twist} in $\wtmap_g$ is a projection of a Dehn twist $T_\delta\in
\map_{g,1}$.

\newlemma{lem-h2-circ1} Let $\pr: X\to Y$ be a Lefschetz fibration of a
\emph{closed} manifold $X$ and $\langle \xi_j, \eta_j, \Gamma_i | R \rangle$ the above presentation
of the fundamental group $\pi_1(Y^\circ)$ with a single relation word.

\sli The image of the homomorphism $\sfh^2(X, \zz_2) \to\sfh^2(X^\circ, \zz_2)$ is
generated by the class of some $\zz_2$-section (if such exists) and by
classes from $\sfh^1( Y^\circ, \scrh^1(X_y, \zz_2))$ represented by cocycles
$\bfla^\lor$ of the form
\begin{equation}\eqqno(h2lef)
\textstyle
\bfla^\lor = \sum_j(\xi_j \otimes^1 \lambda_j^\lor + \eta_j \otimes^1 \mu_j^\lor) + \sum_i \Gamma_i \otimes^1 n_i\delta_i^\lor \in \sfh^1( Y^\circ,
\scrh^1(X_y, \zz_2)),
\end{equation}
where $n_i \in \zz_2$, $\lambda_j^\lor, \mu_j^\lor \in \sfh^1(\Sigma, \zz_2)$, and $\delta_i^\lor$ is
the Poincar{\'e} dual of the class $[\delta_i]$.

\slii The obstruction to the existence of a $\zz_2$-section of the fibration $\pr:X\to Y$ is a coset class
$[\vartheta^\lor]$ of the quotient of the group $\scrr^2$ by the subgroup of all coboundaries
$d^{(2)}(\bfla^\lor)$ such that $\bfla^\lor$ has the form \eqqref(h2lef).

\sliii In the case $g\geq2$, the obstruction $[\vartheta^\lor] \in \scrr^2 / \{
d^{(2)}(\bfla^\lor)\}$ vanishes if and only if the monodromy $\calf: \pi_1(Y^\circ) \to
\map_g$ of $\pr: X^\circ \to Y^\circ$ can be lifted to a homomorphism $\wt\calf: \pi_1(Y^\circ)
\to\wtmap_g$ such that $\wt\calf(\Gamma_i)$ is a Dehn twist in $\wtmap_g$.  The
difference between such two lifts $\wt\calf_1, \wt\calf_2: \pi_1(Y^\circ) \to
\wtmap_g$ is a cross-homomorphism $\phi: \pi_1(Y^\circ) \to \sfh_1(\Sigma_1)$ which is
associated with the unique cochain $\bfla^\lor_\phi$ of the form \eqqref(h2lef)
whose Poincar{\'e} dual $\bfla_\phi$ has coefficients $\lambda_j= \phi(\xi_j)$, $\mu_j= \phi(\eta_j)$,
and $n_i \delta_i = \phi(\Gamma_i)$.
\end{lem}

\proof The first two assertions are dual to \lemma{lem-h2-circ}, so the proof
proceeds by dualising the homomorphism and taking into account the duality
between kernels and images. The last assertion follows from the fact that the
natural embedding $\dd_\pi: \pi_1(\Sigma, z_0) \hrar \map(\Sigma, z_0)$
induces an inclusion $\sfh_1(\Sigma, \zz_2) \hrar \wtmap_g$ for any $g\geq2$.
Then one observes that the group-homological meaning of both constructions is the same.
\qed

\medskip%
It should be noticed that in the case $g=1$ the obstruction $[\vartheta^\lor]$ to
the existence of a $\zz_2$-section is \emph{not determined} by the monodromy in
$\map_1 = \map_{1,1}$. This reflects the fact that the homomorphism $\dd_\pi:
\pi_1(T^2) = \zz^2 \to \pi_0(\diff_+(T^2))$ in \eqqref(hmtpy) is trivial.

\newsubsection[lagr-triv]{An example.} We construct an example of a Lagrangian
embedding of the Klein bottle in a projective surface $X$ equipped with a
K{\"a}hler form $\omega$ such that the homology class $[K]$ is trivial.
In the example, $K$ will be fibred over a non-separating circle in the two-torus $T^2$.

Let $Y$ be a two-torus. Pick a flat metric on $Y$. This gives us a
complex structure and a K{\"a}hler form $\omega_Y$ on $Y$. We may assume that the
$\omega_Y$-volume of $Y$ is $1$. Fix a geometric basis of $Y$, represented by
embedded curves $\xi$ and $\eta$ meeting transversally in a single point.

Let $\Sigma$ be another two-torus, realised as the quotient of $\cc$ by the
lattice $\Lambda$ generated by vectors $\alpha := 1$ and $\beta := e^{2\pi\isl /3}$.
Let $\omega_\Sigma$ be the flat K{\"a}hler form on $\Sigma$ such that the
$\omega_\Sigma$-volume of $\Sigma$ is also $1$. Consider the $\cc$-linear homomorphisms
$F_\xi := -\id: \cc \to \cc$ and $F_\eta := e^{\pi\isl /3}\id: \cc \to \cc$. Clearly, they
define holomorphic automorphisms of $\Sigma$ preserving the base point $z_0 := 0 \in \Sigma$,
also denoted by $F_\xi$ and $F_\eta$.

Let $\pr: X \to Y$ be the fibration over the base $Y$ with the fibre $\Sigma$ and with
monodromy $F_\xi$ along $\xi$ and $F_\eta$ along $\eta$. It follows that there exists
a flat K{\"a}hler structure on $X$ with the K{\"a}hler form $\omega$, in which $X$ is
the product $(Y, \omega_Y) \times (\Sigma, \omega_\Sigma)$ locally near each fibre. In particular,
there exists a fibrewise K{\"a}hler form $\omega_\Sigma$ on $X$ such that $\omega = \pr^*\omega_Y +
\omega_\Sigma$.

We claim that $\omega$ is a polarisation of $X$ corresponding to a certain
holomorphic line bundle $\scrl$. Indeed, since the monodromy preserves the
base point $z_0 \in \Sigma$, we obtain the horizontal section $\sigma_0$ which is constantly
$z_0$. It follows that $\sigma_0$ is holomorphic. Now it is easy to see that $\omega$
represents the Chern class $c_1(\scrl)$ of the holomorphic line bundle
$\scrl:= \scro_X(\sigma_0 + X_y)$ of the divisor $\sigma_0+  X_y$, $X_y$ being the vertical fibre
over some point $y \in Y$. By Kodaira's embedding theorem, see \eg,
\cite{Gr-Ha}, $\scrl$ is ample and hence $X$ is projective algebraic.

Now let us compute $\sfh_2(X, \zz_2)$ using the homology spectral sequence.  Its
$E^2_{p,q }$-term consists of the groups $\sfh_p(Y, \scrh_q(X_y, \zz_2))$. The
groups $\sfh_0(Y, \scrh_2(X_y, \zz_2))$ and $\sfh_2(Y, \scrh_0(X_y, \zz_2))$
are both $\zz_2$'s, generated by the class of the fibre $X_y$ and the section $\sigma_0$, respectively.
To compute $\sfh_1(Y,\scrh_1(X_y, \zz_2))$, observe first that the monodromy actions of $F_\xi$ and
$F_\eta$ are given by the $\zz_2$-matrices $\id$ and $\spmatr{1&1\\1&0}$,
respectively. Using the notation of \refsubsection{hmlg-TLF}, we see that
$\sfh_1(Y, \scrh_1(X_y, \zz_2))$ is represented by  the cycles $\bfla = \xi \otimes \lambda_\xi + \eta
\otimes \lambda_\eta$ with $\lambda_\xi,\lambda_\eta \in \sfh_1(\Sigma, \zz_2)$ satisfying $\dd_1(\bfla) = (F_\eta
-\id)\lambda_\eta = 0$. Since $F_\eta -\id$ is non-degenerate, $\bfla$ must be of the form
$\xi \otimes \lambda_\xi$. The calculation from \eqqref(exp-comm) shows that the
boundaries can be written as
\[
\dd_2([\xi, \eta] \otimes \nu) = \xi \otimes (F_\eta -\id) F_\xi\inv F_\eta\inv \nu + \eta\otimes (\id- F_\xi
)  F_\xi\inv F_\eta\inv \nu  = \xi \otimes (F_\eta -\id) F_\xi\inv F_\eta\inv \nu
\]
and hence cover the entire group of chains. Thus $\sfh_1(Y, \scrh_1(X_y,
\zz_2))$ is trivial. Clearly, the classes of the fibre and the section survive
in $\sfh_2(X, \zz_2)$ (this means that the spectral sequence degenerates in
the term $E^2_{p,q }$), so that we obtain $\sfh_2(X, \zz_2) = \zz_2\langle [X_y],
[\sigma_0] \rangle$.

Now let us construct two Lagrangian embeddings of the Klein bottle in
$X$. Consider the following curves on $\Sigma$:
\[\textstyle
\beta_0 := \big\{ t\cdot e^{2\pi\isl /3} \in \Sigma = \cc/\Lambda : t \in \rr \big\}
\qquad
\beta_1 := \big\{ \frac12\cdot e^{\pi\isl /3} + t\cdot e^{2\pi\isl /3} \in \Sigma = \cc/\Lambda : t \in \rr \big\}
\]
Both curves are closed geodesics on $\Sigma$ representing the homology class
corresponding to $\beta= e^{2\pi\isl /3} \in \Lambda \cong H_1(\Sigma, \zz)$. Moreover, both curves
are invariant with respect to the transformation $F_\xi = -\id$. For the curve $\beta_1$,
this follows from fact that $2  \cdot \frac12e^{\pi\isl /3} =
e^{\pi\isl /3} = e^{2\pi\isl /3} + 1$ lies in the lattice $\Lambda$.

Now define the surfaces $K_0$ and $K_1$ in $X$ which consist of points $(y,z)
\in X$ such that $y$ lies on $\xi$ and $z$ on $\beta_0$ or $\beta_1$, respectively.
In other words, we take a point $y_0 \in\xi$, realise $\beta_0$ and $\beta_1$ on the fibre $\pr\inv(y_0)$,
and then carry them around $\xi$ using parallel transport. Since $F_\xi(\beta_i) = \beta_i$, $K_0$ and $K_1$
are indeed closed surfaces. From the fact that $F_\xi$ inverts the orientation
on each $\beta_i$ we conclude that $K_i$ are Klein bottles. The local product structure
of the symplectic form $\omega$ ensures that $K_i$ are $\omega$-Lagrangian.

To find the homology classes $[K_i] \in \sfh_2(X, \zz_2)$ we observe the
following. Since $K_i$ are disjoint from the generic fibre of $X$, each $K_i$
must be $\zz_2$-homologous to the fibre or homologically trivial. Since $K_1$ is
disjoint from $\sigma_0$, $[K_1]=0 \in \sfh_2(X, \zz_2)$. On the other hand, the pair
$\beta_0, \beta_1$ separates the torus $\Sigma$, and applying the argument from the
proof of \propo{prop-mu-0} to the union $\beta_0 \cup \beta_1$ instead of the meridian
circle $\mbfm$ we conclude that $[K_0] + [K_1]$ is homologous to the fibre. Thus
$[K_0] = [X_y] \neq 0 \in \sfh_2(X, \zz_2)$.

\newsection[map-mu]{Combinatorial structure of mapping class groups}

\newsubsection[graph-braid]{Coxeter--Weyl and braid groups.} Let us first recall
standard definitions and facts concerning Coxeter--Weyl groups and Artin--Brieskorn braid
groups.

Let $\cals =\{s_1, \ldots, s_r\}$ be any finite set. A \slsf{Coxeter matrix} over
$\cals$ is a symmetric $r\times r$-matrix $M = (m_{ij})$ with entries in $\nn \cup
\{\infty\}$, such that $m_{ii} =1$ and $m_{ij} \geq2$ for $i\neq j$. The pair $(\cals, M)$
is called the \slsf{Coxeter system}. Often we shorten the notation for a
Coxeter system to $\cals$ and understand $M$ as a structure on $\cals$,
denoting it by $M_\cals$. If $s=s_i, s'=s_j \in \cals$, then $m_{ij}$ is also denoted by $m_{s,s'}$.

The \slsf{Coxeter graph} $\Delta = \Delta(\cals, M)$ associated with a Coxeter system
$(\cals, M)$ has $\cals$ as the set of vertices; $s_i$ and $s_j$ are connected
by a wedge iff $m_{ij} \geq3$; if $m_{ij} \geq4$, this wedge is labelled by the entry $m_{ij}=m_{ji}$.
The labelled Coxeter graph $\Delta$ completely determines the associated Coxeter matrix.

For any subset $\cals' \subset \cals$ the induced Coxeter matrix $M'$ is the
restriction of $M= M_\cals$ onto $\cals'$. The associated graph $\Delta':= \Delta(\cals',
M')$ is a full subgraph of $\Delta$, \ie, any two vertices $s_1, s_2 \in \cals'$ are
connected and labelled in $\Delta'$ in the same way as in $\Delta$.

For any two letters $a,b$ and a non-negative integer $m$ we denote by $\langle ab\rangle^m$
the word $abab\ldots$ of the length $m$ consisting of alternating letters $a$ and
$b$. If $a$ and $b$ lie in a group $G$, then $\langle ab\rangle^m$ denotes the
corresponding product in $G$.

\newdefi{def-gra-group}  The \slsf{Artin--Brieskorn braid group}
$\br(\cals)= \br(\cals, M)$ of a Coxeter system $(\cals, M)$ is the
group generated by $\cals$ with \slsf{(generalised) braid relations}
\[
\langle s_is_j\rangle^{m_{ij}} = \langle s_js_i\rangle^{m_{ij}}
\qquad
\text{for each $i\neq j$, such that $m_{ij}<\infty$}
\]
as defining relations. The \slsf{Coxeter--Weyl group} $\sfw(\cals)= \sfw(\cals,
M)$ is obtained from the braid group $\br(\cals)$ by adding the
\slsf{reflection relation} $s^2=1$ for each $s \in \cals$.

The kernel of the natural projection $\pi: \br(\cals) \to \sfw(\cals)$ is called
the \slsf{pure braid group} of the Coxeter system $(\cals, M)$ and is denoted
by $\sfp(\cals)$. Its \slsf{abelianisation} is denoted by $\sfp_\ab(\cals)$.
\end{defi}

To distinguish equalities in different groups, we use the notation $w\equiv w'$
for equality in $\sfw(\cals)$ and $w=w'$ for equality in $\br(\cals)$.
An \slsf{expression} of an element $w$ from $\sfw(\cals)$ or from
$\br(\cals)$ is a word $\mbfs = s_1^{\epsilon_1} \cdot \ldots \cdot s_q^{\epsilon_q}$ in the
alphabet $\cals$, $\epsilon_i=\pm1$, whose evaluation in $\sfw(\cals)$ or,
respectively, in $\br(\cals)$ yields $w$.  Such an expression $\mbfs$ is
called \slsf{positive} if all $\epsilon_i$ equal $+1$. In view of the
reflection relation $s_i\inv = s_i\vph$ we may consider only positive
expressions of $w\in \sfw(\cals)$.  An element $w\in \br(\cals)$ is called
positive if it admits a positive expression.

\newdefi{def-cox1}
The \slsf{length} $\ell(w)$ of $w\in\sfw(\cals)$ is the minimal possible
length of an expression $\mbfs = s_1 \cdot \ldots \cdot s_l$ of $w$. The
\slsf{$\sfw$-length} of $x \in \br(\cals)$ is the length $\ell_\sfw(\bar x)$ of its
projection $\bar x := \pi(x)$ to $\sfw(\cals)$.  An expression $\mbfs$ is
\slsf{reduced} (or \slsf{$\sfw$-reduced}) if it is positive and its length
equals its $\sfw$-length.

A \slsf{quasireflection} in a Coxeter--Weyl group $\sfw(\cals)$ is any element
$t$ conjugated to a generator $s\in \cals$. We denote by $\calt(\cals)$ the set
of quasireflection in $\sfw(\cals)$.
\end{defi}

Note that quasireflection are usual reflections in \emph{finite}
Coxeter--Weyl groups, and transpositions in the symmetric groups $\sym_n$.
We will use the following standard properties of Coxeter--Weyl groups,
see \eg\ \cite{Bou} or~\cite{Hum}.

\newprop{prop-cox-gr} \sli \slsf{(Strong Exchange Property.)} Let $\mbfs =
s_1\cdots s_q \in \sfw(\cals)$ be a (positive) expression of an element $w \in \sfw(\cals)$, $s\in
\cals$ a generator, and $t \in \calt(\cals, M)$ a quasireflection. Then $\ell(ws) =
\ell(w) \pm1$ and $\ell(wt) -\ell(w)$ is odd. Moreover, if $\ell(wt) < \ell(w)$, then there
exists $j =1,\ldots,q$ such that $w\equiv s_1\cdots s_{j-1} \, s_{j+1} \cdots s_q \, t$.
The analogous property holds for the left products $s w$ and $t w$.

\slii \slsf{(Relations between generators.)}  The order of $s_is_j$ in
$\sfw(\cals)$ is $m_{ij}$. Furthermore, for any subset $\cals' \subset \cals$, the
natural homomorphism $\sfw(\cals') \to \sfw(\cals)$ is injective, \ie, there are
no new relations in $\sfw(\cals)$ among the generators of $\sfw(\cals')$.

\sliii \slsf{(Uniqueness of the set of factors.)}  For any $w \in \sfw(\cals)$
and any \emph{reduced} expression $w= s_1\cdots s_l$ the set of factors $\cals_w :=
\{s_1, \ldots, s_l\}$ is independent of the particular choice of the reduced expression;
the set $\cals_w$ generates the unique minimal subsystem $\cals' \subset \cals$ such
that $w$ lies in $\sfw(\cals')$.

\sliv \slsf{(Conjugacy classes of generators.)}  Two generators $s', s'' \in
\cals$ are conjugated if and only if they can be connected by a chain $s_0
=s',s_1, \ldots, s_{n-1}, s_n=s''$ such that each index $m_{s_{i-1},s_i}$ is odd.
In particular, for a connected graph $\Delta(\cals)$ with no labelling (i.\,e.\ such that
each $m_{ij}$ is $2$ or $3$), the set of quasireflections $\calt(\cals)$
consists of elements conjugated to any given $s \in \cals$.
\end{prop}

\newlemma{lem-subword} Suppose that $w\in \sfw(\cals, M)$ admits two factorisations $w \equiv
s_1w_1 \equiv s_2w_2$ with $s_1\neq s_2 \in \cals$ and $w_1,w_2 \in \sfw(\cals)$ such that
$\ell(w_1) = \ell(w_2) = \ell(w)-1$. Set $m := m_{s_1s_2}$. Then $w$ admits a
factorisation $w\equiv \langle s_1s_2\rangle^m w_0 \equiv \langle s_2s_1\rangle^m w_0 $ with $\ell(w_0) = \ell(w)-m$.

A similar property holds for factorisation $w \equiv w'_1s_1 \equiv w'_2s_2$.
\end{lem}

\proof Fix factorisations $w \equiv \langle s_1s_2\rangle^{n_1} v_1$ and  $w \equiv \langle s_2s_1\rangle^{n_2}
v_2$ such that $\ell(v_i) = \ell(w) -n_i$ and such that the lengths $n_i$ are the
maximal possible. Then $1\leq n_i\leq m$. Indeed, if $n_1>m$, then $w$ admits a
reduced expression which starts with $\langle s_1s_2\rangle^{m+1}$. But then
\[
\langle s_1s_2\rangle^{m+1} = s_1 \langle s_2s_1\rangle^m \equiv s_1 \langle s_1s_2\rangle^m = s_1^2
\langle s_2s_1\rangle^{m-1} \equiv \langle s_2s_1\rangle^{m-1}
\]
has smaller $\sfw$-length, in contradiction with $\ell(v_1) = \ell(w) -n_1$.

There is nothing to prove if  $n_1$ or $n_2$ equals $m$, so we assume that $n_1,n_2<m$.
Fix reduced expressions $\mbfv_i$ for $v_i$, so that $\langle s_1s_2\rangle^{n_1} \mbfv_1$
and $\langle s_2s_1\rangle^{n_2} \mbfv_2$ are reduced expressions for $w$. Then by the
Exchange Property applied to the relation $\ell(s_2w) =\ell(w)-1$, $w$ admits an
expression $s_2 \mbfw_2$ such that $\mbfw_2$ is obtained by removing a single
generator $s^*$ from $\langle s_1s_2\rangle^{n_1} \mbfv_1$. It could not be the first
letter $s_1$ since $s_1 \neq s_2$. It could not be the last letter in the subword
$\langle s_1s_2\rangle^{n_1}$ since then $\langle s_1s_2\rangle^{n_1} \equiv \langle s_2s_1\rangle^{n_1}$ in
contradiction with $n_1 < m= m_{s_1s_2}$. It could not also be an inner letter
in the subword $\langle s_1s_2\rangle^{n_1}$ since then we would obtain a squared
generator as a subword of $s_2 \mbfw_2$ and conclude that $\ell(s_2 \mbfw_2) < \ell(w)$. Thus
the letter $s^*$ is removed from $\mbfv_1$. Denoting by $\mbfv'_2$ the obtained
word, we get a \emph{reduced} expression $\langle s_2s_1\rangle^{n_1+1} \mbfv'_2$. This
implies that $n_2 > n_1$. By symmetry, we must have $n_1> n_2$. This contradiction shows
that in fact $n_1 = n_2 =m$, as desired.
\qed

\newlemma{lem-w-red} For any two reduced expressions $\mbfs = s_1\cdot \ldots \cdot s_l$
and $\mbfs' = s'_1\cdot \ldots \cdot s'_l$ of a given $w\in \sfw(\cals)$, their evaluations in
$\br(\cals)$ define the same element $\hat w \in \br(\cals)$. In particular, the
correspondence $w \in \sfw(\cals) \mapsto \hat w \in \br(\cals)$ is a set-theoretic section
of the projection $\pi: \br(\cals) \to \sfw(\cals)$.
\end{lem}

\proof Using induction on the length of $w$ we may assume that the claim holds
for all elements $v \in \sfw(\cals)$ of length $\ell(v) < \ell(w)$. This implies the
lemma in the case $s_1 =s'_1$, so we may assume that $s_1 \neq s'_1$. Set $m:=
m_{s_1,s'_1}$. By the previous lemma, $w$ admits reduced factorisations
$\langle s_1s'_1\rangle^m w_0 \equiv \langle s'_1s_1\rangle^m w_0$. By induction, the lift $\hat w_0 \in
\br(\cals)$ is well-defined. Since $\langle s_1s'_1\rangle^m = \langle s'_1s_1\rangle^m$ in
$\br(\cals)$, the latter two factorisations define the same element
$\langle s_1s'_1\rangle^m \hat w_0 = \langle s'_1s_1\rangle^m \hat w_0 \in \br(\cals)$. On the other hand,
$\hat s := s_1\cdot \ldots \cdot s_l$ equals $\langle s_1s'_1\rangle^m \hat w_0$ in $\br(\cals)$ since
they have the same first letter $s_1$. The same argument for $\hat s' := s'_1\cdot
\ldots \cdot s'_l$ finishes the proof.
\qed

\medskip
Applying the Reidemeister--Schreier theorem (see \eg\ \cite{Co-Zi} or
\cite{Ly-Sch}) we obtain a presentation of $\sfp(\cals)$.

\newprop{pu-br-pres} \sli The pure braid group $\sfp(\cals)$ is generated by
the elements $(s_1\ldots s_l)\cdot s_0^2 \cdot(s_1\ldots s_l)\inv$ such that the product $s_1\ldots s_ls_0$
is $\sfw$-reduced.

\slii The Reidemeister--Schreier defining set of relations in $\sfp(\cals)$
consists of products of the form $\wh w \cdot R \cdot \wh w \inv$, where $w$ varies over the
Coxeter--Weyl group $\sfw(\cals)$, $R$ over the defining set of braid relations
in $\br(\cals)$, and  $\wh w \cdot R \cdot \wh w \inv$ is written as a product
of the above generators $(s_1\ldots)\cdot s_0^2 \cdot(s_1\ldots)\inv$.
\end{prop}

\newdefi{def-rs-gen} The elements produced by the proposition are called
\slsf{Reidemeister--Schreier} generators and relations, or simply
\slsf{RS-elements}. If we denote by $a \star b$ the conjugation of $b$ by $a$
so that $a \star b := a\,b\,a\inv$, then all RS-generators have the form $\hat w \star s_0^2$.
\end{defi}

\remark More precisely, the generators produced by the Reidemeister--Schreier
theorem are elements of the free group $\fr(\cals)$, and their explicit form
depends on the concrete expression $s_1\ldots s_l$. However, by \lemma{lem-w-red}, we
only need to know the corresponding element in $\sfw(\cals)$.

\proof We only indicate the algorithm which realises any $f \in \sfp(\cals)$ as
a product of generators. Let $f$ be represented in the form $f= f_1 \cdot \hat w \cdot
s^\epsilon \cdot f_2$ so that $f_1$ is a product of RS-generators, $\hat w$ is a
$\sfw$-reduced element in $\br(\cals)$, $s^\epsilon = s^{\pm1}$ a letter with $s \in
\cals$, and $f_2$ any element in $\br(\cals)$. Find $w_1 \in \sfw(\cals)$
such that $w_1\equiv w \cdot s$. Then $f= f_1 \cdot (\hat w \cdot s^\epsilon \cdot \hat w_1\inv) \cdot \hat
w_1 \cdot f_2$.  Thus we only need to represent $\hat w \cdot s^\epsilon \cdot \hat w_1\inv$ as a
product of RS-generators. Here we must consider four cases according to the possible values of the
$\sfw$-length of $w_1$ and the sign $\epsilon=\pm1$. Recall that $\ell_\sfw(w_1) =
\ell_\sfw(w)\pm1$. Then we obtain
\begin{itemize}
\item case $\ell(w_1) = \ell(w)+1$ and $\epsilon=+1$: \ $\hat w_1 = \hat w \cdot s$ and $ \hat w \cdot s
 \cdot \hat w_1\inv =1$;
\item case $\ell(w_1) = \ell(w)+1$ and $\epsilon=-1$: \ $\hat w_1 = \hat w \cdot s$ and $ \hat
 w \cdot s \inv \cdot \hat w_1\inv = (\hat w \star s^2) \inv$;
\item case $\ell(w_1) = \ell(w) - 1$ and $\epsilon=+1$: \ $\hat w = \hat w_1 \cdot s$ and $
 \hat w \cdot s \cdot \hat w_1\inv = \hat w_1 \star s^2$;
\item case $\ell(w_1) = \ell(w) - 1$ and $\epsilon=-1$: \ $\hat w = \hat w_1 \cdot s$ and $
 \hat w \cdot s\inv \cdot \hat w_1\inv =1$.
\qed

\end{itemize}

\newprop{prop-ab-gen} \sli The abelianisation $\sfp_\ab(\cals)$ of the pure braid
group $\sfp(\cals)$ is a free abelian group admitting a basis generated by
RS-elements $\hat w * s^2 \in \sfp(\cals)$.

\slii The projection of every element $x* s^2$ with $x\in \br(\cals)$ and $s\in
\cals$ into $\sfp_\ab(\cals)$ is a generator $[\hat w * \ti s^2]_\ab$ for some
$w\in \sfw(\cals)$ and $\ti s \in \cals$, not unique in general. Moreover, $\hat w
* \ti s$ can be obtained by a sequence of the following transformations:
\begin{itemize}
\item[(A0)] replacing $x* s^2$ with $x\in \br(\cals)$ by $x'* s^2$ with
 some $x'\in \br(\cals)$ such that $x$ and $x'$ have equal projections to
 $\sfw(\cals)$, $x\equiv x' \in \sfw(\cals)$;
\item[(A1)] replacing $(x \cdot \langle s_1s_2\rangle^k) * s_3^2$ by $(x \cdot \langle s_2s_1\rangle^{m-k-1})
 * s_3'{}^2$ where
\begin{itemize}
\item $x\in\br(\cals)$, $s_1,s_2 \in \cals$, $m:= m_{s_1s_2}<\infty$, $k=0,\ldots,m-1$,
\item $s_3= s_1$ if $k$ is even and $s_3 = s_2$ if $k$ is odd,
\item $s_3' = s_3$ if $m$ is even and $s_3'$ is the remaining element in the
 pair $s_1,s_2$ otherwise.
\end{itemize}
\end{itemize}

\end{prop}

\proof First let us show that relations (A0) and (A1) hold in
$\sfp_\ab(\cals)$. The first one follows from the equality
\[
(xs_1^2y) * s_2^2 = (xs_1^2y)s_2^2(xs_1^2y)\inv =
xs_1^2x\inv \cdot (xy)s_2^2(xy)\inv  \cdot (xs_1^2x\inv)\inv.
\]
For the second one, we observe that it is sufficient to consider the case
$x=1$. Represent the braid relation $\langle s_1s_2\rangle^m = \langle s_2s_1\rangle^m$ in the form
$\langle s_1s_2\rangle^{m-1}s_3 = s_2 \langle s_1s_2\rangle^{m-1}$. An algebraic manipulation gives
\[
\langle s_1s_2\rangle^{m-1}s_3  (\langle s_1s_2\rangle^{m-1})\inv = s_2,
\]
so that squaring yields the desired relation for $k=m-1$. The conjugation by
$(\langle s_1s_2\rangle^l)\inv$ and application of (A0) gives the rest. In the special case
$s_1=s_2 =:s$ we must have $m=1$ and $k=0$, so that $m-k-1=0$ and the only
relation in (A1) is the trivial equality $(x\cdot s)* s^2 = x*s^2$.

Fix $a \neq b \in \cals$ such that $m:= m_{ab} \neq \infty$. Let $w_0 \in \sfw(\cals)$
be any element. Set $a_{2i-1} := b_{2i} := a$, $a_{2i}:= b_{2i-1} := b$, $w_i
:= w_i \cdot \langle ab\rangle^i$, and $w'_i := w_i \cdot \langle ba\rangle^i$.  Observe that $\langle ab\rangle^{2m}\equiv1 \in
\sfw(\cals)$, so that $w_i$ and $w'_i$ are $2m$-periodic, that is, $w_{i+2m} \equiv w_i$ and
$w'_{i+2m} \equiv w'_i$. Besides, $b_i = a_{2m-i}$ and $w'_i \equiv w_{2m-i}$. The
induced RS-relation reads
\[\textstyle
\prod_{i=1}^m (\hat w_{i-1} \cdot a_i \cdot \hat w_i\inv) =
\prod_{i=1}^m (\hat w'_{i-1} \cdot b_i\vph \cdot \hat w'_i{}\inv).
\]

Let $w_k$ be an element of maximal $\sfw$-length. Then $\ell(w_{k-1}) =
\ell(w_{k+1}) = \ell(w_k) -1$. In this situation \lemma{lem-subword} implies that
$v:= w_{k\pm m}$ is the shortest element with $\ell(v) = \ell(w_k) -m$ and for every
$i=0,\ldots,2m$ one has either $w_{k\pm m \pm i} \equiv v \langle ab\rangle^i$ or $w_{k\pm m \pm i} \equiv v \langle
ba\rangle^i$.  It follows that $w_i \equiv v\cdot v_i$ with $\ell(w_i) = \ell(v) + \ell(v_i)$ such
that $v_0$ is either $\langle ab\rangle^l$ or $\langle ba\rangle^l$ for some $l=0,\ldots,m$. This means
that the relation we consider is obtained by conjugation by $v$ from the
relation
\begin{equation}\eqqno(rel-v)\textstyle
\prod_{i=1}^m (\hat v_{i-1} \cdot a_i \cdot \hat v_i\inv) =
\prod_{i=1}^m (\hat v'_{i-1} \cdot b_i\vph \cdot \hat v'_i{}\inv),
\end{equation}
where $v_i$, $v'_i$ are defined in the same way as $w_i$, $w'_i$.

In view of the symmetry $a\leftrightarrow b$, we assume that $v_0 = \langle ab\rangle^l$ with $l\leq m$.
If $l$ is even, we obtain the following sequences of words $v_i$ and $v'_i$:
\[
v_i \equiv\begin{cases}
 \langle ab\rangle^{l+i} & i=0,\ldots,m-l\\
\langle ab\rangle^m \equiv\langle ba\rangle^m & i=m-l\\
\langle ba\rangle^{m+l-i} & i= m-l,\ldots, m
\end{cases}
\qquad
v'_i \equiv\begin{cases}
 \langle ab\rangle^{l-i} & i=0,\ldots,l\\
\langle ab\rangle^0 \equiv\langle ba\rangle^0 \equiv1 & i=m-l\\
\langle ba\rangle^{i-l} & i= l,\ldots, m
\end{cases}
\]
If $l$ is odd, these two sequences are interchanged. The element $\hat v_{i-1}
\cdot a_i \cdot \hat v_i\inv$ is non-trivial if and only if $\ell(v_{i-1}) > \ell(v_i)$,
and the same is true for $\hat v'_{i-1} \cdot b_i \cdot \hat v'_i{}\inv$. The corresponding values of $i$
are $i=m-l+1,\ldots,m$ for $v_i$ and $i=1,\ldots,l$ for $v'_i$. Comparing the
corresponding factors, we see that the induced relation in $\sfp_\ab(\cals)$
is the sum of relations (A1) over $k=1,\ldots,l$.
\qed

\newthm{thm-trans} Two quasireflections $t = w* s$ and $t' = w' * s'$ in
$\sfw(\cals)$ with $w,w' \in \sfw(\cals)$ and $s,s'\in \cals$ are equal if and
only if the second can be obtained from the first by a sequence of the following
transformations:
\begin{itemize}
\item[(T)] replacing $(x \cdot \langle s_1s_2\rangle^k) * s_3$ by $(x \cdot \langle s_2s_1\rangle^{m-k-1})
 * s_3'$ where $x\in \sfw(\cals)$ and $s_1,s_2, s_3,s_3'$ have the same meaning
 as in part {\rm(A1)} of \propo{prop-ab-gen}.
\end{itemize}
In particular, the correspondence
\[
t \equiv w* s \in \calt(\cals) \ \leftrightarrow \ A_t := \hat w * s^2 \in \sfp_\ab(\cals)
\qquad\text{for $w\in \sfw(\cals)$ and $s\in \cals$}
\]
defines a bijection between the set of quasireflections in $\sfw(\cals)$ and the
RS-basis of\/ $\sfp_\ab(\cals)$.
\end{thm}

\proof Denote by $\cala(\cals)$ the RS-basis of $\sfp_\ab(\Delta)$ constructed above and
by $[\hat w*s^2]_\ab $ the image of $\hat w*s^2 = \hat w \cdot s^2 \cdot \hat w\inv$
in $\cala(\cals)$. Observe that the transformation (A1) from
\propo{prop-ab-gen} is obtained by squaring the transformation (T) and that
the analog of (A0) in $\calt(\cals)$ is the trivial transformation. This
implies the second claim of the theorem modulo the first claim.

In any case, we obtain a well-defined surjective map $\mbft: \cala(\cals) \to
\calt(\cals)$ given by $\mbft: [\hat w*s^2]_\ab \mapsto w*s$, and the first claim of
the theorem is equivalent to the injectivity of $\mbft$. Since $\mbft$ is
$\sfw(\cals)$-invariant it is sufficient to show that for every $s\in \cals$ and
every $w\in \sfw(\cals)$ the condition $w*s \equiv s$ implies the equality $[\hat
w*s^2]_\ab = [s^2]_\ab$. We proceed by induction on the length $\ell(w)$. The
case $\ell(w)=0$ is trivial. Fix a pair $(w,s)$ with the property $w*s \equiv s$.
Assume that $\ell(ws) = \ell(w)-1$. Then by the Exchange Property $w=w's$ with
$\ell(w')= \ell(w)-1$. But then
\[
w*s \equiv (w's)*s \equiv (w's)\cdot s\cdot (sw'{}\inv) \equiv w'\cdot s\cdot w'{}\inv \equiv
w'*s,
\]
and hence $[\hat w*s^2]_\ab = [\hat w'*s^2]_\ab = [s^2]_\ab$ by the inductive
assumption. In the case $\ell(ws) = \ell(w)+1$ we conclude from \lemma{lem-w-red}
the equality $\wh{ws} = \hat w\, s$ in $\br(\cals)$. On the other hand,
the equality $w*s \equiv s$ means that $w$ and $s$ commute in $\sfw(\cals)$, and
hence $\ell(sw) = \ell(ws)= \ell(w) +1$. As above, we obtain $\wh{sw} = s\,\hat w$ and
hence $\hat w\, s = s\,\hat w$. This implies the desired identity $[\hat w*s^2]_\ab =
[s^2]_\ab$. The proof follows.
\qed

\medskip
Summing up the results obtained so far, we have

\newthm{thm-gen-pure} For any quasireflection $t \in \calt(\cals)$
and any two representing expressions $t\equiv w*s$ and $t\equiv w'*s'$ with $w,
w' \in \sfw(\cals)$ and $s,s' \in \cals$, the evaluations $\hat w *s$ and
$\hat w' *s'$ define the same element $\ti t \in \br(\cals)$.

The set $\wt\calt^2 := \{ \ti t^2 : t\in \calt(\cals) \}$ is a minimal
generator set for the pure braid group $\sfp(\cals)$.
\end{thm}

\proof Let $s_1,s_2 \in \cals$ be generators and $k< m_{s_1s_2}$ a
non-negative integer. Then the relation $\langle s_1s_2\rangle^k * s_3 = \langle
s_2s_1\rangle^{m-k-1} * s_3'$ with $s_3,s_3'$ as in the transformation (T) holds
in $\br(\cals)$. Moreover, even if $s_1\langle s_1s_2\rangle^k$ is not
$\sfw$-reduced, we obtain the equality
\[
(s_1\langle s_1s_2\rangle^k) * s_3 \buildrel{\text{(T)}}\over{=}
(s_1\langle s_2s_1\rangle^{m-k-1}) * s_3' = \langle s_1s_2\rangle^{m-k}* s_3'
\buildrel{\text{(T)}}\over{=} \langle s_2s_1\rangle^{k-1} * s_3
\]
in $\br(\cals)$. Similarly,
\[
(s_2\langle s_1s_2\rangle^k) * s_3 \buildrel{\text{(T)}}\over{=}
(s_2\langle s_2s_1\rangle^{m-k-1}) * s_3' = \langle s_1s_2\rangle^{m-k-2} * s_3'
\]
for $k>0$. This implies that the relation
$$
(x \cdot \langle s_1s_2\rangle^k)\wh{\ } * s_3 = (x \cdot \langle s_2s_1\rangle^{m-k-1})\wh{\ }
 * s_3',
\leqno{(\wh{\text{T}})}
$$
where $x\in \sfw(\cals)$ and $s_1,s_2, s_3,s_3'$ have the same
meaning as in part {\rm(A1)} of \propo{prop-ab-gen},
holds in $\br(\cals)$ even if $x \cdot \langle s_1s_2\rangle^k$ or $x \cdot \langle
s_2s_1\rangle^{m-k-1}$ is not $\sfw$-reduced. It follows that
$t= w*s\in \calt(\cals) \mapsto \ti t:= \hat w *s$ induces a well-defined
bijection between $\calt(\cals)$ and the sets $\wt\calt^2(\cals)$ and
$\wt\calt(\cals) := \{ \ti t : t\in \calt(\cals) \}$ whose inverse can
be obtained as the composition $\ti t \in \wt\calt(\cals)\mapsto \ti t^2 \in
\wt\calt^2(\cals) \mapsto A_t = [\ti t^2]_\ab \in \sfp_\ab(\cals)$.

Now let $x \in \br(\cals)$ be an element which is either not
$\sfw$-reduced or non-positive. Then there exists a factorisation of
the form $x= x'\cdot s'{}^2 \cdot x''$ or respectively $x= x'\cdot s'{}\inv \cdot x''$
in which the sum of $\br$-lengths of $x'$ and $x''$ is smaller than
that of $x$. Let $s \in \cals$ be any generator. Then $x*s^2$ can be
factorised as
\[
x*s^2 = (x'\cdot s'{}^2 \cdot x'') *s^2 =
\big(x'*s'{}^2 \big) \cdot \big((x'\cdot x'') *s^2\big) \cdot \big(x'*s'{}^{-2} \big)
\]
in the first case, and as
\[
x*s^2 = (x'\cdot s'{}\inv \cdot x'') *s^2 =
\big(x'*s'{}^{-2} \big) \cdot \big((x'\cdot s'\cdot x'') *s^2\big) \cdot \big(x'*s'{}^2 \big)
\]
in the second case. Using this procedure one can express any
RS-generator $x*s^2$ as a word in the alphabet $\wt\calt^2(
\cals)$. Consequently, $\wt\calt^2(\cals)$ generates $\sfp(\cals)$.

The minimality of the generator system $\wt\calt^2(\cals)$ follows
from \refthm{thm-trans}. Indeed, every proper subset of
$\wt\calt^2(\cals)$ generates a proper subgroup of $\sfp_\ab(\cals)$.

Finally, let us observe that if $\ell(w\cdot s) = \ell(w) -1$ for some $w\in
\br(\cals)$ and  $s \in \cals$, then $w= w' \cdot s$ with $w'\equiv w\cdot s$ and
hence $\hat w *s = w'*(s*s) = w'*s$. Thus every RS-generator $\ti t^2
\in \wt\calt^2(\cals)$ can be represented in the form $\ti t^2 = \hat w*
s^2$ with a $\sfw$-reduced word $w\cdot s$.
\qed

\bigskip
Let us now present a procedure that will allow us to describe the possible factorisation
of a given $a \in \sfp(\cals)$ as a product of RS-generators $\ti t^2_i =
\hat w_i * s_i^2 \in \wt\calt^2(\cals)$.

\newdefi{def-hurw} A \slsf{factorisation} of \slsf{length} $l$ of an element $x$
of a group $G$ is an expression $\mbff= f_1 \cdot f_2 \cdot \ldots \cdot f_l$ whose evaluation
in $G$ gives $x$. A \slsf{Hurwitz move} is a transformation of such an $\mbff=
f_1 \cdot f_2 \cdot \ldots \cdot f_l$ whereby a pair $f_i\cdot f_{i+1}$ is replaced by $f_{i+1}\cdot
(f_{i+1}\inv f_i\vph f_{i+1}\vph)$ or by $(f_i\vph f_{i+1}\vph f_i\inv) \cdot f_i$
and the remaining factors remain unchanged. We say either that $f_i$ is
shifted to the right and $f_{i+1}$ is shifted with conjugation to the left or,
respectively, that $f_i$ is  shifted with conjugation.

A \slsf{Hurwitz transformation} is a sequence of Hurwitz moves. Two
factorisations $\mbff$ and $\mbff'$ connected by a Hurwitz transformation are
called \slsf{Hurwitz equivalent}.

{\rm Further details see \eg\ in \cite{Kh-Ku}}.
\end{defi}

\newthm{thm-hurw-W} \slsf{(Hurwitz Problem in Coxeter--Weyl Groups)}. Let
$\mbft = t_1 \cdot t_2 \cdot \ldots \cdot t_l$ be a factorisation of the identity $1\in
\sfw(\cals)$ into quasireflections $t_i\in \calt(\cals)$. Then the length $l$ is
even, $l=2l'$, and $\mbft$ is Hurwitz equivalent to a factorisation $\mbft'$
into squares of quasireflections, \ie, $\mbft'= t_1' \cdot t_2' \cdot \ldots \cdot t_l'$ with
$t'_{2i-1} \equiv t'_{2i}\in \calt(\cals)$.
\end{thm}

For the special case of finite Weyl groups $W$ corresponding to simple complex
Lie algebras, this result was obtained by Kanev \cite{Kan}, \slsf{Proposition 2.3}.
His proof exploits another technique, namely, the
geometry and combinatorics of the corresponding root systems.

\proof Since the set of quasireflections $\calt(\cals)$ is invariant under
conjugations, any Hurwitz transformation of $\mbft$ changes $\mbft$ into a
factorisation into quasireflections.

For a quasireflection $t\in \calt(\cals)$ define its \slsf{width} as the smallest
$k$ such that $t$ can be represented in the form $t\equiv w*s$ with $w \in
\sfw(\cals)$, $s \in \cals$, and such that $k=\ell(w)$.

Shifting the first factor $t_1$ with conjugation to the right end, we obtain a
factorisation $t_2 \cdot t_3 \cdot \ldots \cdot t_l\cdot \ti t_1$. Comparing it with the original
factorisation $\mbft = t_1 \cdot t_2 \cdot \ldots \cdot t_l \equiv1$, we obtain the equality $\ti
t_1\equiv t_1$. Consequently, a cyclic permutation of the factors of $\mbft$ gives a
Hurwitz equivalent factorisation. Thus we may assume that the last factor
$t_l$ has the largest width among all $t_i$.

Now let us represent each $t_i$, $i=1,\ldots ,l-1$, in the form $t_i \equiv w_i *s_i$
with the smallest possible length $\ell(w_i)$. Fix a reduced expression $\mbfw_i$
for each $w_i$. Observe that inverting an expression $\mbfw$ of any element $w
\in \sfw(\cals)$ produces an expression for the inverse element $w \inv$. In view
of this, we denote by $\mbfw\inv$ the inversion of an expression $\mbfw$.

Now consider the expression $\mbfw := \prod_{i=1}^{l-1} \mbfw_i\vph \cdot s_i\vph  \cdot \mbfw_i\inv$. Its
evaluation in $\sfw(\cals)$ gives $t_l\inv \equiv t_l$, and its product with $t_l$
gives $1$. By the Strong Exchange Property there exists a letter $s$ in
$\mbfw$ whose removal from $\mbfw$ gives an expression for $1$.

Consider first the case when $s$ is the middle letter $s_j$ of the expression
$\mbfw_j\vph \cdot s_j\vph \cdot \mbfw_j\inv$ for some $t_j$. Then the remaining pieces
$\mbfw_j\vph \cdot \mbfw_j\inv$ cancel, and hence
\[
t_1 \cdot \ldots \cdot t_{j-1} \cdot t_j \cdot  t_{j+1} \cdot \ldots \cdot t_{l-1} \equiv
t_1 \cdot \ldots \cdot t_{j-1} \cdot t_{j+1} \cdot \ldots \cdot t_{l-1}  \cdot t_l
\]
On the other hand, shifting with conjugation the factor $t_j$ in $t_1 \cdot \ldots \cdot
t_{j-1} \cdot t_j \cdot t_{j+1} \cdot \ldots \cdot t_{l-1}$ to the right end, we obtain a
factorisation $t_1 \cdot \ldots \cdot t_{j-1} \cdot t_{j+1} \cdot \ldots \cdot t_{l-1} \cdot \ti t_j$.
Consequently, $\ti t_j \equiv t_l$ and $t_1 \cdot \ldots \cdot t_{j-1} \cdot t_j \cdot t_{j+1} \cdot \ldots \cdot
t_{l-1} \cdot t_l$ is Hurwitz equivalent to the factorisation $t_1 \cdot \ldots \cdot t_{j-1} \cdot
t_{j+1} \cdot \ldots \cdot t_{l-1} \cdot t_l \cdot t_l$. Since in this case $t_1 \cdot \ldots \cdot t_{j-1} \cdot
t_{j+1} \cdot \ldots \cdot t_{l-1} \equiv1$, we can complete the proof using induction on~$l$.

\smallskip
It remains to consider the case when the letter $s$ appears, say, in the
initial subword $\mbfw_j$ in the expression $\mbfw_j\vph \cdot s_j\vph \cdot
\mbfw_j\inv$ for some $t_j$. (The case when the letter $s$
appears in the final subword $\mbfw_j\inv$ can be treated similarly.) Write
$\mbfw_j$ in the form $\mbfw'\,s\, \mbfw''$. Then $\mbfw'\mbfw'' \equiv (\mbfw's
\mbfw'{}\inv) \cdot \mbfw'\,s\, \mbfw''$ and hence
\begin{equation}\eqqno(w-prim)
\mbfw'\mbfw''\cdot s_j \cdot \mbfw_j\inv \equiv
(\mbfw' s \mbfw'{}\inv) \cdot (\mbfw_js_j\mbfw_j\inv ).
\end{equation}
Now shift $t_l$ with conjugation to the left in-between $t_{j-1}$ and $t_j$.
This gives us a factorisation of the form $t_1 \cdot \ldots \cdot t_{j-1} \cdot \ti t_l \cdot t_j
\cdot t_{j+1} \cdot \ldots \cdot t_{l-1}$. Comparing it with \eqqref(w-prim) we obtain the
equality $\ti t_l \equiv \mbfw' s \mbfw'{}\inv$. Observe that the width of $\ti t_l
\equiv \mbfw' s \mbfw'{}\inv$ is at most $\ell(\mbfw') < \ell(\mbfw' s \mbfw'')$ and
hence less than the width of $t_l$. Now we can use induction on the sum of
the widths of $t_i$.
\qed

\medskip
Now consider the following situation. Let $G$ be a subgroup of $\sfw(\cals)$.
Then $G$ acts on $\sfp_\ab(\cals)$ by conjugation and we denote by
$\sfp_\ab(\cals)_G$ the group of \slsf{coinvariants} (see \cite{Bro}).
Recall that $\sfp_\ab(\cals)_G$ is the quotient of $\sfp_\ab(\cals)$ by the
subgroup generated by the elements of the form $w* A -A$ with $w\in G$ and $A \in
\sfp_\ab(\cals)$. Since $\sfw(\cals)$ permutes basis elements in $\sfp_\ab(\cals)$,
the group $\sfp_\ab(\cals)_G$ is a free abelian group with a basis given by the
quotient set $\calt(\cals) / G$. Another description of $\sfp_\ab(\cals)_G$ is
as the quotient of $\sfp(\cals)$ by the commutator group $[\wt G, \sfp(\cals)]$
where $\wt G$ is the pre-image of $G$ in $\br(\cals)$. We denote the elements
of $\calt(\cals) / G$ by $G\cdot t$ and the elements of the induced basis
in $\sfp_\ab(\cals)_G$ by $A_{G\cdot\theta}$ or $A_{G\cdot t}$.

\newthm{thm-subgr} Let $G$ be a subgroup of $\sfw(\cals)$ and $x$ an element
of $\sfp(\cals)$ which can be represented as the product $x = \prod_i \hat
t_i^{\epsilon_i} \cdot \prod_j [x_{2j-1}, x_{2j}]$ such that $\hat t_i$ are quasigenerators,
$\epsilon_i=\pm1$, and $[x_{2j-1}, x_{2j}]$ are commutators of some $x_j \in \br(\cals)$.
Assume that the projections of $x_i$ and $\hat t_j$ to $\sfw(\cals)$ lie in
$G$.  Then the projection $[x]_G$ of $x$ to $\sfp_\ab(\cals)_G$ lies in the
free abelian group generated by basis elements $A_{G\cdot t}$ with $t \in G \bigcap
\calt(\cals)$.
\end{thm}

\proof First, we reduce the general case to the special one with no
commutators. For this purpose we represent each $x_j$ in the form $x_j = x'_j
\cdot x''_j$ so that $x'_j$ is a product of quasigenerators projecting to
$G$, and $x''_j$ is a product of squares of quasigenerators.  Using the
commutation relations $[x, y\cdot z] = [x, y] \cdot (y*[x, z])$ and $[x\cdot y, z] =
(x*[y, z]) \cdot [x, z]$ we expand each commutator $[x_{2i-1}, x_{2i}]$ as a
product of quasigenerators projecting to $G$ and commutators $[y, \hat
t^2]$ such that $y$ projects to $G$ and $\hat t^2$ is a squared
quasigenerator.  Since such $[y, \hat t^2]$ lie in $\sfp(\cals)$ and project
to zero in $\sfp_\ab(\cals)_G$, we obtain the desired reduction.

\smallskip%
Observe that Hurwitz moves do not destroy the properties of factors $\hat
t_i^{\epsilon_i}$ listed in the hypothesis. Thus applying \refthm{thm-hurw-W}, we
transform the original factorisation into a new one, still denoted in the same way,
in which $t_{2i-1} \equiv t_{2i} \in \sfw(\cals)$. This mean that it is enough to
consider the case when $x$ is the product of \emph{two} factors, $x= \hat
t_1^{\epsilon_1} \cdot \hat t_2^{\epsilon_2}$ with $t_1\equiv t_2$. In the case when the signs $\epsilon_i$
are equal, say, $\epsilon_1 = \epsilon_2 = +1$, we can rewrite this product in the form $\hat
t_1\vph \hat t_2\inv \hat t_2^2$. Since $t_2 \in G \bigcap \calt(\cals)$ by our assumptions,
the square $\hat t_2^2$ has the desired form.  Thus we may
additionally assume that $x= \hat t_1\vph \cdot \hat t_2\inv$.

We claim that the condition $t_1\equiv t_2 \in \sfw(\cals)$ implies that $\hat t_2$
can be obtained from $\hat t_1$ by conjugation with some $z \in \sfp(\cals)$.
For this purpose we write $\hat t_i = y_i * s_i$ with $s_i\in \cals$ and $y_i \in
\br(\cals)$ and conjugate both $\hat t_i$ by $y_1\inv$. This reduces the
situation to the special case when $\hat t_1 = s_1$ is a usual generator of
$\br(\cals)$. Adjusting the notation, we still have $\hat t_2 = y_2 * s_2$.
Let $w$ be the image of $y_2$ in $\sfw(\cals)$. Our claim would follow from
the equality $\hat w * s_2 = s_1$ in the braid group. Let us prove it.
Assume first that $\ell(w s_2) < \ell(w)$. Then by the Exchange Property $w\equiv w'
s_2$ with $\ell(w') = \ell(w)-1$. In this case $\ell(w's_2) > \ell(w')$, $\hat w * s_2 =
(\hat w' s_2) * s_2 = \hat w' * s_2$, and still $\hat w' * s_2 \equiv s_1$ in
$\sfw(\cals)$. Thus we may assume that $\ell(w s_2) > \ell(w)$. Observe that the
equality $w * s_2 \equiv s_1$ is equivalent to $w s_2 \equiv s_1 w$. As $\ell(w s_2) >
\ell(w)$, the expressions $\mbfw s_2$ and $s_1\mbfw$ are reduced for any reduced expression $\mbfw$ of~$w$.
Consequently, $\hat w \cdot s_2 = s_1 \cdot \hat w$ in the braid group.
This is the desired identity $\hat w * s_2 = s_1$, which proves the claim.

\smallskip
Summing up, it remains to consider a commutator $\hat t \,z\, \hat t \inv \, z\inv$
where $z \in \sfp(\cals)$ and $\hat t \in G$. By definition, such elements project to $0$
in $\sfp_\ab(\cals)_G$. This finishes the proof.
\qed

\newdefi{def-gars} A Coxeter system $(\cals, M)$ is \slsf{irreducible} if its
Coxeter graph is connected, and has \slsf{finite type} if the group
$\sfw(\cals)$ is finite. It is known that every irreducible Coxeter system
$(\cals, M)$ of finite type has the unique \slsf{longest element} $w_\circ=
w_\circ(\cals)\in \sfw(\cals)$. The canonical lift $\hat w_\circ \in \br(\cals)$ is called
the \slsf{Garside element} of $\br(\cals)$ and is denoted by $\idel(\cals)$.
\end{defi}

The classification theory of irreducible Coxeter groups of finite type (see
\cite{Bou} or \cite{Hum}) says that such a group is either the Weyl group of a
simple complex Lie group, or the dihedral group $\cald_{2m}$ with the Coxeter
system on $\cals =\{s_1, s_2\}$ given by relation $m_{s_1,s_2} =m$, $m=5,7,
8,9,\ldots$, or else one of the groups $\rmh_3$, $\rmh_4$. It is known that the
longest element $w_\circ \in \sfw(\cals)$ has the property $\ell(w_\circ) = \ell(w_\circ w) +
\ell(w)$ for any $w \in \sfw(\cals)$. For the properties of the Garside element we
refer the reader to \cite{Br-Sa} and \cite{Del}.

\newlemma{lem-gars-sq} The square $\idel^2(\cals)$ of the Garside element lies
in\/ $\sfp(\cals)$ and is equal in $\sfp_\ab(\cals)$ to the sum $\sum_{t\in
 \calt(\cals)} A_t$.
\end{lem}

\proof The property $\ell(w_\circ w) = \ell(w_\circ) - \ell(w)$ applied to $w=w_\circ$ implies that
$w_\circ$ is an idempotent. Thus $\idel^2 =\hat w_\circ^2$ lies in $\sfp(\cals)$.

Take any generator $s \in\cals$ and any reduced expression $\mbfw$ of $w_\circ$.
Then $\ell(sw_\circ) = \ell(w_\circ) -\ell(s) < \ell(w_\circ)$. By the Strong Exchange Property,
removing an appropriate letter $s'$ from $\mbfw$ we obtain a word $\mbfw_1$
which is an expression for $sw_\circ$. Comparing the lengths we see that $\mbfw_1$
is reduced. Writing $\mbfw_1$ as $s_2\cdot \ldots \cdot s_l$ with $l:= \ell(w_\circ)$, we see that
$ss_2\cdot \ldots \cdot s_l$ is a reduced expression for $w_\circ$. Since $w_\circ$ is an
idempotent, $s_l\cdot \ldots \cdot s_2s$ is also a reduced expression for $w_\circ$.
Consequently, $\idel^2 = \hat s\hat s_2\cdot \ldots \cdot \hat s_l \cdot \hat s_l\cdot \ldots \cdot\hat s_2\hat
s$. Denoting $s_1:=s$, transform the latter expression into the product
\[\textstyle
\idel^2 = \prod_{i=1}^l  (\hat s_1\inv\hat s_2\inv\ldots\hat s_{i-1}\inv) * \hat s_i^2.
\]
It shows that $\idel^2$ is equal in $\sfp_\ab(\cals)$ to a sum $\sum_{t\in \calt(\cals)} n_t A_t$
with non-negative integers $n_t$ such that $n_s >0$ for any $s\in \cals$. Since $\idel^2$ lies
in the centre of $\br(\cals)$, the sum $\sum_{t\in \calt(\cals)} n_t A_t$ is invariant under
the action of $\sfw(\cals)$. Since $\sfw(\cals)$ permutes the basis elements $A_t$ of $\sfp_\ab(\cals)$
and any quasireflection $t\in \calt(\cals)$ is conjugated to some $s\in \cals$, we obtain
the inequality $n_t \geq1$ for all  $t\in \calt(\cals)$.

On the other hand, it is shown in \cite{Br-Sa} and \cite{Del} that $\idel^2= \wh \Pi^h$ where $\Pi$ 
is the so-called \slsf{Coxeter element} in $\sfw(\cals)$ and $h$ is the \slsf{Coxeter number} of the system
$\cals$, see \eg \cite{Bou} or \cite{Hum} for definitions. The formula for the
Coxeter number can be reformulated as the equality of the length $l= \ell(w_\circ)$ and
the number of quasireflection. This shows that $\sum_{t\in \calt(\cals)} n_t = l$.
Consequently, all $n_t$ are equal to $1$, and the lemma follows.
\qed

\newsubsection[map-mu-pp]{Combinatorial structure of $\map_g$.}
In this paragraph we solve certain factorisation problems in the mapping class
group, which is the key ingredient in our proof of the \slsf{Main Theorem}.

We use special finite presentations of $\map_g$ and $\map_{g,[1]}$ due to
Wajnryb \cite{Waj} and Matsumoto \cite{Ma} which realise these groups as
quotients of braid groups corresponding to certain Coxeter systems $\cals_g$ with
additional relations given in terms of Garside elements of appropriate
subsystems $\cals'$ of $\cals_g$. Let us give a geometric description of these
relations. For further details see \eg\ \cite{Bi-1} and \cite{Ger}.

\begin{figure}
\includegraphics[height=1.7in]{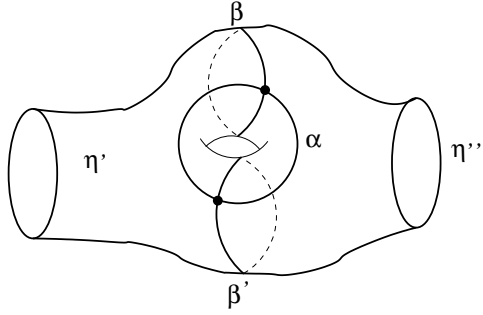}
\caption{Curves in a chain configuration.} \label{fig-chain}
\end{figure}

\newdefi{def-map-pres} {\rm (}Chain and lantern relations.{\rm )} \rm Consider a surface
$\calc$ which is a torus with 2 holes. Denote its boundary circles
by $\eta',\,\eta''$. Consider the curves $\alpha, \beta,\beta'$ on
the surface as shown on \reffig{fig-chain}. Note that if
$\alpha,\beta, \beta'$ are embedded circles on a surface $\Sigma$
such that $\beta$ and $\beta'$ are disjoint and each meets $\alpha$
transversally at a single point, then a collar neighbourhood $U$ of
the graph $\alpha\cup \beta \cup \beta'$ is a torus with 2 holes and
the whole configuration is diffeomorphic to the one on
\reffig{fig-chain}.  We call such a surface $\calc \subset\Sigma$
and the whole configuration $(\calc, \alpha,\beta, \beta')$ a
\slsf{chain} in $\Sigma$ defined by $\alpha,\beta, \beta'$. For any
chain configuration, the \slsf{chain relation element} is defined by the formula
\begin{equation}\eqqno(chain)
C(\alpha,\beta, \beta') := \big(T_{\beta} T_\alpha T_{\beta'} \big)^4
\big(T_{\eta'} T_{\eta''}\big)\inv .
\end{equation}

\begin{figure}
\includegraphics[height=2.5in]{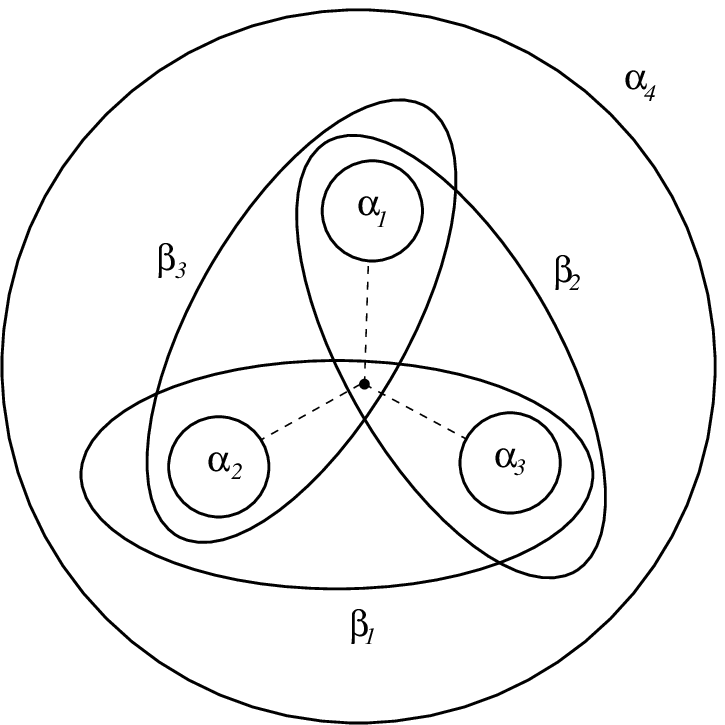}
\caption{Curves in a lantern configuration.} \label{fig-lant}
\end{figure}

\medskip
Consider a surface $\call$ which is a sphere with 4 holes. Denote
its boundary circles by $\alpha_1,\ldots,\alpha_4$. Realise it as a
disc with $3$ holes and consider the curves
$\beta_1,\beta_2,\beta_3$ on it as shown on \reffig{fig-lant}.
Observe that if $\beta_1, \beta_2$ are embedded circles on a surface
$\Sigma$ which meet at two points, then a collar neighbourhood $U$
of the graph $\beta_1 \cup \beta_2$ is a disc with $3$ holes and the
whole configuration is diffeomorphic to the one on
\reffig{fig-lant}. We call such a surface $\call \subset\Sigma$ and
the whole configuration $(\call, \alpha_i,\beta_j)$ a \slsf{lantern}
in $\Sigma$ defined by $\beta_1, \beta_2$. For any lantern
configuration, the \slsf{lantern relation element} is defined by
the formula
\begin{equation}\eqqno(lant)
L(\beta_1, \beta_2) := \big(T_{\alpha_1} T_{\alpha_2} T_{\alpha_3} T_{\alpha_4} \big)
\big( T_{\beta_1} T_{\beta_2} T_{\beta_3} \big)\inv.
\end{equation}

A \slsf{chain} or, respectively, \slsf{lantern relation} is the
equality $C(\alpha,\beta, \beta') = 1$ or, respectively, $L(\beta_1,
\beta_2) = 1$ for the corresponding relation elements. A
configuration (or relation) is called \slsf{non-separating} if all
the circles involved in it are non-separating.

\end{defi}

\smallskip
As suggested by the terminology, the relations above hold in the mapping class
group $\map_{g,k,[l]}$ of the surface of genus $g$ with $k$ marked points and
$l$ boundary circles. (Our notation uses the fact that $\map_{g,k,[l]}$ can be defined
as the mapping class group of the surface of genus $g$ with $k+l$ marked points
such that $l$ of them, say, $z_1,\ldots,z_l$, are \slsf{framed}, \ie,
equipped with a trivialisation of the tangent plane $T_{z_i}\Sigma$.) We will only consider
the case with at most one marked point, \ie, $k+l\leq1$, and shorten our notation
to $\map_{g,1}$ or $\map_{g,[1]}$, dropping the vanishing index $k$ or $l$.

\newdefi{def-S-g}
For $g=1$, we set $\cals_1 := A_2$ so that the graph of $\cals_1$ consists of
two vertices $s_1$ and $s_2$ connected by an edge. For $g\geq2$, the graph
of $\cals_g$ is defined as the extension of the Dynkin linear graph $A_{2g}$
by a single vertex $s_0$ connected to $s_4$. Thus, the graph of $\cals_g$ looks
like
\[
\xymatrix@C-10pt@R-11pt{
& & & s_0\ar@{-}[d]
\\
 s_1 \ar@<2pt>@{-}[r] & s_2\ar@<2pt>@{-}[r] & s_3 \ar@<2pt>@{-}[r] &
s_4\ar@<2pt>@{-}[r]  & s_5 \ar@<2pt>@{-}[r] & s_6 \ar@<2pt>@{-}[r] &
\ldots \ar@<2pt>@{-}[r] & s_{2g-1} \ar@<2pt>@{-}[r] & s_{2g}
}
\]
\begin{figure}
\includegraphics[width=16cm]{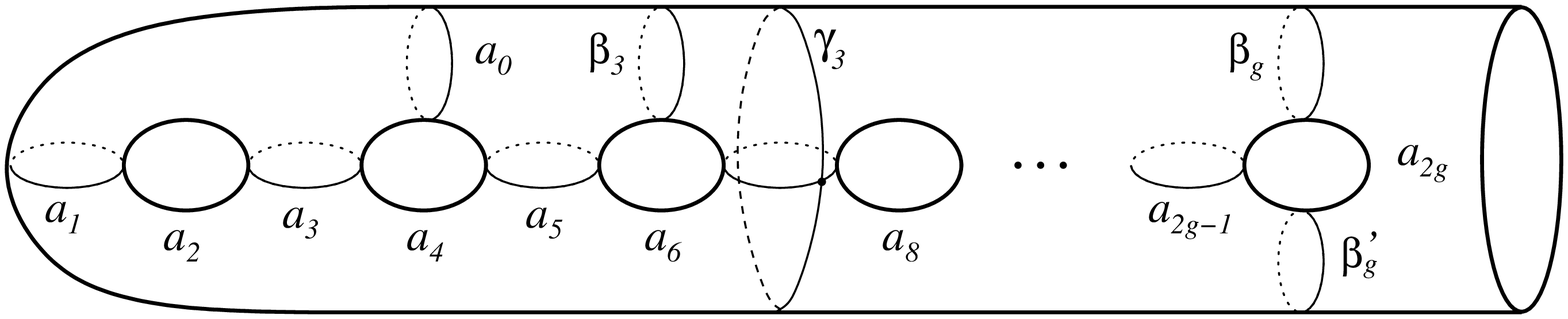}
\caption{Geometric generators of $\map_{g,[1]}$.}
\label{fig-Delta-g}
\end{figure}

Each element $\hat t \in \br(\cals_g)$ conjugated to a generator $s_i$ is called
a \slsf{quasigenerator}. The set of quasigenerators is denoted by $\wh\calt(\cals_g)$.
\end{defi}

\newthm{thm-map-br} The group $\map_{g,[1]}$ is isomorphic to the quotient
of $\br(\cals_g)$ obtained by adding the following relations:
\begin{itemize}
\item In the case $g=1\!:$ No additional relations.
\item In the case $g=2\!:$ A single \emph{non-separating} chain relation
 $C(\alpha,\beta,\beta') =1$.
\item In the case $g\geq3\!:$ A single \emph{non-separating} chain relation
 $C(\alpha,\beta,\beta') =1$ and a single \emph{non-separating} lantern relation
 $L(\beta',\beta'') =1$.
\end{itemize}
The braid generators $s_i$, $i=0;1,\ldots, 2g$, can be realised as Dehn
twists along the curves $a_i$ shown on \reffig{fig-Delta-g}.
In particular, $a_0 = \beta_2$, $a_1 = \beta_1$, $a_{2i} = \alpha_i$ for $i=1,\ldots,g$, and
each $a_{2i-1}$ is $\zz$-homologous to $[\beta_{i-1}] + [\beta_i]$ for $i=2,\ldots,g$.
\begin{itemize}
\item The kernel of the homomorphism $\map_{1,[1]} \to \map_1$ is the central
 free abelian group generated by $(T_{\alpha_1}T_{\beta_1})^6 = (s_1s_2)^6$. In the
 case $g\geq2$, the kernel of the homomorphism $\map_{g,[1]} \to \map_g$ is
 normally generated by the product $T_{\beta_g}\vpht T_{\beta'_g} \inv $ where
 $T_{\beta_g}, T_{\beta'_g}$ are the Dehn twists along the curves $\beta_g, \beta'_g$.
\qed
\end{itemize}
\end{thm}

This presentation was found by Wajnryb \cite{Waj}, with an error corrected in
\cite{Bi-Waj}. In \cite{Ma}, Matsumoto has found simple explicit words expressing
the relations in terms of Garside elements of certain Coxeter subsystems of $\cals_g$.

\newthm{thm-matsu} \begin{itemize}
\item[\slip]  The chain relation can be given by the element $C_0 := \idel^{-4}(\rma_4)
 \idel^2(\rma_5)$ where $\rma_4$ and $\rma_5$ are the Coxeter subsystems
 $\{s_2, s_3,s_4,s_0\}$ and $\{s_1, s_2, s_3, s_4,s_0\}$, respectively.
\item[\sliip] The lantern relation can be given by the element $C_0\cdot
 \idel^{-2}(\rme_6) \idel(\rme_7)$ where $C_0$ is the chain relation element
 above and $\rme_6$, $\rme_7$ are the Coxeter subsystems $\{s_0,s_2, \ldots, s_6\}$
 and $\{s_0,s_1, s_2, \ldots, s_6\}$, respectively.
\item[\sliiip] The relation $T_{\beta_g} =T_{\beta'_g}$ can be replaced by the
 commutator relation $[T_{\beta_g}, \idel^2(\rma_{2g})] =1$ where $\rma_{2g}$ is
the Coxeter subsystem $\{s_1,s_2, \ldots, s_{2g}\}$.
\qed
\end{itemize}
\end{thm}

\remark It is known that any two non-separating chain (or lantern)
configurations are conjugate by some $f \in \diff_+(\Sigma, z_0)$.
This implies the fact, used implicitly in the above theorem, that
the subgroup normally generated by a single non-separating chain
relation $C(\alpha,\beta,\beta')$ is independent of a specific
choice of the chain configuration, and the same holds for lantern
configurations. For the proof of this and further details, we refer
to \cite{Ger}.

\smallskip
We call $\idel^2(\rma_{2g})$ the
\slsf{hyperelliptic element}, the commutator $[T_{\beta_g}, \idel^2(\rma_{2g})]$
the \slsf{(basic) hyperelliptic relation element}, and $[T_{\beta_g},
\idel^2(\rma_{2g})]=1$ the \slsf{hyperelliptic relation}. The reason for this
terminology is that the element $\idel^2(\rma_{2g}) \in \br(\cals_g)$ represents
the hyperelliptic involution of $\Sigma$ which can be realised on
\reffig{fig-Delta-g} as the rotation by $180^\circ$ about the horizontal axis going
along the chain of curves $a_1,a_2,\ldots ,a_{2g}$.  Besides, we call the relation
$\idel^2(\rme_6) =\idel(\rme_7)$ the \slsf{modified lantern relation} and the
product $\idel^{-2}(\rme_6) \idel(\rme_7)$ the \slsf{(basic) modified lantern
relation element}. The same names will be used for any conjugates of these
elements.

\medskip
Observe that all $s_i$ are conjugate in $\br(\cals_g)$, so that quasigenerators $\hat t\in \wh\calt(\cals_g)$
are preferred lifts of Dehn twists along non-separating curves~$\delta$. The $\zz_2$-homology class
$[\delta] \in \sfh_1(\Sigma, \zz_2)$ is determined by the image of $\hat t$ in $\sfw(\cals_g)$.
Hence, we obtain a well-defined set-theoretic $\sfw(\cals_g) $-equivariant map
$\calt(\cals_g) \to \calh_g := \sfh_1(\Sigma, \zz_2) \bs \{ 0\}$.
(This map is bijective if $g=1,2,3$, but if $g\geq4$, then the set $\calt(\cals_g)$ is infinite,
and the map is only surjective.) Furthermore, the projection of a quasigenerator $\hat t$ to $\Sp(2g, \zz_2)$
has order $2$. Therefore there exists a well-defined homomorphism $\sfw(\cals_g) \to \Sp(2g, \zz_2)$.

\newdefi{def-iw-group} The kernel $\IW_g := \ker\big( \sfw(\cals_g) \to
\Sp(2g, \zz_2) \big)$ is called the \slsf{Weyl--Torelli group} of the surface
$\Sigma$.
\end{defi}

It follows that $\calh_g$ is the quotient set $\calt(\cals_g) / \IW_g$.
Consequently, the group $\zz\langle \calh_g \rangle$ is naturally isomorphic to the
coinvariant group $\sfp_\ab(\cals_g)_{\IW_g}$. This allows us to apply
\refthm{thm-subgr} with the group $\zz\langle \calh_g \rangle$ instead of $\zz\langle
\calt(\cals_g) \rangle$. Moreover, we can reduce the extension $1 \to \sfp(\cals_g) \to
\br(\cals_g) \to \sfw(\cals_g) \to 1$ to an extension $1 \to \zz\langle \calh_g \rangle \to E \to
\sfw(\cals_g) \to 1$. However, the group $\sfw(\cals_g)$ is bigger than $\Sp(2g, \zz_2)$
for all $g\geq3$,, and their difference --- the Weyl--Torelli group $\IW_g$ ---
is generated by lantern relation elements and is therefore invisible in $\map_g$.
Worse still, $\IW_g$ is infinite for all $g\geq4$,  which makes effective calculation
in $\IW_g$ difficult. In order to see which part of $\IW_g$ can be factored out,
we need a description of the kernel $\ker \big(\br(\cals_g) \to \Sp(2g,\zz_2) \big)$.

\newprop{prop-map-gen} The kernel of the natural homomorphism $\br(\cals_g)
\to \Sp(2g,\zz_2)$ is\break generated by squared quasigenerators $\hat t^2$ and by
modified lantern relation elements\break $x* \big( \idel^{-2}(\rme_6)
\idel(\rme_7) \big)$, $x\in \br(\cals_g)$.
\end{prop}

\proof By \lemma{lem-gars-sq}, every element of the kernel of the homomorphism
$\br(\cals_g) \to \map_g$ can be represented as a product of squared
quasigenerators $\hat t^2$ and lantern relation elements.  Consequently, it is
sufficient to find the corresponding presentation in the mapping class group
$\map_g$.

Denote by $\bfg$ the kernel $\ker \big(\map_g \to \Sp(2g,\zz_2) \big)$ and by $\bfg_0$
the subgroup generated by products $T_{\delta_1}^{\epsilon_1} T_{\delta_2}^{\epsilon_2}$ with $\epsilon_i=\pm1$
and $[\delta_1] = [\delta_2] \neq0 \in \sfh_1(\Sigma,\zz_2)$.  We use $\equiv$ to denote equality
modulo $2$, \ie, the equality in $\zz_2$, or in $\sfh_1(\Sigma, \zz_2)$, or in
$\Sp(2g,\zz_2)$.  Let $\alpha_1,\ldots,\alpha_g; \beta_1,\ldots,\beta_g$ be the geometric basis of $\Sigma$
fixed above.

Fix some $f \in \bfg$ and write $f$ as a product of Dehn twists along
non-separating curves, $f= \prod_{i=1}^n T_{\delta_i}^{\epsilon_i}$. The braid relations $T_\gamma
T_\delta^\epsilon = T_{\delta'}^\epsilon T_\gamma$ with $\epsilon=\pm1$ and $\delta'= T_\gamma(\delta)$
allow us to use Hurwitz moves to re-organise the product by moving any single Dehn twist $T_\gamma$
in both directions without changing it.  Note that the generators $T_{\delta_1}^{\epsilon_1}
T_{\delta_2}^{\epsilon_2}$ of $\bfg_0$ are stable under conjugation.

\begin{stepping}

\item\label{st1} Assume that there exist factors $T_{\delta_i}$ with $\delta_i\cap \beta_1 \equiv0$.
Collect all such factors $T_{\delta_i}$ on the right using the braid relation. Henceforth,
we assume that $\delta_i\cap \beta_1 \equiv1$ for $i=1,\ldots,n_1$ and $\delta_i\cap \beta_1 \equiv0$ for
$i=n_1+1,\ldots,n$.

\item\label{st2} Assume that $\delta_i \cap \delta_j\equiv1$ for some distinct $i,j\leq n_1$. Bring
them together, so that $j=i+1$, and then apply the braid relation. It
transforms $T_{\delta_i}T_{\delta_{i+1}}$ into $T_\gamma T_{\delta_i}$ with $\gamma \equiv \delta_i +\delta_{i+1}$.
Consequently, $\gamma\cap \beta_1 \equiv0$, and we can shift $T_\gamma$ to the right.
This allows us to diminish the number $n_1$ of $\delta_i$ with $\delta_i\cap \beta_1 \equiv1$.
Repeating this procedure, we come to the case when $\delta_i \cap \delta_j\equiv0$ for all
$i,j\leq n_1$.

\item\label{st3} Assume that $n_1>0$, \ie, there exist factors $T_{\delta_i}$ with $\delta_i\cap
 \beta_1 \equiv1$. Take the initial subproduct $\prod_{i=1}^{n_1} T_{\delta_i}$ and
 multiply it by $T_{\beta_1}\inv T_{\beta_1}\vph$ on the left. Apply the braid  relation twice
 to the initial subword $T_{\beta_1}\inv T_{\beta_1}\vph T_{\delta_1}\vph$:
\[
T_{\beta_1}\inv T_{\beta_1}\vph T_{\delta_1}\vph \Rightarrow
T_{\beta_1}\inv T_\delta\vph T_{\beta_1}\vph \Rightarrow
T_{\beta_1}\inv T_\gamma\vph T_\delta\vph,
\]
where $\delta := T_\beta(\delta_1)$ and $\gamma:= T_\delta(\beta_1)$. Then $\delta \equiv \delta_1 +\beta_1$ and $\delta \cap \delta_i \equiv1$
for all $i=2,\ldots,n_1$. Then shift $T_\delta$ to the right behind all $T_{\delta_i}$,
$i=2,\ldots,n_1$:
\[\textstyle
T_{\beta_1}\inv T_\gamma\vph T_\delta\vph \prod_{i=2}^{n_1} T_{\delta_i}\vph
\Rightarrow
T_{\beta_1}\inv T_\gamma\vph \prod_{i=2}^{n_1} T_{\delta_i'}\vph T_\delta\vph,
\]
where $\delta'_i := T_\delta(\delta_i)$. Then $\delta'_i \equiv \delta_i + \delta$ and $\delta'_i \cap \beta_i\equiv0$. Finally,
shift $T_{\beta_1}\inv$ and the newly obtained $T_{\delta_i'}$, $i=2,\ldots,n_1$, to the right. We
obtain a new decomposition $f= \prod_{i=1}^{n+2} T_{\delta_i''}$ of the original $f$
in which $\delta''_i \cap \beta_1\equiv0$ for $i=3,\ldots,n+2$ and $\delta''_1 \cap \beta_1\equiv1$.

\item\label{st4} If $\delta''_2 \cap \beta_1\equiv0$, then $\prod_{i=1}^{n+2} T_{\delta_i''} (\beta_1)\equiv \beta_1
 + \delta''_1 \not\equiv \beta_1$ which is in contradiction with the assumption that $f_* \equiv\id$ in
 $\zz_2$-homology. If $\delta''_1 \cap \delta''_2\equiv1$, transform $T_{\delta''_1}T_{\delta''_2}$
 into $T_{\delta'''_2}T_{\delta''_1}$ with $\delta'''_2 \cap \beta_1\equiv0$. This leads to the same
 contradiction. Thus $\delta''_1 \cap \delta''_2\equiv0$ and $\delta''_1 \cap \beta_1\equiv\delta''_2 \cap \beta_1\equiv1$. But
 then $\beta_1\equiv \prod_{i=1}^{n+2} T_{\delta_i''} (\beta_1)\equiv \beta_1 + \delta''_1 + \delta''_2$ and hence
 $\delta''_1 \equiv \delta''_2$.  This means that $T_{\delta_1''}T_{\delta_2''}$ lies in $\bfg_0$. Shift it
 to the right. In this way we have represented $f$ in the form $\prod_i T_{\delta_i}
 \cdot f_1$ with $f_1 \in \bfg_0$ and $\delta_i \cap \beta_1\equiv0$.

\item\label{st6} Repeat \slsf{Steps \ref{st1}--\ref{st4}} subsequently for
 $\beta_2,\ldots,\beta_g$ instead of $\beta_1$. For each $\beta_k$, the previously reached
 relations $\delta_i \cap \beta_j\equiv0$ with $i=1,\ldots,n$ and $j=1,\ldots,k-1$ remain undestroyed.
 As a result, we represent $f$ in the form $\prod_i T_{\delta_i} \cdot f_2$ with $f_2 \in
 \bfg_0$ and $\delta_i \cap \beta_j\equiv0$ for all $j=1, \ldots,g$. It follows that after this step
 each $\zz_2$-class $[\delta_i]$ is a linear combination of $\zz_2$-classes
 $[\beta_1],\ldots,[\beta_g]$, and hence $\delta_i \cap \delta_j\equiv0$ for each $i,j=1,\ldots,n$.  Consequently,
 the application of the braid relation does not change the $\zz_2$-homology
 class of the remaining $\delta_i$, for if $T_{\delta_i} T_{\delta_j} T_{\delta_i}\inv = T_{\delta'_j}$,
 then $\delta_j \equiv \delta'_j$.

\item\label{st7} Consider $\alpha_g$. Let $\Delta_{\alpha_g}$ be the set of those $\delta_i$ for which
$\alpha_g \cap \delta_i\equiv1$ and let $V_{\alpha_g}$ be the $\zz_2$-vector space spanned by
$\Delta_{\alpha_g}$.  Assume that $\dim_{\zz_2} V_{\alpha_g}\geq3$. Find three circles in
$\Delta_{\alpha_g}$, say, $\delta_1,\delta_2, \delta_3$, which are $\zz_2$-linearly independent. Then
$[\delta_2] \equiv [\delta_1] + [\gamma_2]$ and $[\delta_3] \equiv [\delta_1] +[\gamma_3]$ for some classes $[\gamma_2],
[\gamma_3] \in \sfh_1(\Sigma, \zz_2)$ such that $[\delta_1], [\gamma_2], [\gamma_3]$ are
$\zz_2$-linearly independent.  Note that $\gamma_1 \cap \alpha_g \equiv (\delta_i-\delta_1)\cap \alpha_g \equiv 0$.
Realise the classes $[\gamma_2], [\gamma_3] \in \sfh_1(\Sigma, \zz_2)$ by embedded curves $\gamma_2,
\gamma_3 \subset\Sigma$ disjoint from $\delta_1$ and from each other. Pick a point $z_0 \in
\Sigma$ disjoint from $\delta_1,\gamma_2, \gamma_3$ and choose embedded arcs $a_1, a_2, a_3$
connecting $z_0$ with $\delta_1,\gamma_2, \gamma_3$, respectively, and disjoint
(except for the end points) from each other and from $\delta_1,\gamma_2, \gamma_3$. Then a
collar neighbourhood $U$ of the graph formed by $\delta_1,\gamma_2, \gamma_3$ and $a_1, a_2,
a_3$ is a disc with three holes bounded by $\delta_1,\gamma_2, \gamma_3$, see \reffig{fig-lant}.
Multiply the product $\prod_{i=1}^n T_{\delta_i}$ by the corresponding lantern relation
\[
T_{\delta_1 +\gamma_2}\inv T_{\delta_1+ \gamma_3}\inv T_{\gamma_2+ \gamma_3}\inv
T_{\delta_1} T_{\gamma_2} T_{\gamma_3} T_{\delta_1 +\gamma_2+ \gamma_3} = \id,
\]
see \reffig{fig-lant}. The factors $T_{\delta_1 +\gamma_2}\inv T_{\delta_2}$ and $T_{\delta_1
 +\gamma_3}\inv T_{\delta_3}$ cancel out.
Shift the factor $T_{\delta_1}^2 \in \bfg_0$ to the right. Now we
obtain the product $\prod_{i=1}^{n+1} T_{\delta'_i}$ which contains \emph{less}
factors $T_{\delta'_i}$ with $\delta'_i \cap \alpha_g \equiv1$.
Repeating this procedure, we come to the case when $\dim_{\zz_2} V_{\alpha_g}\leq2$.

\item\label{st8} We maintain the notation $\Delta_{\alpha_g}$ and $V_{\alpha_g}$ of
 \slsf{Step \ref{st7}}. In the case $\dim_{\zz_2} V_{\alpha_g}=0$ we have $\delta_i\cap \alpha_g
 \equiv0$ for every factor $T_{\delta_i}$ and proceed to the next step. Assume that
 $\dim_{\zz_2} V_{\alpha_g}=1$. Then $\delta_i\equiv\delta_j$ for every $\delta_i,\delta_j \in \Delta_{\alpha_g}$.
 Consequently, $\prod_{i=1}^n T_{\delta_i}(\alpha_g) \equiv \alpha_g + n_{\alpha_g} \delta_1$ where $\delta_1$ is
 any element of $\Delta_{\alpha_g}$ and $n_{\alpha_g}$ is the cardinality of $\Delta_{\alpha_g}$.
 It follows that $n_{\alpha_g}$ is even and we can collect the factors $T_{\delta_i}$
 in pairs $T_{\delta_{2i-1}} T_{\delta_{2i}}$ lying in $\bfg_0$.

\item\label{st9} Assume that $\dim_{\zz_2} V_{\alpha_g}=2$. Take two elements from
$\Delta_{\alpha_g}$, say $\delta_1$ and $\delta_2$, which are linearly independent. Then every
$\delta_i \in \Delta_{\alpha_g}$ is $\zz_2$-homologous to $\delta_1$, or to $\delta_2$, or to $\delta_1 +\delta_2$.
Find $\gamma_1, \gamma_2 \in \sfh_1(\Sigma, \zz_2)$ such that $\delta_1 \cap \gamma_1 \equiv \delta_2 \cap \gamma_2\equiv 1$, $\delta_1
\cap \gamma_2 \equiv \delta_2 \cap \gamma_1\equiv \gamma_1\cap \gamma_2 \equiv 0$, and $\delta_i \cap \gamma_j \equiv0$ for every $\delta_i$
\emph{not} lying in $\Delta_{\alpha_g}$. (This can be done by choosing an
appropriate $\zz_2$-Darboux basis in $\sfh_1(\Sigma, \zz_2)$.) Then
\[\textstyle
\prod_{i=1}^n T_{\delta_i}(\gamma_1) -\gamma_1 \equiv \sum_{i=1}^n (\gamma_1 \cap \delta_i)\delta_i
\]
and hence
\[\textstyle
\gamma_2 \cap \Big( \prod_{i=1}^n T_{\delta_i}(\gamma_1) -\gamma_1 \big)  \equiv n_{12}
\]
where $n_{12}$ is the number of $\delta_i \in \Delta_{\alpha_g}$ with $\delta_i \equiv \delta_1 + \delta_2$.
Consequently, $n_{12}$ is even and we can collect the factors $T_{\delta_i}$ with
$\delta_i \equiv \delta_1 + \delta_2$ in pairs as at \slsf{Step \ref{st7}}. Similarly one
shows that the numbers $n_1$ and $n_2$ of $\delta_i \in \Delta_{\alpha_g}$ with $\delta_i \equiv \delta_1$ or,
respectively, $\delta_i \equiv \delta_2$ are also even. Thus all factors $T_{\delta_i}$ with $\delta_i \in
\Delta_{\alpha_g}$ can be collected in pairs $T_{\delta_{2i-1}} T_{\delta_{2i}}$ lying in $\bfg_0$.

\item\label{st10} After \slsf{Step \ref{st9}}, the remaining factors $T_{\delta_i}$
 fulfil the relation $\delta_i \cap \alpha_g \equiv0$. We repeat \slsf{Steps
  \ref{st7}--\ref{st9}} for $\alpha_{g-1}, \alpha_{g-2}$ and so on. Notice that the
 previously achieved $\zz_2$-orthogonality with $\beta_1,\ldots,\beta_g; \alpha_g, \ldots$
 remains unaffected. Therefore, we conclude that $\bfg_0=\bfg$.

 Thus we have shown that the kernel $\ker \big(\map_g \to \Sp(2g,\zz_2) \big)$
 is generated by products $T_{\delta_1}^{\epsilon_1} T_{\delta_2}^{\epsilon_2}$ with $\epsilon_i=\pm1$ and
 $[\delta_1] = [\delta_2] \neq0 \in \sfh_1(\Sigma,\zz_2)$.

\item\label{st6b} Multiplying a generator $T_{\delta_1}^{\epsilon_1} T_{\delta_2}^{\epsilon_2}$
 by $T_{\delta_1}^2$ and $T_{\delta_2}^{-2}$ if necessary, we take it to the form
 $T_{\delta_1}\vpht T_{\delta_2}\inv$. Let $\delta_1 =: \alpha$, $\delta_2 =: \alpha'$, so that our
 generator is $T_\alpha \vpht T_{\alpha'}\inv$.

 Write this product in the form $T_{\alpha}\vpht T_\gamma^{-2} T_\gamma^2
 T_{\alpha'}\inv$. Shift $T_\gamma^2$ to the right conjugating $T_{\alpha'}\inv$ and then
 shift $T_\gamma^{-2}$ with conjugation to the right. This transforms $T_{\alpha}\vpht
 T_\gamma^{-2} T_\gamma^2 T_{\alpha'}\inv$ into $T_{\alpha}\vpht T_{\alpha''}\inv T_{\gamma'}^{-2} T_\gamma^2$,
 with $ \alpha''= T_\gamma^2(\alpha')$.  This allows us to change the homology class
 $[\alpha']$ to $[\alpha''] = [\alpha''] + 2(\alpha' \cap \gamma)\cdot [\gamma]$.

 We claim that there exists a sequence of such ``moves'' which transforms the
 integer homology class $[\alpha']$ into the class $[\alpha]$. Clearly, it is sufficient
 to find an inverse transformation of $[\alpha]$ into $[\alpha']$. Fix a curve $\beta$ such
 that $\alpha \cap \beta =1$. Then $[\alpha'] = (2k+1)[\alpha] + 2l[\beta] + 2m[\gamma]$ for some
 non-separating curve $\gamma$ with $\gamma \cap \alpha = \gamma\cap \beta =0$, where both $l$ and $m$
 could be zero. Applying $T_\beta^{\pm2}$ to $\alpha$, we can change $2l$ into $2l \pm2(2k+1)$.
 Iterating this we can transform $l$ into $l'$ with $|l'| \leq |2k+1|$. On the
 other hand, we can replace $\beta$ by a new curve $\beta'$ in the homology class $[\alpha]
 \pm[\beta]$. Then $l$ remains unchanged and $2k+1$ changes to $2k+1 \pm 2l$.
 Consequently, these operations allow us to cancel $l$ out.

 A similar procedure is applied to eliminate $m$. Indeed, for a curve $\gamma'$ in
 the homology class $[\gamma] +[\beta]$ the map $T^{\pm2}_{\gamma'} T^{\mp2}_\gamma$ transforms
 $m$ into $m \pm 2(2k+1)$. Thus we can replace $m$ by $m'$ with $|m'| \leq |2k+1|$.
 However the equality $m' = \pm (2k+1)$ is impossible since $[\alpha']$ is a
 primitive cohomology class. To change $k$, we fix an embedded curve $\delta$ such
 that $\delta \cap \gamma =1$ and $\delta \cap \alpha = \delta \cap \beta =0$ and choose an embedded curve $\delta'$ in
 the class $[\delta] + [\alpha]$. Then $T^{\pm2}_{\delta'}$ transforms the class $(2k+1)[\alpha] +
 2m[\gamma]$ into $(2k+1\pm4m)[\alpha] + 2m([\gamma]\mp2[\delta])$. So we could also make $|2k+1|$
 smaller than $2m$, possibly changing the class $[\gamma]$. The procedure
 terminates at $l=m=0$ and $2k+1 = \pm1$. Since a Dehn twist $T_\delta$ is independent
 of the choice of the orientation on $\delta$, we obtain the equality $[\alpha] = [\alpha']$
 of integral homology classes.

\item\label{st7b} Let $\alpha, \alpha'$ be the simple curves with $[\alpha] =
 [\alpha'] \in \sfh_1(\Sigma, \zz)$ obtained above. It follows that $\alpha' = F_1(\alpha)$ for some $F_1 \in
 \map_g$. As above, let $\beta$ be a curve such that $\alpha \cap \beta =1$
 Set $\beta' := F_1(\beta)$. Then $\beta' \cap \alpha =1$ and
 hence $[\beta'] = [\beta] + l[\alpha] + m[\gamma]$ in $\sfh_1(\Sigma, \zz)$ with some primitive $[\gamma]
 \in \sfh_1(\Sigma, \zz)$ such that $\gamma \cap \alpha = \gamma \cap \beta =0$. Then $T^{-l}_{\alpha'}$ preserves
 the curve $\alpha'$ (up to isotopy) and $\big[T^{-l}_{\alpha'}(\beta') \big] = [\beta] + m[\gamma]$.
 Further, find a curve $\gamma'$ in the homology class $[\alpha] + [\gamma] = [\alpha'] + [\gamma]$
 which is disjoint from $\alpha'$. Then $T^l_{\alpha'} T^{-l}_{\gamma'}$ preserves the curve
 $\alpha'$ (up to isotopy) and transforms the class $[\beta] + m[\gamma]$ into $[\beta]$.

 In this way we have constructed $F_2 \in \map_g$ which takes $\alpha$ to $\alpha'$ and
 preserves the integral class $[\beta]$. Observe that the action of $F_2$ on the
 $\cap$-orthogonal complement to $\zz\langle [\alpha], [\beta]\rangle$ can be realised as a product of
 Dehn twists along curves \emph{disjoint from $\alpha'= F_2(\alpha)$} and from $\beta':= F_2(\beta)$.
 After an appropriate correction of $F_2$ we obtain an $F$ lying in the Torelli
 group $\cali_g := \ker \big(\map_g \to \Sp(2g,\zz) \big)$ such that $\alpha' =
 F(\alpha)$.

\item\label{st12} Now let us apply the explicit description of the Torelli group. It
 is known that $\cali_2$ is generated by Dehn twists $T_\delta$ along
 \emph{separating} curves (\cite{Po}, see also \cite{Jo2}). Every such curve
 cuts $\Sigma$ into two pieces, say $\Sigma'$ and $\Sigma''$, each of them being a surface of
 genus $1$ with one hole. The Dehn twist along $T_\delta$ is given by $(T_\alpha T_\beta)^6$
 for any geometric basis $\alpha$ and $\beta$ of $\Sigma'$, \ie, two curves on $\Sigma'$ meeting
 transversally at a single point. But then $T_\alpha$ and $T_\beta$ are conjugated to the
 Coxeter subsystem $\rma_2:= \{s_1, s_2\}$, and hence $(T_\alpha T_\beta)^6$ is
 conjugated to $\idel^4(\rma_2)$ and is a product of squared quasigenerators
 $\hat t^2$.

 By \cite{Jo1} (see also \cite{Jo2}), the Torelli group  $\cali_g$, $g\geq3$,
 is generated by products $T_\eta\vpht T_{\eta'}\inv$ where $\eta$ and $\eta'$
 are disjoint non-separating curves such that $\eta \sqcup \eta'$ cuts $\Sigma$ into two
 pieces.  Denote these pieces by $\Sigma'$ and $\Sigma''$. Each of them is a surface
 with two boundary circles, and their genera $g'$ and $g''$ are related by $g' +  g'' = g-1$.
 If $g'=0$ or $g''=0$, the curves $\eta, \eta'$ are isotopic and
 the product $T_\eta\vpht T_{\eta'}\inv$ is trivial. It follows that such products
 $T_\eta\vpht T_{\eta'}\inv$ with additional condition $g'=1$ also generate the
 Torelli group $\cali_g$. However, in the case $g'=1$ the piece $\Sigma'$ is a
 chain surface, see  \reffig{fig-chain}. Using the definition of the chain
 relation and the curves on \reffig{fig-chain}, we obtain
\[
T_\eta\vpht T_{\eta'}\inv = T_\eta^2 T_\eta\inv T_{\eta'}\inv =
T_\eta^2 (T_\beta\vpht T_\alpha\vpht T_{\beta'})^4.
\]
Now observe that the chain configuration $\{\beta, \alpha, \beta'\}$
is conjugated to the configuration $\{a_1, a_2, a_3\} =: \rma_3$ on
\reffig{fig-Delta-g}. Consequently, $(T_\beta\vpht T_\alpha\vpht
T_{\beta'})^4$ is conjugated to the squared Garside element
$\idel^2(\rma_3)$ and can be represented as a product of squared
quasigenerators $\hat t^2$.
\qed
\end{stepping}

\newcorol{cor-IW}
\sli The kernel $\ker\big( \map_{g,1} \to \Sp(2g, \zz_2) \big)$
is generated by the squares $T^2_\delta$ of Dehn twists along non-separating curves
$\delta \subset \Sigma\bs \{ z_0 \}$.

\slii The Weyl--Torelli group $\IW_g = \ker\big( \sfw(\cals_g) \to \Sp(2g, \zz_2)
\big)$ is generated by elements conjugated to $w_\circ(\rme_7)$.
\end{corol}

\smallskip
Recall that \propo{prop-K-van} reduces the topological statement about the (non-)existence
of special embeddings of the Klein bottle into topological Lefschetz
fibrations  to certain algebraic relations in the mapping class group $\map_g$.
It turns out that these relations hold in $\Sp(2g,\zz_2)$. Lifting the
corresponding elements to the braid group $\br(\cals_g)$, we can project them
to $\zz\langle \calh_g\rangle$. So we need tools to distinguish the combinatorial
structure of $\calh_g$ and $\zz\langle \calh_g\rangle$.

\smallskip
Let us describe the first such tool. Let $R$ be any commutative ring and $H$
a free $R$-module of finite rank. Denote by $\sftt^\bullet H$ the tensor algebra of
$H$ over $R$. For each degree $d$, let $\sftt^d_\sym H \subset \sftt^dH$ be the
submodule consisting of tensors invariant with respect to the natural action of the
symmetric group $\sym_d$ permuting the tensor factors $H \otimes \cdots \otimes H$. Define the
\slsf{shuffle product} in $\sftt^\bullet H$ as follows. For $A \in \sftt^k H$ and $B \in
\sftt^lH$ set $A \bullet B := \sum_{\sigma \in \sym_{k+l}}' \sigma(A \otimes B)$ where the sum $\sum'$ is
taken over all permutations $\sigma \in \sym_{k+l}$ of tensor factors which preserve
the order of the first $k$ and of the last $l$ factors.  In other words, $\sigma \in
\sym_{k+l}$ must satisfy the condition $\sigma(i) <\sigma(j)$ if $i<j\leq k$ and if $k<
i<j$. In particular, $v \bullet w = v\otimes w + w\otimes v$ and $u \bullet(v\otimes w) = u\otimes v\otimes w + v\otimes u \otimes w
+ v\otimes w\otimes u$ for $u, v,w \in H$.  One verifies immediately the following
properties of the introduced structures:
\begin{itemize}
\item[--] \ $v^{\otimes d} \in \sftt^d_\sym H$ for any $v\in H$; \ \ $A \bullet B \in \sftt^{k+l}_\sym H$
 if $A \in \sftt^k_\sym H$ and $B \in \sftt^l_\sym H$;
\item[--] the $\bullet\,$-product satisfies the associativity, commutativity, and
 distributivity laws;
\item[--]  \ for $v,w \in H$ one has the \slsf{binomial formula\/}:
\begin{equation}\eqqno(bull-bin)
(v+w)^{\otimes d} =
v^{\otimes d} + v^{\otimes(d-1)} \bullet w + v^{\otimes(d-2)} \bullet w^{\otimes2}+\cdots + w^{\otimes d};
\end{equation}
\item[--] \ $v^{\otimes k} \bullet v^{\otimes l}= \binom{k+l}{k} v^{\otimes(k+l)}$ for $v\in H$.
\end{itemize}
The last property shows that in the case when $H\cong R^{\oplus r}$ is the free module of rank
$r$, the algebra $(\sftt^\bullet_\sym H, \bullet)$ is isomorphic to the $r$-th tensor
power $\mathbf{A}^\bullet \otimes \cdots \otimes \mathbf{A}^\bullet$ of the so-called \slsf{algebra of
divided powers} $\mathbf{A}^\bullet$.

In our case, we set $R: =\zz_2$ and $H: = \sfh_1(\Sigma, \zz_2)$.  Let us define
homomorphisms $\wp^d: \zz_2\langle\calh_g\rangle \to \sftt^d H$ by setting $\wp^d(A_v) := v^{\otimes d}$
for each $v\in \calh_g$ and use the same notation $\wp^d: \sfp(\cals_g) \to \sftt^d H$ for the composition.
Clearly, each $\wp^d$ takes values in $\sftt^d_\sym H \cong \sym^d(H)$. The meaning of $\wp^1$ is obvious:
it maps each $A_v$ to the vector $v \in \sfh_1(\Sigma, \zz_2)$.  To describe $\wp^2$, we
observe that the space $\sym^2(H)$ is naturally isomorphic to the $\zz_2$-Lie
algebra $\sp(2g, \zz_2)$. Explicitly, the isomorphism is given by $M \in
\sym^d(H) \leftrightarrow J \cdot M \in \sp(2g,\zz_2)$, where $J$ is the ``symplectic'' matrix $J
:= \spmatr{ 0 &\id \\ \id & 0}$ with square $g\times g$ blocks, and
$\sp(2g,\zz_2)$ is realised as a matrix Lie subalgebra of $\mat(2g, \zz_2)$
with respect to a symplectic basis $\alpha_1,\ldots,\alpha_g; \beta_1, \ldots, \beta_g$ of $H$. Further,
$\sp(2g,\zz_2)$ admits a natural group extension
\[
1 \to \sp(2g,\zz_2) \to \Sp(2g,\zz_4)\to \Sp(2g,\zz_2)\to 1,
\]
and an explicit calculation with any $T^2_\delta$ shows that the homomorphism
$J\cdot\wp^2 : \sfp(\cals_g) \to \sp(2g,\zz_2) \hrar \Sp(2g,\zz_4)$ coincides with
the composition $\sfp(\cals_g) \hrar \br(\cals_g) \to \map_g \to \Sp(2g,\zz)\to
\Sp(2g,\zz_4)$. Observe also that the matrix $J$ can be given by $ \sum_{i=1}^g
\alpha_i \bullet \beta_i $ for any symplectic basis $\alpha_1,\ldots,\alpha_g; \beta_1, \ldots, \beta_g$ of $H$.
Besides, recall the definition of the hyperelliptic relation element $[T_{\beta_g},
\idel^2(\rma_{2g})]$ given in \refthm{thm-matsu}.

\newlemma{lem-theta3} \begin{itemize}
\item[\slip] $\wp^3\big( \idel^2(\rma_5) \big) = \wp^3\big(
 \idel^2(\rme_7) \big) =0$.
\item[\sliip] $\wp^3\big( [T_{\beta_g}, \idel^2(\rma_{2g})] \big) = \beta_g \bullet \left(
  \sum_{i=1}^g \alpha_i \bullet \beta_i \right) $.
\end{itemize}
\end{lem}

\proof \sli Note that $\calt(\rma_5) = \calh_2$ and
$\calt(\rme_7) = \calh_3$. So in both cases $\idel^2(\ldots)$ is the sum $\sum _{v \neq0
 \in V} A_v$ where the subspace $V \subset \sfh_1(\Sigma, \zz_2)$ is of dimension $4$ and $6$,
respectively. Consequently, $\wp^3\big( \idel^2(\ldots) \big)$ is the sum $\sum _{v
 \in V} v^{\otimes3}$. Fixing a basis $e_1, \ldots, e_r$ of $V$ ($r=4$ or $6$ according to
the case) and expanding the sum $\sum _{v \in V} v^{\otimes3}$ using the binomial formula
\eqqref(bull-bin) we obtain the sum of monomials $e^{\mathbf{n}} := e_1^{\otimes
 n_1} \bullet\cdots \bullet e_r^{\otimes n_r}$ with $n_1+\cdots +n_r=3$. Since at least one $n_i$
vanishes, each monomial $e^{\mathbf{n}}$ appears an even number of times, and the
sum vanishes.

\medskip\noindent
\slii For simplicity, we use the same notation $a_i, \alpha_j, \beta_k \in \calh_g$ for the
$\zz_2$-homology classes of the corresponding curves. Set $v_{ij} :=
\sum_{k=i}^{j} a_k$ and $v_i := v_{i,2g}$. It follows that the projection of the
squared Garside element $\idel^2(\rma_{2g})$ to $\zz_2\langle\calh_g\rangle$ is the sum
$\sum_{1\leq i\leq j\leq2g} A_{v_{ij}}$. Since $\beta_g$ is disjoint from each $a_j$ with
$j<2g$, the hyperelliptic relation element $[T_{\beta_g}, \idel^2(\rma_{2g})]$
projects to the sum $\sum_{i=1}^{2g} \big( A_{v_i} + A_{v_i+ \beta_g} \big)$. Hence
\[\textstyle
\wp^3\big( [T_{\beta_g}, \idel^2(\rma_{2g})] \big) =
\sum_{i=1}^{2g} \big( v_i^{\otimes3} + (v_i+\beta_g)^{\otimes3} \big).
\]
Expanding and using the identity $2=0 \in \zz_2$ we obtain
\[\textstyle
\wp^3\big( [T_{\beta_g}, \idel^2(\rma_{2g})] \big) =
\sum_{i=1}^{2g} \big( v_i \bullet \beta_g^{\otimes2}  +  v_i^{\otimes2} \bullet \beta_g \big).
\]
Plugging $v_i = a_i + a_{i+1}+\cdots + a_{2g}$ into the sum $\sum_{i=1}^{2g} v_i$
we see that every $a_j$ appears $j$ times. Thus
\[\textstyle
\sum_{i=1}^{2g} v_i = a_1 +a_3+\cdots
+ a_{2g-1} = \beta_1 +(\beta_1+\beta_2) +\cdots + (\beta_{g-1} + \beta_g) = \beta_g.
\]
Similarly, plugging $v_i^{\otimes2} = \sum_{j=i}^{2g} a_j^{\otimes2} + \sum_{i\leq j <k \leq2g} a_j \bullet
a_k$ into the sum $\sum_{i=1}^{2g} v_i^{\otimes2}$ we obtain
\[\textstyle
\sum_{i=1}^{2g} v_i^{\otimes2} = \sum_{j=1}^{2g} j \cdot a_j^{\otimes2}
+ \sum_{1\leq j <k \leq2g} j\cdot a_j \bullet a_k.
\]
The first sum gives
\begin{align*}
a_1^{\otimes2} + a_3^{\otimes2} + \cdots + a_{2g-1}^{\otimes2} &=
\beta_1^{\otimes2} + (\beta_1 + \beta_2)^{\otimes2} + \cdots + (\beta_{g-1} + \beta_g)^{\otimes2}\\
& = \beta_1 \bullet \beta_2 + \beta_2 \bullet \beta_3 + \cdots + \beta_{g-1} \bullet \beta_g + \beta_g^{\otimes2}.
\end{align*}
The second sum can be rewritten as $\sum_{1\leq j <k \leq2g} j\cdot a_j \bullet a_k = \sum_{j=1}^g
\sum_{k=2j}^{2g} a_{2j-1} \bullet a_k$. The interior summation yields $\sum_{k=2j}^{2g} a_k
= (\alpha_j +\ldots+ \alpha_g) + (\beta_j + \beta_g)$. So, setting $\beta_0 :=0$, we can expand
the second sum as
\[\textstyle
\sum_{1\leq j <k \leq2g} j\cdot a_j \bullet a_k =  \sum_{j=1}^g (\beta_{j-1}  + \beta_j)\bullet
\big (  (\alpha_j +\ldots+ \alpha_g) + (\beta_j + \beta_g)\big).
\]
Let us analyse the terms in the last sum.
Firstly, $\sum_{j=1}^g (\beta_{j-1} + \beta_j)\bullet \beta_g = \beta_g\bullet \beta_g =2 \beta_g^{\otimes2} =0$.
Similarly, $\beta_j\bullet \beta_j=0$ and $\sum_{j=1}^g (\beta_{j-1} +
\beta_j)\bullet \beta_j = \beta_1 \bullet \beta_2 + \beta_2 \bullet \beta_3 + \cdots + \beta_{g-1} \bullet \beta_g$.
The terms containing a given $\alpha_k$ are
\[\textstyle
\sum_{j=1}^k (\beta_{j-1} + \beta_j)\bullet \alpha_k = (\beta_0 + \beta_k)\bullet\alpha_k = \beta_k\bullet\alpha_k.
\]
So finally we can compute
\[\textstyle
\wp^3\big( [T_{\beta_g}, \idel^2(\rma_{2g})] \big) =
\beta_g \bullet \beta_g^{\otimes2} +(\beta_g^{\otimes2} +  \sum_{j=1}^g \beta_j \bullet \alpha_j)  \bullet \beta_g =
\beta_g \bullet  \sum_{j=1}^g \alpha_j \bullet \beta_j,
\]
as desired.
\qed

\newcorol{cor-p3} For $g\geq2$ the homomorphism $\wp^3$ extends to the group
$\ker \big( \map_{g,1} \to \Sp(2g,\zz_2) \big)$ so that the
composition $\wp^3 \circ \dd_\pi : \pi_1(\Sigma) \to \sftt^3$ can be realised as the composition
of the projection $\pi_1(\Sigma) \to \sfh_1(\Sigma, \zz_2)$ with the embedding $\sfh_1(\Sigma,
\zz_2) \hrar \sftt^3$ given by $\gamma \mapsto \gamma \bullet \sum_{i=1}^g \alpha_i\bullet \beta_i$.
\end{corol}

\proof Apply \refcorol{cor-IW}, \slip. \qed

\newprop{prop-IW-ab}
\sli The abelianisation $\IW_{g,\ab}$ of the
Weyl--Torelli group $\IW_g$ is a $\zz_2$-vector space naturally
isomorphic to $\land^6 \sfh_1(\Sigma, \zz_2)$.

\slii  There exists a natural lattice extension
\[
0 \to \zz\langle \calh_g \rangle \to \Lambda_g \to \IW_{g,\ab} \to 0
\]
such that $\Lambda_g$ can be realised as a sublattice in $\frac12 \zz\langle \calh_g \rangle$
and the image of the induced embedding $\IW_{g,\ab} \subset \frac12 \zz\langle
\calh_g \rangle /\zz\langle \calh_g \rangle \cong \zz_2\langle \calh_g \rangle$ is generated by the sums
$L_V :=\sum_{v \neq 0 \in V} A_v$ in which $V$ is a symplectic $6$-dimensional subspace of\/
$\sfh_1(\Sigma, \zz_2)$.

\sliii There exists a group extension
\begin{equation}\eqqno(latt-ext)
0 \to \Lambda_g \to \wh\Sp(2g, \zz_2)  \to \Sp(2g, \zz_2) \to 1
\end{equation}
and a homomorphism of the extension $1 \to \sfp(\cals_g) \to \br(\cals_g) \to
\sfw(\cals_g) \to1$ onto this extension such that the homomorphism
$\sfw(\cals_g) \to \Sp(2g, \zz_2)$ has the usual meaning and $\sfp(\cals_g)
\to\Lambda_g$ is the composition $\sfp(\cals_g) \to \sfp_\ab(\cals_g) \to \zz\langle \calh_g \rangle \subset
\Lambda_g $.
\end{prop}

A subspace $V \subset \sfh_1(\Sigma, \zz_2)$ is called symplectic if the restriction of the
intersection form to $V$ is non-degenerate.

\proof Let $L_0:= C\cdot \idel^{-2}(\rme_6) \idel(\rme_7)$ and $L'_0:=
\idel^{-2}(\rme_6) \idel(\rme_7)$ be the ``basic'' non-modified and
modified lantern relation elements. Then their projections to $\sfw(\cals_g)$
are equal and give the longest element $w_\circ(\rme_7)$ of the Coxeter subgroup
$\sfw(\rme_7) \subset\sfw(\cals_g)$. Since $w_\circ(\rme_7)$ has order $2$ and its
conjugates generate $\IW_g$, the abelianisation $\IW_{g,\ab}$ is a
$\zz_2$-vector space.

Now consider the natural extension
\[
0 \to \zz\langle \calh_g \rangle \to Q \to \IW_g \to 1
\]
so that $Q$ is the quotient of the kernel $\ker\big( \br(\cals_g) \to \Sp(2g,
\zz_2) \big)$ by the kernel\break $\ker\big( \sfp(\cals_g) \to \zz\langle \calh_g \rangle\big)$.
Every element in $Q$ is the product of conjugates $x*L_0$ with $x\in
\br(\cals_g)$ and elements of $\zz\langle \calh_g \rangle$. Since $\zz\langle \calh_g \rangle$ lies in
the centre of $Q$ and the chain relation $C = L_0 (L'_0)\inv$ lies in $\zz\langle
\calh_g \rangle$, the commutators $[x*L_0, y*L_0]$ and $[x*L'_0, y*L_0]$ (and so on)
are equal and generate the commutator group $Q' := [Q,Q]$.

We claim that the intersection $Q' \bigcap \zz\langle \calh_g \rangle$ is trivial.
Suppose that $z\in Q' \bigcap \zz\langle \calh_g \rangle$. The condition $z\in Q'$ means that $z$
is represented by an element $\hat z= \prod_j [x_j * L_0, y_j * L_0]$ from
$\br(\cals_g)$. Since $x_j * L_0$ and $y_j * L_0$ project to $\id \in \Sp(2g,\zz_2)$,
we can apply to them the algorithm from \propo{prop-map-gen}. However, we use only the steps
that do not involve the lantern relation, i.\,e.\ omit \slsf{Steps \ref{st7}--\ref{st10}}.
As a result, each $x_j * L_0$ and $y_j * L_0$ is represented as a product of elements
from $\zz\langle \calh_g \rangle$ and Dehn twists $T_{\delta_i}$ such that the $\zz_2$-homology
classes $[\delta_i]$ are linear combinations of $[\beta_1],\ldots ,[\beta_g]$. Since $\zz\langle
\calh_g \rangle$ lies in the centre, it follows that the commutators $[x_j * L_0,
y_j * L_0]$ are products of Dehn twists $T_{\delta_i}$ with $[\delta_i] \in \zz_2\langle[\beta_1],\ldots
,[\beta_g] \rangle$.

Now consider $z$ as an element of $\sfp(\cals_g)$. Write the projection of $z$ to
$\zz\langle \calh\rangle$ in the form $z= \sum_{v\in \calh_g} n_v A_v$ with $n_v \in \zz$.  Let
$G$ be the subgroup of $\sfw(\cals_g)$ obtained from $\IW_g$ by adding
quasireflections $t\in \calt(\cals_g)$ whose projections to $\Sp(2g,\zz_2)$ are
Dehn twists $T_{\delta_i}$ with $[\delta_i] \in \zz_2\langle[\beta_1],\ldots ,[\beta_g] \rangle$.
Since $G$ contains $\IW_g$, there exists a $G$-equivariant projection $\calg_g =
\calt(\cals_g)/ \IW_g \to \calt(\cals_g)/G$ such that coset classes in
$\calt(\cals_g)/G$ can be described as $G$-orbits in $\calg_g$. Clearly,
each $[\delta] \in \zz_2\langle[\beta_1],\ldots ,[\beta_g] \rangle$ is fixed by $G$. The remaining orbits are
as follows. Set $V:= \zz_2\langle[\beta_1],\ldots ,[\beta_g] \rangle$ and $U := \zz_2\langle[\alpha_1],\ldots ,[\alpha_g]
\rangle$. Then for every $u \in \calh_g$ with non-trivial projection to $U$ parallel
to $V$ the $G$-orbit $G \cdot u$ consists of all vectors $u+v$ with $v$ varying
over $V$. In particular, $\sum_{v\in V} n_{u+v} =0$ by \refthm{thm-subgr}.  Set
$V_1 := \zz_2\langle[\beta_2],\ldots ,[\beta_g] \rangle$. Then for the special case $u=[\alpha_1]$ we obtain
$\sum_{v\in V_1} \left( n_{[\alpha_1] +v} + n_{[\alpha_1 +\beta_1]+v} \right) =0$. The
possibility to choose an arbitrary $\zz_2$-symplectic basis provides two
further relations:
\[\textstyle
\sum_{v\in V_1} \left( n_{[\beta_1] +v} + n_{[\alpha_1 + \beta_1]+v} \right) =0
\text{\ \ and\ \ }
\sum_{v\in V_1} \left( n_{[\alpha_1] +v} + n_{[\beta_1]+v} \right) =0.
\]
Comparing them, we obtain the relation
\[\textstyle
\sum_{v\in V_1} n_{w +v}   =0
\]
first for $w= [\alpha_1], [\beta_1], [\alpha_1 + \beta_1]$, and then for all $w \in \calh_g$ which
are $\zz_2$-orthogonal to $\beta_2,\ldots, \beta_g$. Now set $V_k := \zz_2\langle[\beta_{k+1}],\ldots
,[\beta_g] \rangle$ and write the latter equality for $w=[\beta_1]$ in the form
\[\textstyle
\sum_{v\in V_2} \left( n_{[\beta_1] +v } + n_{[ \beta_1+ \beta_2]+v} \right) =0.
\]
The same argument as above provides the relation $\sum_{v\in V_2} n_{w +v}   =0$,
again first for $w= [\beta_1]$ and then for all $w \in \calh_g$ that
are $\zz_2$-orthogonal to $V_2$. Repeating this argument, we obtain the
relation $n_w =0$ for all $w \in \calh_g$. This means the desired triviality of
$z$ in $\zz\langle \calh_g \rangle$.

The property $Q' \bigcap \zz\langle \calh_g \rangle =1$ implies that the abelianisation $\Lambda :=
Q_\ab = Q / Q'$ includes into the extension \eqqref(latt-ext). Since the chain
relation element $C$ and $\idel^2(\rme_6)$ lie in $\zz\langle \calh_g \rangle$, in the
description of the extension \eqqref(latt-ext) we can replace $L_0$ by
$\idel(\rme_7)$. The square $\idel^2(\rme_7)$ is the sum $\sum_{t\in \calt(\rme_7)}
A_t$. The quasireflection set $\calt(\rme_7)$ coincides with the set $\calh_3$
of non-zero vectors in $\zz_2\langle[\alpha_1],\ldots ,[\beta_3] \rangle$. Conjugating by $\Sp(2g,
\zz_2)$ we obtain the description of generators of $\IW_{g,\ab} = \Lambda_g /
\zz\langle \calh_g \rangle$ claimed in the lemma.

To obtain the isomorphism $\IW_{g,\ab}\cong \land^6 \sfh_1(\Sigma, \zz_2)$, we
compute $\wp^6 \big(\idel^2(\rme_7) \big)$. An explicit calculation using
the binomial formula \eqqref(bull-bin) shows that $\wp^6
\big(\idel^2(\rme_7) \big) = \alpha_1 \bullet \alpha_2 \bullet \alpha_3 \bullet \beta_1 \bullet \beta_2 \bullet \beta_3$.
Further, note that over the coefficient ring $\zz_2$ we have the identity $v \bullet v =0$. It
follows that the algebra generated by $\sfh_1(\Sigma, \zz_2)$ and by the shuffle
product is simply the $\zz_2$-Grassmann algebra over $\sfh_1(\Sigma, \zz_2)$.
Consequently, the shuffle product $\alpha_1 \bullet \alpha_2 \bullet \alpha_3 \bullet \beta_1 \bullet \beta_2 \bullet \beta_3$ can be
identified with the wedge product $\alpha_1 \land \alpha_2 \land \alpha_3 \land \beta_1 \land \beta_2 \land \beta_3 \in \land^6
\sfh_1(\Sigma,\zz_2)$. Vice versa, any symplectic vector subspace $V \subset \sfh_1(\Sigma,
\zz_2)$ of dimension $6$ is of the form $V = x\cdot \zz_2\langle[\alpha_1],\ldots ,[\beta_3] \rangle$
for some $x \in \Sp(2g, \zz_2)$. For such $V$, the conjugate $x * \idel^2(\rme_7)$ equals the sum
$\sum _{v \neq0 \in V} A_v =: L_V \in \zz_2\langle \calh_g \rangle$, and $\wp^6 \big(L_V)$ is the wedge
product of the vectors in any basis of $V$.  So the isomorphism $\IW_{g,\ab}\cong \land^6 \sfh_1(\Sigma,
\zz_2)$ follows from the fact that the elements $\Sp(2g, \zz_2)$-conjugate to $\alpha_1
\land \alpha_2 \land \alpha_3 \land \beta_1 \land \beta_2 \land \beta_3$ span the space $\land^6 \sfh_1(\Sigma,\zz_2)$.
\qed

\medskip
Now we introduce the second tool used to explore the structure of the group
$\zz\langle \calh_g\rangle$. To any $\varphi_\mu=(\mu\cap\cdot):\sfh_1(\Sigma, \zz_2) \to \zz_2$
we associate a homomorphism $\hat\varphi_\mu: \zz\langle \calh_g\rangle \to \zz$ by setting
$\hat\varphi_\mu(A_v) := 1\in \zz$ if $\varphi_\mu(v)\equiv 1\in \zz_2$ and
$\hat\varphi_\mu(A_v) := 0$ if $\varphi_\mu(v)\equiv 0$ for any element $v\in \calh_g \subset \sfh_1(\Sigma, \zz_2)$.
Thus $\hat\varphi_\mu$ counts the algebraic number of generators $A_v$ with $\varphi_\mu(v)\equiv 1$. Extend this
homomorphism to a homomorphism $\hat\varphi_\mu: \sfp(\cals_g) \to \zz$ in the obvious
way.

\newlemma{lem-chain-eval} Let $\hat \varphi_\mu: \zz\langle\calh_g\rangle \to \zz$ be the
homomorphism induced by $\varphi_\mu: \sfh_1(\Sigma, \zz_2) \to \zz_2$.

\sli For any chain relation element $C$ the value $\hat \varphi_\mu(C)$ is $0$ or $-4$.

\slii For any squared modified lantern relation element $L^2$ the value $\hat
\varphi_\mu(L^2)$ is $0$ or $-8$.

\sliii For any element $x$ representing a Dehn twist $T_\delta$ along a
\emph{separating} curve $\delta$ one has $\hat \varphi_\mu(x) \equiv0 \mod 4$.

\sliv $\hat \varphi_\mu\big( [T_{\beta_g} , \idel^2(\rma_{2g}) ] \big) \equiv (2g-2)\varphi_\mu(\beta_g)
\mod 4$.
\end{lem}

\newcorol{cor-lant-ext} The homomorphism $\hat \varphi_\mu: \zz\langle\calh_g\rangle \to \zz$
extends naturally to a homomorphism  between $\zz$ and the stabiliser group
$\br(\cals_g)_\mu$ of the element $\mu$. The latter induces a homomorphism $\hat
\varphi_\mu: \map_{g,1,\mu} \to \zz_4$.
\end{corol}

\proof First we show the claims of the corollary. It follow immediately from
part \sliii of the lemma that there exists a natural extension $\hat \varphi_\mu: \Lambda_g
\to \zz$, such that $\hat \varphi_\mu(L) = 0$ or $-4$ for any \emph{modified} lantern
relation element $L$. To extend $\hat \varphi_\mu$ onto $\br(\cals_g)_\mu$, we observe
that the group $\Sp(2g, \zz_2)_\mu$ is generated by (the images of) Dehn twists
$T_\delta$ with $\delta \cap \mu \equiv0 \mod2$. Consequently, $\br(\cals_g)_\mu$ is generated by
such Dehn twists, $\sfp(\cals_g)$, and lantern relation elements. We set $\hat
\varphi_\mu(T_\delta) = 0$ for such Dehn twists. To see that this definition makes sense,
consider an element $x$ which is a product of $T_{\delta_j}$ with $\delta_j \cap \mu \equiv0 \mod2$,
and which lies in the subgroup generated by $\sfp(\cals_g)$ and lantern relation
elements. Then $x$ can be projected to $\Lambda_g$, and hence its square $x^2$ lies
in $\sfp(\cals_g)$. Now using \refthm{thm-subgr}, we obtain $\hat \varphi_\mu(x^2) = 0$.
This implies that $\hat \varphi_\mu(x) = 0$ and shows that the homomorphism
$\hat \varphi_\mu: \br(\cals_g)_\mu \to \zz$ is indeed well-defined. The existence
of the induced homomorphism $\hat \varphi_\mu: \map_{g,1,\mu} \to \zz_4$
is trivial (modulo the lemma, of course).

\medskip%
Now we prove the assertions of the lemma.

\smallskip\noindent%
\sli Obviously, we have the following \slsf{conjugation property}:
$\hat\varphi_{x* \mu}(x*A) = \hat\varphi_\mu(A)$ for any $x \in \Sp(2g,\zz_2)$ and any $A \in
\zz\langle\calh_g\rangle$. Hence we may assume that $C$ is the basic chain relation
element $C_0 :=\idel^{-4}(\rma_4) \idel^2(\rma_5)$. Further, restricting $\varphi_\mu$
to $\zz_2\langle \alpha_1, \alpha_2, \beta_1, \beta_2\rangle$ we can also suppose that $g=2$ and that this
restriction is a non-trivial homomorphism.

We make use of the standard facts on root systems, see \cite{Bou}.  In
particular, we have the natural identifications $\cals_2 = \rma_5$ of the
Coxeter systems and $\Sp(4,\zz_2) = \sfw(\rma_5)$ of the groups. Besides, the
set $\calh_2$ is naturally identified with the set $\Phi^+(\rma_5)$ of
\slsf{positive roots} of the system $\rma_5$.

We claim that there are exactly two orbits of the action of $\sfw(\rma_4)$ on
$\Phi^+(\rma_5)$.  The first one is clearly the root system $\Phi^+(\rma_4)$.
The explicit description of the root systems gives the five remaining elements: $\{
v_1, v_1+v_2, v_1+v_2+v_3, \ldots, v_1+\cdots +v_4+ v_0\}$. (Here we are using the
notation from \refdefi{def-S-g} and denote by $v_i$ the simple root
corresponding to $s_i$, $i=1,\ldots,4,0$.) We see explicitly that $v_1+v_2$ is
obtained from $v_1$ with the help of the reflection $s_2 = T_{a_2}$,
$v_1+v_2+v_3$ from $v_1+v_2$ with the help of $s_3 = T_{a_3}$ and so on.  This
gives us the desired description of the orbits.

It follows that it is sufficient to calculate $\hat\varphi_\mu\big(\idel^{-4}(\rma_4)
\idel^2(\rma_5) \big)$ for the cases in which $\mu=\beta_1=a_1$ and $\mu=\beta_2=a_0$. By definition,
$\hat\varphi_\mu\big( \idel^2(\rma_5) \big)$ is the number of vectors $v \in \sfh_1(\Sigma,
\zz_2)$ satisfying $\varphi_\mu(v)=1$. So there are $2^3=8$ such vectors, and
$\hat\varphi_\mu\big(\idel^2(\rma_5) \big) =8$.

To simplify notation, we replace the system $\{s_2,s_3,s_4,s_0\}$ by $\{s_1,
s_2,s_3,s_4\}$. This does not change anything in the case $g=2$ because of the
existence of an appropriate isomorphism. Then $\idel^2(\rma_4) = \sum_{1\leq j\leq k \leq4}
A_{v_{jk}}$ with $v_{jk} := \sum_{i=j}^k a_i \in \sfh_1(\Sigma, \zz_2)$. Then
$\hat\varphi_{\beta_1}\big(\idel^2(\rma_4) \big)$ is the number of vectors $v_{jk}$
``containing'' the summand $\alpha_1=a_2$. For such vectors we must have $1\leq j \leq2
\leq k\leq4$, so that $\hat\varphi_{\beta_1}\big(\idel^2(\rma_4) \big) =2\times3=6$ and
$\hat\varphi_{\beta_1}\big(\idel^4(\rma_4) \big) =12$. A similar consideration for
$\mu=\beta_2$ leads to the vectors $v_{jk}$ with $1\leq j \leq k=4$, and hence
$\hat\varphi_{\beta_2}\big(\idel^4(\rma_4) \big) =4\times2 =8$. So $\hat\varphi_\mu(C)$ equals
either $8-8=0$ or $8-12=-4$.

\medskip\noindent
\slii Now we consider the squared modified lantern relations $L^2$. We claim that
after replacing $\rma_4$ by $\rme_6$ and $\rma_5$ by $\rme_7$, the argumentat
remains essentially the same. In particular, we have the natural identifications $\cals_3
= \rme_7$ of the Coxeter systems and $\Phi^+(\rme_7) = \calh_3$ of positive
roots, and the natural isomorphism $\sfw(\rme_7) =\Sp(6,\zz_2) \times \zz_2\langle
w_\circ(\rme_7)\rangle$ of the groups.  After conjugation we come to the case of the
basic lantern relation $\idel^{-4}(\rme_6) \idel^2(\rme_7)$.  Similar to the
previous situation, we may assume that $g=3$ and that the restricted $\varphi_\mu$ is
non-trivial. As above, $\hat\varphi_\mu\big( \idel^2(\rme_7) \big)$ is the number of
vectors $v \in \sfh_1(\Sigma, \zz_2)$ satisfying $\varphi_\mu(v)=1$.  This time this number
is $2^5=32 =\hat\varphi_\mu\big( \idel^2(\rme_7) \big)$.

Now let us describe the $\sfw(\rme_6)$-orbits in $\calt(\rme_7)$. Since the
quasireflections are in bijection with positive roots of the corresponding Lie
algebra, we have $63$ quasireflections in $\calt(\rme_7)$ and $36$ in
$\calt(\rme_6)$. Next, the embedding $\rme_6 \subset \rme_7$ of the Coxeter systems
induces an embedding of the set of quasireflections. Thus $\calt(\rme_6)$ is a
$\sfw(\rme_6)$-orbit in $\calt(\rme_6)$, and there are $63-36 =27$ remaining
quasireflections.

One of these remaining quasireflections is the element $s_1$. (Here we
are using the notation from \refdefi{def-S-g}.) To determine its
$\sfw(\rme_6)$-orbit let us consider the complex Lie algebra $\frg$ of type
$\rme_7$. Fix a Cartan subalgebra $\frh$ in $\frg$ and a compatible system of
simple roots, the latter being in natural bijection with the system $\rme_7=\{
s_0;s_1,\ldots,s_6\}$.  For each $t\in \calt(\rme_7)$ denote by $\alpha_t \in \frh^*$ the
corresponding positive root, and by $\frg^\pm_t$ the root subspace corresponding
to $\pm\alpha_t$. Further, let $\frg' \subset \frg$ be the Lie subalgebra of type $\rme_6$
defined by the embedding $\rme_6 \subset \rme_7$ and $\frh'\subset \frh$ the compatible
embedding of the Cartan subalgebra. Denote by $B : \frh^* \times \frh^*
\to\cc$ the \slsf{canonical bilinear form} of $\frh^*$, normalised by the
condition $B(\alpha_{s_i}, \alpha_{s_i}) = +2$.

Denote by $V^+$ (respectively, by $V^-$) the sum of root spaces $\frg_t^+$ (respectively,
$\frg_t^-$) over quasireflections $t \in \calt(\rme_7) \bs \calt(\rme_6)$. Since
for every root $\alpha_t$ with $t \in \calt(\rme_7) \bs \calt(\rme_6)$ the
coefficient of $\alpha_{s_1}$ is positive, $V^+$ is invariant with respect to the adjoint
action of $\frg'$. This means that the entire orbit $\sfw(\rme_6)\cdot \alpha_{s_1}$ lies
in $\Phi^+(\rme_7)$ (this is different from the action of the full Weyl group
$\sfw(\rme_7)$ which inverts every root).

From the Dynkin diagram we see that $B(\alpha_{s_1}, \alpha_{s_2}) =-1$ and $\alpha_{s_1}$ is
orthogonal to the remaining simple roots $\alpha_{s_0}; \alpha_{s_3},\ldots, \alpha_{s_6}$ in
$\rme_6$. This means that for every $v \in \frh'$ we have $\alpha_{s_1}(v)
= -\omega_{s_2}(v)$, where $\omega_{s_2}$ is the \slsf{fundamental weight} of the system
$\rme_6$ dual to $\alpha_{s_2}$.  Consequently, $V^+$ contains an irreducible
$\frg'$-submodule with the minimal weight $-\omega_{s_2}$.  It is known that the
dimension of this submodule is $27$. Thus $\dim V^+ \geq27$. By the same
argument, $\dim V^- \geq27$. Comparing dimensions, we conclude that the
$\frg'$-irreducible decomposition of $\frg$ is $\frg = \frg' \oplus V^+ \oplus V^- \oplus
\cc\langle \alpha^\lor_{s_1} \rangle$, where $\alpha^\lor_{s_1}$ is the coroot dual to $\alpha_{s_1}$. Now
observe that by \cite{Bou}, Ch.VI, {\S}1, \slsf{{\'E}x{\'e}rcise 23}, the weight
$-\omega_{s_2}$ is \slsf{minuscule}.%
\footnote{I would like to thank W.~S{\"o}rgel for this reference.}
This property is equivalent to the assertion that the
$\frg'$-weights of $V^+$ form a single $\sfw(\rme_6)$-orbit.

\medskip
Let us turn back to the calculation of the possible values of $\hat \varphi_\mu(\idel^2 (\rme_6))$.
For this purpose we use the following explicit construction of the
root system of type $\rme_6$. Let $\frh'_\rr{}^*$ (the real form of the dual
Cartan subalgebra $\frh'{}^*$) be the space spanned by vectors $\varepsilon_1,\ldots,\varepsilon_6;\varepsilon$
satisfying the linear relation $\varepsilon_1+ \cdots+\varepsilon_6 =0$. Define the bilinear form on
$\frh'_\rr{}^*$ by setting $B(\varepsilon_i, \varepsilon_i) = \frac56$, $B(\varepsilon_i, \varepsilon_j) = -\frac16$,
$B(\varepsilon, \varepsilon) = \frac12$, and $B(\varepsilon_i, \varepsilon) = 0$. Then the set
\[
\Phi := \{\; \pm2\varepsilon,\  \varepsilon_i + \varepsilon_j + \varepsilon_k \pm\varepsilon\;(i>j>k);\  \ \varepsilon_i -\varepsilon_j \,(i\neq j)\; \}
\]
is a root system of type $\rme_6$, and the set
\[
\bfpi := \{\; \alpha_{s_i} := \varepsilon_i -\varepsilon_{i-1}\;(i=2,\ldots,6);\   \ \alpha_{s_0} := \varepsilon_1 + \varepsilon_2 + \varepsilon_3 +\varepsilon\;\}
\]
is the system of simple roots with respect to the appropriate Weyl chamber.
The corresponding positive roots are
\[
\Phi^+ := \{\; 2\varepsilon,\;\   \ \varepsilon_i + \varepsilon_j + \varepsilon_k +\varepsilon \;(i>j>k), \   \ \varepsilon_i -\varepsilon_j \;(i>j)\; \}.
\]
Denote by $\Lambda$ the integer lattice generated by $\bfpi$. Since the curves $a_0;
a_2,\ldots,a_6$ form a basis of the integer homology group of $\Sigma$, we can identify $\Lambda$
with $\sfh_1(\Sigma, \zz)$. As we have shown, every homomorphism $\varphi_\mu: \sfh_1(\Sigma,
\zz_2) \to \zz_2$ is obtained by $\zz_2$-reduction from either a homomorphism
$\lambda \in \Lambda \mapsto B(\alpha_t, \lambda) \in \zz$ with some $\alpha_t \in \Phi$,
or a homomorphism $\lambda \in \Lambda \mapsto B(\gamma, \lambda) \in \zz$ with some weight $\gamma$
of the $\frg'$-module $V^+$. Moreover, all homomorphisms $\varphi_\mu$ of the same type
yield the same value of $\hat \varphi_\mu(\idel^2(\rme_6))$.

As a representative of the first $\sfw(\rme_6)$-orbit we take the root $2\varepsilon$.
Then $\hat \varphi_\mu(\idel^2(\rme_6))$ is the number of roots
$\varepsilon_i + \varepsilon_j + \varepsilon_k +\varepsilon$ with $i>j>k$.
Thus $\hat \varphi_\mu(\idel^2(\rme_6)) = \binom{6}{3} =20$ in this case.
This gives $\hat \varphi_\mu(L^2) = 32 -2\times20 =-8$.

As a representative of the second $\sfw(\rme_6)$-orbit we take the weight $\varepsilon_1
+\varepsilon_2$. To count the roots in question we use the following observations: First,
since $\varepsilon$ is orthogonal to $\varepsilon_1$ and $\varepsilon_2$, we can replace $\varepsilon_i + \varepsilon_j + \varepsilon_k +\varepsilon$
by $\varepsilon_i + \varepsilon_j + \varepsilon_k$. Second, since $\varepsilon_1+\cdots + \varepsilon_6=0$ and we are interested only
in the parity of $B(\varepsilon_1 +\varepsilon_2, \alpha_t)$, we can replace each $\varepsilon_i + \varepsilon_j + \varepsilon_k$ by
$\varepsilon_1+\cdots + \varepsilon_6 - \varepsilon_i + \varepsilon_j + \varepsilon_k$. This yields twice each $\varepsilon_i + \varepsilon_j + \varepsilon_k$ with
$i>j>k\geq2$. The explicit combinatorics is:
\begin{itemize}
\item $\varepsilon_2 -\varepsilon_1$ is orthogonal to $\varepsilon_1 +\varepsilon_2$;
\item $B(\varepsilon_1 +\varepsilon_2, \varepsilon_i -\varepsilon_1) = - 1$ for $i\geq3$, and there are $4$ such roots;
\item $B(\varepsilon_1 +\varepsilon_2, \varepsilon_i -\varepsilon_2) = - 1$ for $i\geq3$, and there are $4$ such roots;
\item $B(\varepsilon_1 +\varepsilon_2, \varepsilon_i + \varepsilon_j + \varepsilon_k ) = - 1$ for $i>j>k\geq3$,
and there are $2\times4=8$ such roots;
\item $B(\varepsilon_1 +\varepsilon_2, \varepsilon_2 + \varepsilon_i + \varepsilon_j) = 0$ for $i>j\geq3$.
\end{itemize}
Thus $\hat \varphi_\mu(\idel^2(\rme_6)) = 16$ and $\hat \varphi_\mu(L^2) = 32 -2\times16 =0$.

\medskip\noindent
\sliii Every separating curve $\gamma \subset \Sigma\bs \{ z_0 \}$
divides $\Sigma$ into two pieces $\Sigma'$ and $\Sigma''$ one of
which, say $\Sigma''$, contains $z_0$. Let $p\geq1$ be the genus of the
other piece $\Sigma'$. Then the whole configuration is conjugated to
the one in which $\Sigma'$ contains the curves $a_1,\ldots, a_{2p}$
and $\Sigma''$ the curves $a_{2p+2},\ldots, a_{2g}$ (see the curve
$\gamma_3$ on \reffig{fig-Delta-g}). In particular, in the case
$p=g$ the curve $\gamma$ surrounds the base point $z_0$ and
corresponds to the boundary curve $\dd$ on \reffig{fig-Delta-g}. By
the chain relation (this time \slsf{separating}),
$T_{\gamma_p}T_{\gamma_{p-1}} = (T_{\beta_p} T_{\alpha_p}
T_{\beta'_p})^4$, where the curve $\beta'_p$ is like the one shown on
\reffig{fig-chain} and where we set $T_{\gamma_0} =1 \in \map_{g,1}$
in the special case $p=1$. Thus $T_\gamma$ is conjugated to the
alternating product of the elements $(T_{\beta_p} T_{\alpha_p}
T_{\beta'_p})^4$. The incidence relations in the configuration
$\beta_p,\alpha_p, \beta'_p$ correspond to the Coxeter system
$\rma_3$, whereas the product $(T_{\beta_p} T_{\alpha_p}
T_{\beta'_p})^4$ is the squared Garside element $\idel^2(\rma_3)$.
Thus $(T_{\beta_p} T_{\alpha_p} T_{\beta'_p})^4$ lies in
$\sfp(\cals_g)$ and equals in $\zz\langle\calh_g\rangle$ to the sum
$\sum_{1\leq i\leq j\leq 3} A_{v_{ij}}$ in which $v_{11} = [\beta_p]
\in \sfh_1(\Sigma, \zz_2)$, $v_{22} = [\alpha_p]$, $v_{33} =
[\beta_p']= [\beta_p]$, and $v_{ij}= \sum_{k=i}^j v_{kk}$. An
explicit calculation shows that $\sum_{1\leq i\leq j\leq 3}
A_{v_{ij}} = 2A_{[\beta_p]} + 2A_{[\alpha_p]} + 2A_{[\alpha_p +
\beta_p]}$. It follows that $\hat \varphi_\mu$ takes value $0$ or
$4$ on every such product $(T_{\beta_p} T_{\alpha_p}
T_{\beta'_p})^4$. This implies assertion \sliiip.

\medskip\noindent
\sliv The subgroup $\dd_\pi\big( \pi_1(\Sigma, z_0) \big) \subset \map_{g,1}$ is normally
generated by the basic hyperelliptic relation element $\big[ T_{\beta_g},
\idel^2(\rma_{2g}) \big]$. Its geometric realisation is $\dd_\pi(\beta''_g)$ where
the element $\beta''_g\in \pi_1(\Sigma, z_0)$ is represented by an embedded curve,
still denoted by $\beta''_g$, which is isotopic to $\beta_g$ and passes through
the base point $z_0$. On \reffig{fig-Delta-g} this curve corresponds to the arc $\beta''_g$,
and the correct picture on the \emph{closed} surface is obtained by contracting
the boundary circle $\partial$ to the base point $z_0$. In particular, after this
contraction the curves $\beta_g$ and $\beta'_g$ will cut a regular neighbourhood of the curve
$\beta''_g$. Note that the element $\idel^2(\rma_{2g})$ is represented by a
certain hyperelliptic involution of $\Sigma$ which maps $\beta_g$ to $\beta'_g$.

In view of this algebraic description, it is sufficient to find the
possible values of $\hat \varphi_\mu \big( \big[ T_{\beta_g}, \idel^2(\rma_{2g}) \big] \big)$.
First, we observe that in the group $\zz\langle\calh_g\rangle$ we have the
equality $\idel^2(\rma_{2g}) = \sum_{1\leq i\leq j\leq 2g} A_{v_{ij}}$ where $v_{ij} :=
\sum_{k=i}^j [a_k] \in \sfh_1(\Sigma, \zz_2)$. Then $ T_{\beta_g}(v_{ij}) = v_{ij}$ for
$i\leq j<2g$, so that in the group $\zz\langle\calh_g\rangle$ we obtain
\[\textstyle
\big[ T_{\beta_g}, \idel^2(\rma_{2g}) \big]=
\sum_{i=1}^{2g} \Big( A_{v_{i,2g}+ [\beta_g]} - A_{v_{i,2g}} \Big).
\]
It is easy to see that if $\varphi_\mu[\beta_g]=0$, then
$\hat\varphi_\mu\big(A_{v_{i,2g}+ [\beta_g]} - A_{v_{i,2g}} \Big) =0$ for all $i=1,\ldots,2g$.
So it remains to consider the case $\varphi_\mu[\beta_g]=1$.

Here we observe that $[a_1],\ldots,[a_{2g}]$ form a basis of $\sfh_1(\Sigma, \zz_2)$ in
which $[\beta_g] = [a_1] + [a_3] +\cdots+ [a_{2g-1}]$. Thus the homomorphism $\varphi_\mu:
\sfh_1(\Sigma, \zz_2) \to \zz_2$ is completely determined by its values $\varphi_\mu[a_1] ,\ldots,
\varphi_\mu[a_{2g}]$. In other words, we can obtain all homomorphisms $\varphi_\mu$
by varying these values. Clearly, the group of such transformations of $\varphi_\mu$
preserving the value $\varphi_\mu[\beta_g]=1\in \zz_2$ is generated by the following two types
of \slsf{simple transformations}: either we ``switch'' a single value $\varphi_\mu[a_{2k}]$
on an even curve $a_{2k}$, or the values $\varphi_\mu[a_{2k-1}]$ and $\varphi_\mu[a_{2k+1}]$
on two consecutive odd curves $a_{2k-1}$ and $a_{2k+1}$. Since $\varphi_\mu[\beta_g]=1$, every
value $\hat\varphi_\mu\big( A_{v_{i,2g}+ [\beta_g]} - A_{v_{i,2g}} \Big)$ is either $+1$
or $-1$. The crucial observation is that for any simple transformation we get an
\emph{even} number of sign changes for the values $\hat\varphi_\mu\big( A_{v_{i,2g}+ [\beta_g]} - A_{v_{i,2g}} \Big)$.
Namely, $2k$ signs change when we ``switch'' the value $\varphi_\mu[a_{2k}]$,
and $2$ signs change when we ``switch'' the values $\varphi_\mu[a_{2k-1}]$ and $\varphi_\mu[a_{2k+1}]$.
It follows that the value $\hat\varphi_\mu\big( \big[ T_{\beta_g}, \idel^2(\rma_{2g}) \big] \mod4$
remains unchanged, so that it depends only on $\varphi_\mu([\beta_g])$. An explicit calculation
in the case when $\varphi_\mu$ vanishes on all $[\alpha_1], \ldots, [\alpha_g];\, [\beta_1],\ldots,
[\beta_{g-1}]$ shows that $\hat \varphi_\mu \big( \big[ T_{\beta_g}, \idel^2(\rma_{2g}) \big] \big) = 2g-2$,
as claimed.
\qed

\smallskip%
\newsubsection[map-g-fac]{Factorisation problems in mapping class groups.}%
In this paragraph we consider several factorisation problems in the group
$\map_{g,1}$ and in its subgroup $\map_{g,1,\mu}$ stabilising a given non-zero
homology class $\mu \in \sfh_1(\Sigma, \zz_2)$. The topological meaning of these
problems is the existence or non-triviality of certain special homology classes
in Lefschetz fibrations.

Recall that $\wtmap_g$ denotes the group in the extension
\[
1 \to \sfh_1(\Sigma, \zz_2) \to \wtmap_g \to \map_g\to1,
\]
so that $\wtmap_g$ is the quotient of $\map_{g,1}$ by the image with respect
to $\dd_\pi$ of the kernel of the homomorphism $\pi_1(\Sigma, z_0)\to \sfh_1(\Sigma, \zz_2)$.
For our purposes, only the projections of $F \in \map_{g,1}$ to $\wtmap_g$ will be
relevant. Let $\wtmap_{g,\mu}$ be the stabiliser of $\mu$ in $\wtmap_g$.

\newprop{prop-fi-ti} There exists a homomorphism $\ti\varphi_\mu : \wtmap_{g,\mu} \to
\zz_4$ with the following properties:
\begin{itemize}
\item[\sli] $\ti\varphi_\mu$ vanishes on the subgroup $\wtmap'_{g,\mu} \subset \wtmap_{g,\mu}$
 generated by Dehn twists $T_\delta$ with $\delta \cap \mu \equiv0\mod2$ and commutators
 $[F,F']$ with $F,F'\in\wtmap_{g,\mu}$.
\item[\slii]  The restriction of $\ti\varphi_\mu$ to the kernel $\ker \big( \wtmap_{g,\mu} \to
 \Sp(2g, \zz_2) \big)$ is the $\mod4$ reduction of the homomorphism $\hat\varphi_\mu$.
\item[\sliii]  In particular, $\ti\varphi_\mu( \dd_\pi (\gamma)) \equiv (2g-2) \varphi_\mu(\gamma) \mod4$
for any $\gamma \in \pi_1(\Sigma)$.
\end{itemize}
\end{prop}

\proof Combine \refthm{thm-subgr} for the group $G:= \sfw(\cals)_\mu$
with \lemma{lem-chain-eval} and \refcorol{cor-lant-ext}.
\qed

\medskip

\newthm{thm-mu0} Let $g= 2g'\geq2$ be even and let\/ $\id = \prod_i [F_{2i-1}, F_{2i}] \circ \prod_j
T_{\delta_j}$ be a factorisation in $\map_g$ in which $T_{\delta_j}$ are Dehn twists along
embedded circles $\delta_1,\ldots,\delta_n \subset \Sigma\bs \{z_0\}$ and $[F_{2i-1}, F_{2i}]$
denotes the commutator of $F_{2i-1}, F_{2i} \in\map_g$. Then there exist lifts $\wt F_{2i-1},
\wt F_{2i} \in \map_{g,1}$ and curves $\delta'_j$ on $\Sigma\bs \{z_0\}$ such that each $\delta'_j$
is isotopic to $\delta_i$ on $\Sigma$ and $\prod_i [\wt F_{2i-1}, \wt F_{2i}] \cdot \prod_j
T_{\delta_j'} = \dd_\pi \gamma \in \map_{g,1}$ for some $\gamma \in \pi_1(\Sigma, z_0)$ with \emph{trivial}
$\zz_2$-homology class $[\gamma]=0 \in \sfh_1(\Sigma, \zz_2)$.
\end{thm}

\proof Fix some lifts $\wt F_i \in \map_{g,1}$ of $F_i \in\map_g$ and lift
$T_{\delta_j}$ to $\map_{g,1}$ in the natural way. For any $v \in \sfh_1(\Sigma, \zz_2)$
we denote by $\vartheta_v$ the corresponding element in $\wtmap_g$.  Then $\prod_i [\wt
F_{2i-1}, \wt F_{2i}] \cdot \prod_j T_{\delta_j} = \vartheta_{w_0} \in \wtmap_g$ for some $w_0 \in
\sfh_1(\Sigma, \zz_2)$.

Denote by $W$ the set of all $w \in \sfh_1(\Sigma, \zz_2)$ for which $\vartheta_w$ can
be represented as the product $\prod_i [\wt F_{2i-1}, \wt F_{2i}] \cdot \prod_j T_{\delta_j'}$
for some possible lifts $\wt F_i$ and $T_{\delta_j'}$. Then $W = w_0 +V$ for an
appropriate $\zz_2$-subspace $V \subset \sfh_1(\Sigma, \zz_2)$, and the assertion of
the theorem is equivalent to the claim that $w_0 \in V$. Assuming the contrary, $V$ must be a
proper subspace of $\sfh_1(\Sigma, \zz_2)$.

Let $\mu \in \sfh_1(\Sigma, \zz_2)$ be an element such that the functional $\varphi_\mu: v \in
\sfh_1(\Sigma, \zz_2) \mapsto \mu \cap v \in \zz_2$ vanishes on $V$ but is non-zero on $w_0$.
We claim that $\mu \in \sfh_1(\Sigma, \zz_2)$ is invariant under the action of all
factors $F_i, T_{\delta_j}$.  It follows from \lemma{lem-h2-circ1}\ that $V$
consists of (the Poincar{\'e} duals of) coboundaries $d^{(2)}(\bfla^\lor)$ with
$\bfla^\lor$ having the form \eqqref(h2lef). The equivalent dual condition is
that, considering the factorisation $\prod_i [F_{2i-1}, F_{2i}] \cdot \prod_j
T_{\delta_j}$ as a single relation word $R$, the boundary $\dd_2(R \otimes_2 \mu)$
cancels all such $\bfla^\lor$. The calculation done in the proof of
\propo{prop-K-van} shows that
\[\textstyle
\dd_2(R \otimes_2 \mu) = \sum_j T_{\delta_j} \otimes_1 w^T_j \mu +
\sum_i \big( F_{2i-1} \otimes_1 (F_{2i} - \id)w^F_i \mu + F_{2i} \otimes_1 (\id - F_{2i-1})w^F_i \mu\big)
\]
in which $w^T_j, w^F_i$ are certain final subwords of $R$. More precisely,
$w^T_j$ consists of all letters of $R$ after $T_{\delta_j}$, and $w^F_i$ of all
letters of $R$ after $F_{2i}$. Now, consider the conditions $w^T_j \mu \cap \delta_j \equiv0
\mod2$ starting from the last letter $T_{\delta_l}$ and going backwards.
Recursively, we obtain the desired congruence $\mu \cap \delta_j \equiv0 \mod2$.
A similar argument for the conditions $(F_{2i} - \id)w^F_i \mu \cap \lambda_i \equiv 0$ and $(F_{2i-1} -
\id)w^F_i \mu \cap \nu_i \equiv 0$ gives the desired result for $F_i$.

Now, $\ti \varphi_\mu(\vartheta_{w_0}) \equiv 0\mod4$ by \propo{prop-fi-ti}~{\sl i)}.
On the other hand, $\ti \varphi_\mu(\vartheta_{w_0}) \equiv (2g-2)\varphi_\mu(w_0) \equiv2(\mu\cap w_0) \equiv2$
by \propo{prop-fi-ti}~{\sl iii)}. This contradiction shows that the desired lifts exist.
\qed

\smallskip\noindent{\bf Remark.}
From this theorem and~\lemma{lem-h2-circ1} we obtain \refthm{thm10} and thence
the non-triviality of the $\zz_2$-homology class of the fibre for any
topological Lefschetz fibration with even fibre genus stated in the remark
following \propo{prop-mu-0}.

\medskip
From the proof of \refthm{thm-mu0} we can conclude the following:

\newcorol{cor-mu1} Let $g \geq3$ be odd and let\/ $\id = \prod_i [F_{2i-1}, F_{2i}] \circ
\prod_j T_{\delta_j}$ be a factorisation in $\map_g$ as in \refthm{thm-mu0}. Assume
that for every lift $\prod_i [\wt F_{2i-1}, \wt F_{2i}] \cdot \prod_j T_{\delta_j'} = \dd_\pi \gamma \in
\map_{g,1}$ with $\gamma \in \pi_1(\Sigma, z_0)$ the $\zz_2$-homology class $[\gamma]_{\zz_2}$ is
non-trivial. Then there exists $\mu \in \sfh_1(\Sigma, \zz_2)$ stabilised by
all $F_i$ and all $T_{\delta_j}$ and such that $\mu \cap [\gamma] \equiv1 \mod2$ for every $\gamma \in
\pi_1(\Sigma, z_0)$ as above.
\end{corol}


\newthm{thm-mu-prim} \sli Let $g\geq1$ be odd. Then the element $(T_{\alpha_1}T_{\beta_1}
)^3$ can not be represented in $\map_g$ as a product of Dehn twists and
commutators $\prod_i T_{\delta_i}^{\epsilon_i}\cdot \prod_j[F_{2j-1}, F_{2j}]$ such that $\delta_i \cap \beta_1 \equiv0
\mod2$ and such that every $F_j$ fixes the homology class $[\beta_1]_{\zz_2}
\in \sfh_1(\Sigma, \zz_2)$.

\slii Let $g\geq2$ be even, let $\prod_i T_{\delta_i}^{\epsilon_i}\cdot \prod_j [F_{2j-1}, F_{2j}]$ be a
factorisation in $\map_g$ of the element $(T_{\alpha_1}T_{\beta_1} )^3$ satisfying
the properties above, and let $\wh F_j$ be some lifts of the involved factors
to $\map_{g,1}$. Then $(T_{\alpha_1} T_{\beta_1} )^{-3} \prod_i T_{\delta_i}^{\epsilon_i}\cdot \prod_j [\wh
F_{2j-1}, \wh F_{2j}] = \dd_\pi(\gamma)$ for some $\gamma \in \pi_1(\Sigma, z_0)$ satisfying $[\gamma] \cap
\beta_1 \equiv1 \mod2$.
\end{thm}

\proof \sli We consider the special case $g=1$ first, abbreviating the notation
$\alpha_1, \beta_1$ to $\alpha,\beta$. Here $\cals_1 = \{ s_1, s_2\} = \rma_2$, and hence
$\sfw(\cals_1) = \sym_3$ and $\br(\cals_1) = \br_3$.  Moreover, we can
identify $\map_{1,[1]}$ with $\br_3$ setting $T_\alpha, T_\beta$ to be the standard
generators of $\br_3$. Besides we have $\calt(\cals_1) = \calh_1 = \{[\alpha], [\beta],
[\alpha+\beta] \}$, all three are $\zz_2$-homology classes. Let us apply \refthm{thm-subgr}
to the subgroup $G := \zz_2\langle T_\beta \rangle \subset \sym_3$. It has two orbits
in $\calt(\cals_1)$, namely, $\{[\beta]\}$ and $\{[\alpha], [\alpha+\beta] \}$.
Let $A_{G\cdot [\beta]}$, $A_{G\cdot [\alpha]}$ be the corresponding basis of $\sfp_\ab(
\cals_1)_G$.  As in the proof above, the projection from $\sfp_\ab(\cals_1) =
\zz\langle A_{[\alpha]}, A_{[\beta]}, A_{[\alpha+\beta]} \rangle$ onto the component
$\zz\langle A_{G\cdot [\alpha]}\rangle$ is given by the homomorphism
$\hat \varphi_{[\beta]}: \sfp_\ab(\cals_1) \to \zz$. Then \refthm{thm-subgr} shows that
$\hat \varphi_{[\beta]}(F) =0$ for every factorisation $F$ as in the theorem, whereas
$\hat \varphi_{[\beta]}( (T_\alpha T_\beta)^3) = \hat \varphi_{[\beta]}(A_{[\alpha]}+ A_{[\beta]} + A_{[\alpha+\beta]})=2$.
This contradiction excludes the existence of such a factorisation.

\medskip
Now consider the general case of an odd $g\geq3$ and assume that a factorisation
$\prod_i T_{\delta_i}^{\epsilon_i}\cdot \prod_j [F_{2j-1}, F_{2j}]$ as in the hypotheses exists.
Lifting it to $\br(\cals_g)$ we obtain an element $\wh F$ lying in the kernel
$\ker\big( \br(\cals_g) \to \Sp(2g, \zz_2) \big)$. Repeating the argument from the proof
of \refthm{thm-mu0} we obtain that $\hat \varphi_{[\beta_1]}(\wh F) \equiv0 \mod 4$.
Note that such a lift $\wh F$ differs from $(T_{\alpha_1}T_{\beta_1} )^3$
(considered now as an element from $\br(\cals_g)$) by a product of chain, lantern,
and hyperelliptic relation elements. Since $g$ is odd, $\hat \varphi_{[\beta_1]}(x) \equiv0 \mod 4$
for every relation element $x$ (including the hyperelliptic one), so the value
$\hat \varphi_{[\beta_1]} \big((T_{\alpha_1}T_{\beta_1} )^3 \big)$ must be a multiple of $4$.
On the other hand, the calculation done in the case $g=1$ shows that
$\hat \varphi_{[\beta_1]} \big((T_{\alpha_1}T_{\beta_1} )^3 \big)=2$, a contradiction.

\medskip\noindent
\slii As above, we apply the argument used in the proof of \refthm{thm-mu0}.
This yields $\hat \varphi_{[\beta_1]} \big( (T_{\alpha_1}T_{\beta_1} )^3 \cdot \dd_\pi(\gamma)
\big) \equiv0 \mod 4$. Since $\hat \varphi_{[\beta_1]} \big((T_{\alpha_1}T_{\beta_1} )^3 \big)=2$,
the assertion follows from~\propo{prop-fi-ti}, \sliiip.
\qed

\newsection[proof]{Proof of the main result}

In this section we give the proof of the \slsf{Main Theorem}.  We maintain the
notation of \refsubsection{hmlg-TLF}. In particular, $\mbfm$ denotes the
meridian circle of $K$, considered as a curve on $\Sigma$, and $\mu =
[\mbfm]_{\zz_2}$ its homology class on $\Sigma$.

In the case when $\mbfm$ separates $\Sigma$, the theorem follows from the remark
immediately after \propo{prop-mu-0}. So we assume henceforth that $\mbfm$ is non-separating on
$\Sigma$. Furthermore, we may assume that the hypothesis of \propo{prop-K-van} is
fulfilled and the projection of $[K]$ to $\sfh_1( Y^\circ, \scrh_1(X_y, \zz_2))$
vanishes, since otherwise there is nothing to prove.

Let us make the following observation about the monodromy $F_\Gamma$ along $\Gamma$.
Realise the meridian $\mbfm$ as the curve $\beta_1$ in some geometric basis $\alpha_1,
\beta_1, \ldots, \alpha_g,\beta_g$ of $\Sigma$. Then we can isotopically deform $F_\Gamma$ so that $\beta_1 =
\mbfm$ is $F_\Gamma$-stable. The induced map $F_\Gamma: \beta_1 \to \beta_1$ inverts the
orientation. Observe that the map $(T_{\alpha_1} T_{\beta_1})^3$ has the same property
in the homotopy, namely, it maps the free homotopy class of $\beta_1$ onto itself
inverting the orientation. It follows that the homotopy class $F'_\Gamma := F_\Gamma
(T_{\alpha_1} T_{\beta_1})^{-3}$ admits a representative (still denoted by $F'_\Gamma$)
which fixes $\beta_1$ pointwise. In particular, $F'_\Gamma$ is a diffeomorphism of $\Sigma
\bs \beta_1$. The classification of diffeomorphisms of surfaces implies that
$F'_\Gamma$ can be represented as a product of Dehn twists $T_\delta$ along curves
disjoint from $\beta_1$. Fix such a decomposition of $F'_\Gamma$ and lift it to the
group $\map_{g, 1}$. The obtained elements of $\map_{g, 1}$ are denoted by
$\wh F'_\Gamma$ and $\wh F_\Gamma := \wh F'_\Gamma \cdot (T_{\alpha_1} T_{\beta_1})^3$.

\newsubsection[fine-mndr]{Fine structure of the monodromy.} In this paragraph
we refine the result of \propo{prop-K-van}.

\newprop{prop-fine} Under the hypotheses of \propo{prop-K-van} \slip, assume
in addition that $[K]$ vanishes in $\sfh_2(X, \zz_2)$. Then the
decomposition $\mu=\mu_+ +\mu_-\in \sfh_1(\Sigma, \zz_2)$ constructed in \propo{prop-K-van}
satisfies $\mu \cap \mu_+ \equiv \mu \cap \mu_- \equiv 1 \mod2$.

Moreover, the class $\mu$ is {\bf\em not} invariant either with respect to the
monodromy action of $\pi_1(Y^\circ_+)$ or with respect to the action of $\pi_1(Y^\circ_-)$.
\end{prop}

\proof We must exclude the following two possibilities:
\begin{itemize}
\item One of the classes $\mu_\pm$ vanishes, say, $\mu_-=0$ and so $\mu_+ =\mu$.
\item Both  $\mu_+, \mu_-$ are non-trivial, but  $\mu \cap \mu_+ \equiv \mu \cap \mu_- \equiv 0 \mod2$.
\end{itemize}

\smallskip%
Consider first the case $\mu_-=0$. Then the monodromy $F_\Gamma$ admits a
representation as a product $F_+ := \prod_i [F_{2i-1}, F_{2i}] \cdot \prod_j T_{\delta_j}$ such
that $\mu$ is invariant with respect to all $F_i$ and all $T_{\delta_j}$.  This gives
the equality $(T_{\alpha_1} T_{\beta_1})^3 = F_+ \cdot (F'_\Gamma)\inv$ with the same property
for the right hand side.  \refthm{thm-mu-prim} excludes such a possibility for
odd $g$ and ensures the following in the case of even $g\geq2$. For any lift $\wh
F_+$ of $ F_+$ to the group $\map_{g,1}$ we obtain the relation $(T_{\alpha_1}
T_{\beta_1})^3 = \dd_\pi(\gamma) \cdot\wh F_+ \cdot (\wh F'_\Gamma)\inv $ with some $\gamma \in \pi_1(\Sigma)$ whose
homology class $[\gamma] \in \sfh_1(\Sigma, \zz_2)$ satisfies $[\gamma] \cap \mu \equiv1 \mod2$.  Now
choose such a lift $\wh F_+$ of $F_+$ as a part of a lift of the whole
monodromy using \refthm{thm-mu0}. The lift $\wh F_\Gamma \in \map_{g,1}$ gives us a
section $\sigma_K : \Gamma \to X$ of the fibration $\pr: X \to Y$ over $\Gamma$ which is disjoint
from $K$. On the other hand, the lift $\wh F_+$ defines a $\zz_2$-section
$\sigma|_\Gamma$ of $\pr: X \to Y$ over $\Gamma$ which is the restriction of some
$\zz_2$-section $\sigma$. The condition $[\gamma] \cap \mu \equiv1 \mod2$ means that the
$\zz_2$-intersection index $\sigma \cap [K]$ is non-zero, hence $[K]$ is non-trivial.
This rules out the case $\mu_-=0$.

\smallskip%
This shows also that $\mu$ can not be invariant with respect to the
monodromy action of $\pi_1(Y^\circ_+)$ or $\pi_1(Y^\circ_-)$.

\medskip%
Let us assume that both $\mu_+$ and $\mu_-$ are non-trivial but $\mu_+ \cap \mu \equiv \mu_- \cap \mu
\equiv 0$. Then $g\geq2$.  Set $l := \hat \varphi_{\mu_+}(\wh F_\Gamma)$. In the factorisation of
$\wh F'_\Gamma = \wh F_\Gamma \cdot (T_{\alpha_1} T_{\beta_1})^{-3}$ considered above there are no Dehn
twist $T_{\delta_j}$ with $\delta_j \cap \mu \equiv1$.  So using the definition of $\hat \varphi_{\mu},
\hat \varphi_{\mu_\pm}$ we conclude that $\hat \varphi_{\mu}(\wh F_\Gamma)=2$ and $\hat \varphi_{\mu_-}(\wh
F_\Gamma) =l+2$.

Let $F_+$ and $F_-$ be the factorisations of $F_\Gamma$ of the form $F_\pm := \prod_i
[F_{2i-1}, F_{2i}] \cdot \prod_j T_{\delta_j}$ arising from the monodromy homomorphisms
$\calf_\pm: \pi_1(Y^\circ_\pm) \to \map_g$. Lifting the monodromy homomorphisms $\calf_\pm$
to $\map_{g,1}$ we obtain lifts $\wh F_\pm \in \map_{g,1}$ of $F_\Gamma$, which differ
from $\wh F_\Gamma$ by some elements lying in the image of the homomorphism $\dd_\pi:
\pi_1(\Sigma) \to \map_{g,1}$.

Assume that $g$ is odd. Then $\hat \varphi_{\mu}, \hat \varphi_{\mu_\pm}$ vanish $\mod4$ on the
image of $\dd_\pi$. Applying \refthm{thm-subgr} we conclude that $\hat
\varphi_{\mu_+}(\wh F_\Gamma) \equiv \hat \varphi_{\mu_-}(\wh F_\Gamma) \equiv0\mod4$. This contradicts the
relation $\hat \varphi_{\mu_-}(\wh F_\Gamma) - \hat \varphi_{\mu_+}(\wh F_\Gamma) =2$.

Assume that $g$ is even. Define $\gamma_\pm \in \pi_1(\Sigma)$ from the relations $\wh F_\pm =
\wh F_\Gamma \cdot \dd_\pi(\gamma_\pm)$. By \refthm{thm-mu0} there exists a lift $\wt\calf:
\pi_1(Y^\circ) \to \wtmap_g$ of the monodromy homomorphism $\calf: \pi_1(Y^\circ) \to \map_g$,
see \refsection{map-g-fac}. By \lemma{lem-h2-circ1}, every such lift
$\wt\calf$ corresponds to a $\zz_2$-section $\sigma$. We can suppose that the lifts
$\wh F_\pm \in \map_{g,1}$ of $F_\Gamma$ are compatible with such a section $\sigma$. This
means that the projections of $\wh F_\pm$ to $\wtmap_g$ are equal. This shows that
$[\gamma_+] = [\gamma_-] \in \sfh_1(\Sigma, \zz_2)$. Comparing the values of $\hat
\varphi_{\mu}, \hat \varphi_{\mu_\pm}$ and using $\hat \varphi_{\mu_\pm}(\wh F_\pm) =0$, we obtain $\hat
\varphi_{\mu_+}\big( \dd_\pi(\gamma_+) \big) = -l$ and $\hat \varphi_{\mu_-}\big( \dd_\pi(\gamma_-) \big) =
-l-2$. Now, using \propo{prop-fi-ti} \sliiip, we obtain $[\gamma_+] \cap \mu \equiv1
\mod2$. This means that $[K] \cap [\sigma] \not\equiv0$, which contradicts the
hypotheses of the proposition.
%
%
\qed

\smallskip

\newsubsection[proof-top]{Proof of \rm\slsf{Main Theorem}} Let us change
the notation slightly and denote by $X_0$ the original symplectic manifold
($\cp^2$ or ruled one) and by $X$ the blow-up of $X_0$ constructed in
\lemma{lem-slf}.

Let us assume that the claim of the \slsf{Main Theorem} is false and $K\subset X$ is
homologically trivial. First, we sum up the properties of $\pr:X \to Y$ obtained so
far. By the construction of \lemma{lem-slf}, $Y$ is the sphere $S^2$ and $\Gamma$
separates $Y$ into discs $Y_\pm$. By \propo{prop-mu-0}, the meridian circle
$\mbfm$ does not separates $\Sigma$, and by \propo{prop-fine} the monodromy of
$\pr: X \to Y$ has the following structure: there exist $\mu_+,\mu_- \in
\sfh_1(\Sigma, \zz_2)$ such that $\mu \cap \mu_+ \equiv \mu \cap \mu_- \equiv 1 \mod2$ and such that the
monodromy $\calf_\pm: \pi_1(Y_\pm^\circ) \to \map_g$ of each piece preserves the class
$\mu_\pm$ but does not preserve $\mu$. Finally, the monodromy $F_\Gamma$ along $\Gamma$ can
be realised by a map, still denoted by $F_\Gamma \in \diff_+(\Sigma)$, which maps $\mbfm$
onto itself reversing the orientation on $\mbfm$.

Fix local ``polar'' coordinates $(r, \varphi)$ in a neighbourhood $U_\mbfm$ of $\mbfm$
on $\Sigma$ so that $\mbfm$ is defined by the equation $r=1$, $\varphi$ is an angle
coordinate on $\mbfm$, and the map $F_\Gamma$ is given by $(r,\varphi) \mapsto (r\inv, -\varphi)$. It
follows that there exists a geometric realisation $T_\mbfm \in \diff_+(\Sigma)$
of the Dehn twist along $\mbfm$ which commutes with $F_\Gamma$. Using
this map $T_\mbfm$ and a local trivialisation of the bundle $\pr: X_\Gamma \to\Gamma$ used
in the definition of the monodromy map $F_\Gamma$, we obtain a map $\Psi: X_\Gamma \to X_\Gamma$
which preserves the fibres and acts on each fibre $X_y$ ($y\in \Gamma$) by the
map $T_\mbfm$.

Construct a new manifold $X'$ by gluing together the pieces $X_{Y_+}$ and
$X_{Y_-}$ along $X_\Gamma = \dd X_{Y_+} = \dd X_{Y_-}$ via the map $\Psi: X_\Gamma \to
X_\Gamma$. Obviously, $X'$ admits a Lefschetz fibration $\pr':X' \to Y$ over the
same base $Y= S^2$. Moreover, the monodromy of $\pr':X' \to Y$ remains unchanged
in $Y_+$ and gets conjugated by $T_\mbfm$ in $Y_-$. Since $[\mbfm]_{\zz_2} =\mu$
and $\mu_+ =\mu + \mu_-$, the whole monodromy of $\pr':X' \to Y$ preserves the class
$\mu_+$.

\newlemma{lem-h1-xprim} $\rank\;\sfh_1(X',\zz_2) = \rank\; \sfh_1(X,\zz_2) +1 $.
\end{lem}

\proof Let $V_+$ and $V_-$ be the $\zz_2$-subspaces of $\sfh_1(\Sigma,\zz_2)$
generated by the vanishing classes $[\delta_i]$ of the projections $\pr:X_{Y_+} \to Y_+$
and $\pr:X_{Y_-} \to Y_-$, respectively. Then by \lemma{lem-h1},
\[
\sfh_1(X,\zz_2) = \sfh_1(\Sigma,\zz_2) \big/ ( V_+ + V_-)
\text{\ \ and\ \ }
\sfh_1(X',\zz_2)= \sfh_1(\Sigma,\zz_2) \big/ ( V_+ + T_\mu(V_-))
\]
Set $W_0 := \zz_2\langle \mu_+, \mu_- \rangle^\perp \subset \sfh_1(\Sigma,\zz_2)$ so that $\sfh_1(\Sigma,\zz_2)=
\zz_2\langle \mu_+, \mu_- \rangle \oplus W_0$ is an orthogonal decomposition. It follows from the
second assertion of \propo{prop-fine} and the relation $\mu= \mu_+ + \mu_-$ that
$V_+ \big/ ( V_+ \bigcap W_0)$ has rank $1$ and is generated by the coset class of
$\mu_+$. Similarly, $V_- \big/ ( V_- \bigcap W_0) = \zz_2\langle \mu_- \rangle$. Since $T_\mu(\mu_-) =
\mu_+$, we obtain
{\small\[
( V_+ + V_-) / \big( ( V_+ + V_-) \bigcap W_0 \big) = \zz_2\langle \mu_+, \mu_- \rangle
\text{\ \ and\ \ }
( V_+ + T_\mu(V_-)) / \big( ( V_+ + T_\mu(V_-)) \bigcap W_0 \big) = \zz_2\langle \mu_+\rangle.
\]}%
The lemma follows.
\qed

\remark I am grateful to Stefan Nemirovski for the following observation.
The twisting construction of the manifold $X'$ above is the \slsf{Luttinger
surgery} of $X$ along the Klein bottle $K$, see \refsubsection{0.1a} and
\cite{N-2}. In its turn, the description of the construction shows that the
Luttinger surgery along the Klein bottle $K$ is compatible with the projection $\pr$
of the Morse--Lefschetz fibration $\pr: (X, K) \to (Y,\Gamma)$ provided that the
restricted projection $\pr: K \to \Gamma \cong S^1$ is an $S^1$-bundle without critical
points. Note that this compatibility property holds also if $K$ is replaced by the torus $T^2$
(see \cite{ADK}).

\newlemma{lem-xprim-top} Both $\rank\;\sfh_1(X,\zz_2)$ and $\rank\;
\sfh_1(X',\zz_2)$  must be even.
\end{lem}

\proof The classification of ruled symplectic manifolds \cite{MD-Sa} implies
that a finite sequence of blow-ups and blow-downs transforms $X$ into a
product $S^2 \times Y$ where $Y$ is a closed oriented surface. Since blow-ups do
not change $\pi_1(X)$, we conclude the first part of the lemma.

The second part is obtained in the same way once we show that $X'$ is also a ruled
symplectic manifold. First, we observe that the gluing map $\Psi: X_\Gamma \to X_\Gamma$
can be extended to a symplectomorphism of some neighbourhood of $X_\Gamma$. It follows
that $X'$ carries a symplectic form $\omega'$ that coincides with the original form $\omega$
on $X$ on the pieces $X_{Y_\pm}$.

At this point we use the specific structure of $\omega$ and the monodromy of $\pr:
X \to Y$. Recall that $X$ was constructed as a symplectic blow-up of the
original $X_0$. It follows that there exist symplectic sections $E_1, \ldots, E_N \subset X$
of the projection $\pr: X\to S^2$ such that for the class
$[D] := [\Sigma] + \sum_i [E_i]\in \sfh_2(X,\zz)$ one has $c_1(X) \cdot [D]>0$.
The sections $E_i$ are simply the exceptional spheres resulting from the blow-up construction.
Since $E_i$ are disjoint from $K$, they survive in $X'$, and we obtain symplectic sections
$E'_1,\ldots ,E'_N$ in $X'$. Moreover, for the class $[D'] := [\Sigma] + \sum_i [E_i']\in
\sfh_2(X,\zz)$ we again have $c_1(X') \cdot [D']>0$. In this situation the characterisation
theorem of McDuff and Salamon, see \slsf{Corollary 1.5} in \cite{MD-Sa}, says that $X'$
must indeed be a ruled symplectic manifold.
\qed

\medskip%
The obtained contradiction shows that under the hypothesis of the \slsf{Main Theorem} the monodromy of $\pr: X\to Y$
can not have the structure described in \propo{prop-fine}. This implies the assertion of the \slsf{Main Theorem}.
\qed


\end{document}